\documentclass[a4paper,11pt]{article}

\usepackage{theorem}

\usepackage[all]{xy}

\usepackage{stmaryrd}
\usepackage{graphicx}
\usepackage{amsfonts}
\usepackage[leqno]{amsmath}
\usepackage{amssymb,mathrsfs}
\usepackage{color}
\usepackage{times}
\usepackage{enumerate,vmargin}
\usepackage{tocvsec2}
\usepackage{verbatim}
\RequirePackage[numbers,sort&compress]{natbib}

\setpapersize{A4}
\setmarginsrb{2.5cm}{2cm}{2.5cm}{2cm}{.4cm}{4mm}{0cm}{7.5mm}

\newtheorem{lem}{Lemma}[section]
\newtheorem{thm}[lem]{Theorem}
\newtheorem{prop}[lem]{Proposition}

{\theorembodyfont{\rmfamily}\newtheorem{defi}{Definition}[section]}
{\theorembodyfont{\rmfamily}}
{\theorembodyfont{\rmfamily}}
{\theorembodyfont{\rmfamily}\newtheorem{rem}{Remark}[section]}
{\theorembodyfont{\rmfamily}}

\newcommand{\noi}{\noindent}

\newcommand{\lgeo}{\llbracket}
\newcommand{\rgeo}{\rrbracket}
\newcommand{\bbm}{{\bf m}}
\newcommand{\un}{{\bf 1}}

\newcommand{\bnu}{\boldsymbol \nu}
\newcommand{\bQ}{\mathbf{Q}}
\newcommand{\bC}{\mathbf{C}}
\newcommand{\cE}{\mathcal E}
\newcommand{\ccM}{\mathscr M}
\newcommand{\bbC}{\mathbb C}

\newcommand{\ep}{\varepsilon}

\newcommand{\cT}{\mathcal T}

\newcommand{\cF}{\mathcal F}
\newcommand{\cG}{\mathcal G}

\newcommand{\cL}{\mathcal L}
\newcommand{\cM}{\mathcal M}

\newcommand{\cI}{\mathcal I}
\newcommand{\cJ}{\mathcal J}
\newcommand{\cH}{\mathcal H}

\newcommand{\cO}{\mathcal O}

\newcommand{\bbE}{\mathbb E}
\newcommand{\bbP}{\mathbb P}

\newcommand{\bbR}{\mathbb R}

\newcommand{\bbN}{\mathbb N}

\newcommand{\bE}{\mathbf E}
\newcommand{\bN}{\mathbf N}

\newcommand{\bP}{\mathbf P}

\newcommand{\bH}{\mathbf H}

\newcommand{\ccB}{\mathscr B}
\newcommand{\bD}{\mathbf D}
\newcommand{\Gam}{\Gamma}

\renewcommand{\[}{\left[}
\renewcommand{\]}{\right]}
\renewcommand{\(}{\left(}
\renewcommand{\)}{\right)}

\newcommand{\lam}{\lambda}
\newcommand{\gam}{\gamma}

\newcommand{\indi}{\boldsymbol{1}}

\newcommand{\igam}{\frac{1}{\gam}}

\newcommand{\nr}{\bN_{\! {\rm nr}}}

\def\cq{$\hfill \square$}
\def\cqfd{$\hfill \blacksquare$}
\def\ino{\! \in \! }
\def\epp{\varepsilon}
\def\llcr{\lambda_{{\rm cr}}}

\begin{document}

\title{{\huge \textsc{Decomposition of L\'evy trees along their diameter.}}}
\author{Thomas \textsc{Duquesne}
\thanks{\textbf{Institution}: PRES Sorbonne Universit\'es, UPMC Universit\'e Paris 06,  LPMA (UMR 7599). \textbf{Postal address}: LPMA, Bo\^ite courrier 188, 4 place Jussieu, 75252 Paris Cedex 05, FRANCE. \textbf{Email}: thomas.duquesne@upmc.fr. \textbf{Grant}: ANR-14-CE25-0014 (ANR GRAAL)} \and Minmin \textsc{Wang}
\thanks{\textbf{Institution}: PRES Sorbonne Universit\'es, UPMC Universit\'e Paris 06,  LPMA (UMR 7599). \textbf{Postal address}: LPMA, Bo\^ite courrier 188, 4 place Jussieu, 75252 Paris Cedex 05, FRANCE. \textbf{Email}: wangminmin03@gmail.com. \textbf{Grant}: ANR-14-CE25-0014 (ANR GRAAL)}}
\date{\today}

\maketitle
\begin{abstract}
We study the diameter of L\'evy trees that are random compact metric spaces obtained as the scaling limits of Galton-Watson trees. L\'evy trees have been introduced by 
Le Gall and Le Jan (1998) and they generalise Aldous' Continuum Random Tree (1991) that corresponds to the Brownian case. We first characterize the law of the diameter of L\'evy trees and we prove that it is realized by a unique pair of points. We prove that the law of L\'evy trees 
conditioned to have a fixed diameter $r \! \in \! (0, \infty)$ is obtained by glueing at their respective roots 
two independent size-biased L\'evy trees conditioned to have height $r/2$ and then by uniformly re-rooting the resulting tree; we also describe by a Poisson point measure the law of the subtrees that are grafted on the diameter. As an application of this decomposition of L\'evy trees according to their diameter, we characterize the joint law of the height and the diameter of stable L\'evy trees conditioned by their total mass; we also provide asymptotic expansions of the law of the height and of the diameter of such normalised stable trees, which generalises the identity due to Szekeres (1983) in the Brownian case.

\smallskip


\smallskip

\end{abstract}



\section{Introduction and main results}

L\'evy tree are random compact metric spaces that are the scaling limits of Galton-Watson trees. The Brownian tree, also called the continuum random
tree, is a particular instance of 
L\'evy trees; it is the limit of the rescaled uniformly distributed rooted labelled tree with $n$ vertices. The Brownian tree has been introduced by Aldous in \cite{aldcrt1} and further studied in Aldous 
\cite{aldcrt2, aldcrt3}. L\'evy trees have been introduced by Le Gall \& Le Jan \cite{legalllejan} via a coding function called the height process that is a local time functional of a 
spectrally positive L\'evy process. L\'evy trees (and especially stable trees) have been studied in D.~\& Le Gall \cite{Duquesne02, Duquesne05} (geometric and fractal properties, connection with superprocesses), see D.~\& Winkel \cite{Duquesne07} and Marchal \cite{marchal08} for alternative constructions, 
see also Miermont \cite{Miermont_stable1, Miermont05}, Haas \& Miermont \cite{HaMi04}, Goldschmidt \& Haas \cite{GolHaa} for applications to stable fragmentations, and Abraham \& Delmas \cite{AbDel07, AbDel08}, Abraham, Delmas \& Voisin \cite{AbDeVo} for general fragmentations and pruning processes on L\'evy trees. 

   In this article, we study the diameter of L\'evy trees. As observed by Aldous (see \cite{aldcrt2}, Section 3.4), in the Browian case the law of the diameter has been found by Szekeres \cite{Sz83} by taking the limit of the generating function of the diameter of uniformly distributed rooted labelled tree with $n$ vertices. 
Then, the question was raised by Aldous that whether we can derive the law of the diameter 
directly from the normalised Brownian excursion that codes the Brownian tree (see also Pitman \cite{pitmanstflour}, Exercise 9.4.1). This question is now answered in W.~\cite{Wang14}. 

  In this article we compute the law of the diameter for general L\'evy trees (see Theorem \ref{diamlaw}). We also prove that the diameter of L\'evy trees 
is realized by a unique pair of points. In Theorem \ref{thm: decomp}, we describe the coding function (the height process) of the L\'evy trees tree rerooted at the midpoint of their diameter that plays the role of an intrinsic root. The proof of Theorem \ref{thm: decomp} relies on the invariance of L\'evy trees by uniform rerooting, as proved by D.~\& Le Gall in \cite{Duquesne09}, and on the decomposition of L\'evy trees according to their height, as proved by  
Abraham \& Delmas in \cite{AbDe09} (this decomposition generalizes Williams decomposition of the Brownian excursion). 
Roughly speaking, Theorem \ref{thm: decomp} asserts that a L\'evy tree that is conditioned to have diameter $r$ and that is rooted at its midpoint is obtained 
by glueing at their root two size-biased independent L\'evy trees conditioned to have height $r/2$; 
Theorem \ref{thm: decomp} also 
explains the distribution of the subtrees that are grafted on the diameter. As an application of this theorem, we characterize the joint law of the height and the diameter of stable trees conditioned by their total mass (see Proposition \ref{prop: lp}) and we provide asymptotic expansions for the distribution of the law of the height (Theorem \ref{thm: ht_asy}) and for the law of the diameter (Theorem \ref{thm: dm_asy}). These two asymptotic expansions generalize the identities due to Szekeres in the Brownian case which involves theta functions (these identities are recalled in (\ref{htBronor}) and (\ref{diaBronor})). Theorem \ref{htdm0} also provides precise asymptotics of the tail at zero of the law of the height and that of the diameter of normalised stable trees. 
Before stating precisely our main results we need to recall definitions and to set notations.

\paragraph{Real trees.} 
Real trees are metric spaces extending the definition of graph-trees: 
let $(T, d)$ be a metric space; it is a {\it real tree} iff the following holds true.

\smallskip

\noi
(a)~~For any $\sigma_1, \sigma_2 \! \in\!  T$, there is a unique isometry 
$f\! : [0,d(\sigma_1,\sigma_2)] \! \rightarrow \! T$ such
that $f(0)\!=\! \sigma_1$ and $f(d(\sigma_1,\sigma_2))\! =\! \sigma_2$. Then, we shall use the following notation:   
$\lgeo \sigma_1,\sigma_2\rgeo \! :=\! f([0,d (\sigma_1,\sigma_2)])$. 

\smallskip

\noi
(b)~~For any continuous injective function 
$q: [0, 1] \! \rightarrow \! T$,  $q([0,1]) \! = \! \lgeo q(0), q(1)\rgeo$.

\smallskip

\noi
When a point $\rho \! \in\!  T$ is distinguished, $(T, d, \rho)$ is said to be a \textit{rooted} real tree, $\rho$ being the \textit{root} of $T$. 
Among connected metric spaces, real trees are characterized by the so-called {\it four-point condition} that is expressed as follows: let $(T, d)$ be a connected metric space; then $(T, d)$ is a real tree iff for any $\sigma_1,  \sigma_2,  \sigma_3,  \sigma_4  \in T$, we have 
\begin{equation}
\label{fourpoint}
d(\sigma_1, \sigma_2) + d(\sigma_3, \sigma_4) \leq \big(d(\sigma_1, \sigma_3) + d(\sigma_2, \sigma_4)\big) \vee  \big( d(\sigma_1, \sigma_4) + d(\sigma_2, \sigma_3)  \big) . 
\end{equation}
We refer to Evans \cite{evans05} or to Dress, Moulton and Terhalle \cite{DMT96} for a detailed account on this property. 
Let us briefly mention that the set of (pointed) isometry classes of compact rooted real trees can be equipped with the (pointed) Gromov-Hausdorff distance that makes it a Polish space: see Evans, Pitman \& Winter \cite{EPW}, Theorem 2, for more details on this intrinsic point of view on trees that we shall not use here.

\paragraph{The coding of real tree.}
Let us briefly recall how real trees can be obtained thanks to continuous functions. To that end we denote by $\bC(\bbR_+, \bbR_+)$ the space of $\bbR_+$-valued continuous function equipped with the topology of the uniform convergence on every compact subsets  of $\bbR_+$. We shall denote by $H\! = \! (H_t)_{t\geq 0}$ the canonical process on $\bC(\bbR_+, \bbR_+)$. 
We first assume that $H$ has a compact support, that $H_0 \! =\! 0$ and that $H$ is distinct from the null function: we call such a function a \textit{coding function} and 
we then set $\zeta_H \! = \! \sup \{ t\! >\! 0: H_t \! >\! 0 \}$ that is called the \textit{lifetime} of the coding function $H$. Note that $\zeta_H \! \in \! (0, \infty)$. 
Then, for every $s,t\! \in \! [0, \zeta_H]$, we set
\begin{equation}
\label{pseudometric}
b_H(s,t)=\inf_{r\in[s\wedge t,s\vee t]}H_r \quad {\rm and} \quad d_H(s,t)=H_s+H_t-2b_H(s,t).
\end{equation}
It is easy to check that $d_H$ satisfies the four-point condition: namely, for all $s_1, s_2, s_3, s_4 \! \in \!  [0, \zeta_H]$, 
$d_H(s_1,s_2) + d_H(s_3, s_4) \! \leq \! \big(d_H(s_1, s_3) + d_H(s_2, s_4)\big) \! \vee \!  \big( d_H(s_1, s_4) + d_H(s_2, s_3)  \big) $. By taking $s_3\! = \! s_4$, we see that $d_H$ is a pseudometric on $[0, \zeta_H]$.  We then introduce the equivalence relation
$s\! \sim_H \! t$ iff $d_H(s,t) \! =\! 0$ and we set 
\begin{equation}
\label{codef}
\cT_H= [0, \zeta_H] / \! \sim_H \; . 
\end{equation}
Standard arguments show that $d_H$ induces a true metric on the quotient set $\cT_H$ that we keep denoting by $d_H$. We denote by $p_H: [0, \zeta_H] \rightarrow \cT_H$ the \textit{canonical projection}. Since $H$ is continuous, so is $p_H$ and $(\cT_H, d_H)$ is therefore a compact connected metric space that satisfies the four-point condition: it is a compact real tree. We next set $\rho_H= p_H(0)= p_H(\zeta_H)$ that is chosen as the \textit{root} of $\cT_H$. 

We next define the \textit{total height} and the \textit{diameter} of $\cT_H$ that 
are expressed in terms of $d_H$ as follows: 
\begin{equation}
\label{heigdiam}
\Gamma (H) \! :=\! \! \sup_{\sigma \in \cT_H} \! d_H(\rho_H, \sigma) \! = \! \! \!\!  \sup_{t\in [0, \zeta_H ]}\!\!  \! H_t  
\quad \textrm{and} \quad D(H)\! :=\! \! \!\! \sup_{\sigma, \sigma^\prime\in \cT_H}\!\! d_H(\sigma , \sigma^\prime) \! = \! 
\!\!\!\! \!\!  \sup_{0\leq s<t \leq \zeta_H} \!\! \!\! \big( H_s+H_t \! -\! 2\!\!  \inf_{r\in [s, t]} \!\! H_r \big).
\end{equation}
For any $\sigma \! \in \! \cT_H$, we denote by ${\rm n} (\sigma)$ the number of connected components of 
the open set $\cT_H \backslash \{ \sigma\}$. Note that ${\rm n} (\sigma)$ is possibly infinite. We call this number the \textit{degree} of $\sigma$. We say that $\sigma$ is a \textit{branching point} if 
${\rm n} (\sigma) \! \geq\!  3$; we say that $\sigma$ is a \textit{leaf} if ${\rm n} (\sigma) \! =\!  1$ and we say that $\sigma$ is \textit{simple} if ${\rm n} (\sigma) \! = \! 2$. We shall use the following notation for the set of branching points and the set of leaves of $\cT_H$: 
\begin{equation}
\label{brleafset}
\mathtt{Br} (\cT_H):= \big\{ \sigma \! \in \! \cT_H: {\rm n} (\sigma) \! \geq \! 3 \big\} \quad \textrm{and} \quad \mathtt{Lf} (\cT_H):= \big\{ \sigma \! \in \! \cT_H: {\rm n} (\sigma) \! =\! 1 \big\} \; . 
\end{equation}
In addition to the metric $d_H$ and to the root $\rho_H$, the coding function yields two additional useful features: first, the \textit{mass measure} 
$\bbm_H$ that is the pushforward measure of the Lebesgue measure on $[0, \zeta_H]$ induced by $p_H$ on $\cT_H$; namely, for any Borel measurable function $f: \cT_H\rightarrow \bbR_+$, 
\begin{equation}
\label{masmea}
 \int_{\cT_H} \!\!\!\!  f(\sigma ) \, \bbm_H (d\sigma) = \! \int_0^{\zeta_H} \!\!  \!\!\! f(p_H(t)) \, dt \; .    
\end{equation}
This measure plays an important role in the study of L\'evy trees (that are defined below): in a certain sense, the mass measure is the most spread out measure on $\cT_H$. The coding $H$ also induces a \textit{linear order} $\leq_H$ on $\cT_H$ that is inherited from that of $[0, \zeta_H]$: namely for any $\sigma_1 , \sigma_2 \! \in \! \cT_H$, 
\begin{equation}
\label{linorder} 
\sigma_1 \leq_H  \sigma_2  \quad \Longleftrightarrow \quad \inf \{ t\! \in \! [0, \zeta_H]: p_H (t) \! = \! \sigma_1 \} \leq \inf \{ t\! \in \! [0, \zeta_H]: p_H (t) \! = \! \sigma_2 \} \; . 
\end{equation}
Roughly speaking, the coding function $H$ is completely characterized by $(\cT_H, d_H, \rho_H, \bbm_H, \leq_H)$: see D.~\cite{Duquesne06} for more detail about the coding of real trees by functions.

\paragraph{Re-rooting trees.} Several statements of our article
involve a re-rooting procedure at the level of the coding functions that is recalled here from D.~\& Le Gall \cite{Duquesne05}, Lemma 2.2  (see also D.~\& Le Gall \cite{Duquesne09}). Let $H$ be a coding function as defined above and recall that $\zeta_H \! \in \! (0, \infty)$. For any $t\! \in \! \bbR_+$, denote by $\overline{t}$ the unique element of $[0, \zeta_H)$ such that $t\! -\! \overline{t}$ is an integer multiple of $\zeta_H$. Then for all 
$t_0 \! \in \! \bbR_+$, 
we set 
\begin{equation}
\label{rerootH}
\forall t\! \in \! [0, \zeta_H] , \quad  H^{[t_0]}_t = d_H \big(  \overline{t_0}, \overline{t+t_0}\,  \big) \quad \textrm{and} \quad 
\forall t \geq \zeta_H, \quad   H^{[t_0]}_t =0 \;. 
\end{equation}     
Then observe that $\zeta_H= \zeta_{H^{[t_0]}}$ and that  
\begin{equation}
\label{isoreroot}
\forall t, t' \in [0, \zeta_H], \quad d_{H^{[t_0]}} (t, t') = d_H \big(\, \overline{t+t_0}, \overline{t'+t_0} \, \big) \; . 
\end{equation}
Then, Lemma 2.2 \cite{Duquesne05} asserts that there exists a unique isometry $\phi: \cT_{H^{[t_0]}} \rightarrow \cT_{H}$ such that $\phi (p_{H^{[t_0]}}  (t))= p_H \big( \overline{t+t_0}\big)$ for all $t \! \in \! [0, \zeta_H]$. This allows to \textit{identify canonically $\cT_{H^{[t_0]}}$ with the tree 
$\cT_H$ re-rooted at $p_H (t_0)$}: 
\begin{equation}
\label{identif}
\big( \cT_{H^{[t_0]}}, d_{H^{[t_0]}} , \rho_{H^{[t_0]}} \big) \equiv \big( \cT_H, d_H , p_H(t_0) \big) \; .
\end{equation} 
Note that up to this identification, $\bbm_{H^{[t_0]}}$ is the same as $\bbm_H$. Roughly speaking, the linear order $\leq_{H^{[t_0]}}$ is obtained from $\leq_H$ by a cyclic shift after $p_H(t_0)$.

\paragraph{Spinal decomposition.} The law of the L\'evy tree conditioned by its diameter that is discussed below is described as a Poisson decomposition of the trees grafted along the diameter. To explain such a decomposition in terms of the coding function of the tree, we introduce the following definitions and notations.

\begin{figure}
   \begin{minipage}[c]{.46\linewidth}
      \includegraphics[scale=1.2]{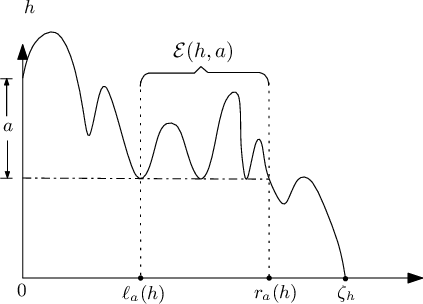}
   \end{minipage} \hfill
   \begin{minipage}[c]{.42\linewidth}
      \includegraphics[scale=0.9]{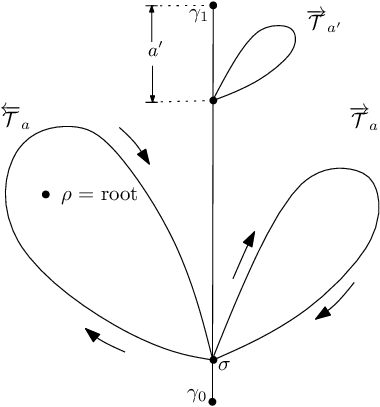}
   \end{minipage}
\caption{\label{spinal}{\small \textit{ \hspace{-3mm}the figure on the left hand side illustrates the definition of $\cE(h,a)$; the figure on the right hand side represents the spinal decomposition of $H$ at times $t_0$ and $t_1$ in terms of the tree $\cT$ coded by $H$.}} }
\end{figure}

  Let $h\! \in \! \bC(\bbR_+, \bbR_+) $ have compact support. Note that $h(0) \! >\! 0$ possibly. 
We first define the excursions of $h$ above its infimum as follows. For any $a\! \in \! [0, h(0)]$, we first set 
$$ \ell_a(h):= \inf \! \big\{ t\! \in \! \bbR_+ : h(t) \! = \! h(0) \! -\! a   \big\} \quad \textrm{and} \quad  r_a(h):= \zeta_h \wedge \inf \! \big\{ t\! \in \! (0, \infty) : h(0)\! -\! a > h(t)  \big\}, $$
with the convention that $\inf \emptyset \! = \! \infty$, so that $r_{h(0)}(h)\! =\!  \zeta_h$. We then set 
$$ \forall s \! \in \! \bbR_+ , \quad \cE_s (h, a) := h \big(  (\ell_a (h) +s)\! \wedge \! r_a(h) \big)-h(0) +a \; .$$
See Figure \ref{spinal}. Note that $\cE (h, a)$ is a nonnegative continuous function with compact support such that $\cE_0(h, a)\! = \! 0$. Moreover, if $\ell_a (h)\!  =\! r_a(h)$, then $\cE(h, a)\! = \! \mathbf{0}$, the \textit{null function}. 

  Let $H$ be a coding function as defined above. Let $t \! \in \! \bbR_+$, we next set 
$$\forall s \! \in \! \bbR_+, \quad  H^-_s= H_{(t-s )_+}   \quad \textrm{and} \quad  H^{+}_s= H_{t+ s}\; .$$  
Note that $H^-_0\! =\! H^+_0\! =\! H_t$. To simplify notation we also set 
$$ \forall  a \in [0, H_t], \quad \overleftarrow{H}^a:= \cE (H^- \! , a) \quad \textrm{and} \quad  \overrightarrow{H}^a:= \cE (H^+\!  , a) $$
and $\cJ_{0,t}:= \big\{ a\! \in \! [0, H_t] :  \textrm{either $\ell_a (H^-)  \! < \! r_a (H^-) $ or $\ell_a (H^+) \! < \! r_a (H^+) $} \big\} $, that is countable. We then define the following point measure on $[0, H_t] \! \times \! \bC(\bbR_+, \bbR_+)^2$:  
\begin{equation}
\label{spindec0t}
 \cM_{0, t} (H)= \sum_{a \in \cJ_{0, t}} \delta_{(a, \overleftarrow{H}^a, \overrightarrow{H}^a)} \; , 
\end{equation} 
with the convention that $\cM_{0, t} (H)\! = \! 0$ if $\cJ_{0, t}\! = \! \emptyset $. In Lemma \ref{mesdec2}, we see that if $\bbm_H $ is diffuse and supported by the set of leaves of $\cT_H$, then there is a measurable way to recover $(t, H) $ from $\cM_{0, t} (H)$.

For all $t_1\! \geq t_0 \! \geq \! 0$, we also set 
\begin{equation}
\label{spint0t1}  
\cM_{t_0, t_1} (H):= \cM_{0, t_1-t_0} \big( H^{[t_0]}\big)=: \sum_{a \in \cJ_{t_0, t_1}} \delta_{(a, \overleftarrow{H}^a, \overrightarrow{H}^a)} \; .
\end{equation} 
This point measure on $[0, d_H(t_0, t_1)] \! \times \! \bC(\bbR_+, \bbR_+)^2$ is 
the \textit{spinal decomposition of $H$ between $t_0$ and $t_1$}.

\begin{rem}
\label{etcxfxs}
Let us interpret this decomposition in terms of the tree $\cT_H$ (more precisely in terms of the tree $\cT_{H^{[t_0]}}$, see Figure \ref{spinal}). Let us set $\gamma_0\! =\! p_H(t_0)$ and $\gamma_1\! = \! p_H(t_1)$; to simplify our explanation, we assume that $\gamma_0$ and $\gamma_1$ are leaves. 
Recall that $\lgeo \gamma_0, \gamma_1 \rgeo$ is the geodesic path joining $\gamma_0$ to $\gamma_1$; then $\cJ_{t_0, t_1} \! =\!  \{ d( \sigma , \gamma_1) ; \sigma \! \in \! \mathtt{Br} (\cT_H) \cap \, \lgeo \gamma_0, \gamma_1\rgeo\}$. For any positive $a \! \in \! \cJ_{t_0, t_1}$, there exists $\sigma \! \in \!  \mathtt{Br} (\cT_H) \cap \, \lgeo \gamma_0, \gamma_1\rgeo$ such that the following holds true.  

\smallskip

\noi
$\bullet$ $\overleftarrow{\cT}_{\!\! a} \! : =\! \{ \sigma\} \cup \big\{\sigma^\prime \! \in \! \cT_H: \gamma_0 \! <_H \! \sigma^\prime \! <_H \! \gamma_1\; \textrm{and} \; \lgeo \gamma_0, \sigma \rgeo \! = \!  \lgeo \gamma_0, \sigma^\prime \rgeo  \cap   \lgeo \gamma_0, \gamma_1 \rgeo   \big\}$ is the tree grafted at $\sigma$ on the left hand side of $\lgeo \gamma_0, \gamma_1\rgeo$ and the tree $(\overleftarrow{\cT}_{\!\! a} , d, \sigma)$ is coded by $\overleftarrow{H}^a$.

\smallskip

\noi
$\bullet$ $\overrightarrow{\cT}_{\!\! a} \! := \! \{ \sigma\} \cup \big\{ \sigma^\prime \! \in \! \cT_H:  \textrm{either $\sigma^\prime \! <_H \! \gamma_0$ or $\gamma_1  \! <_H \! \sigma^\prime$ and} \; \lgeo \gamma_0, \sigma \rgeo \! =\!   \lgeo \gamma_0, \sigma^\prime \rgeo   \cap   \lgeo \gamma_0, \gamma_1 \rgeo    \big\}$ is the tree grafted at $\sigma$ on the right hand side of $\lgeo \gamma_0, \gamma_1 \rgeo$ and the tree $(\overrightarrow{\cT}_{\!\! a} , d, \sigma)$ is coded by $\overrightarrow{H}^a$. \cq 
\end{rem}

\paragraph{Height process and L\'evy trees.} The Brownian tree (also called Continuum Random Tree) has been introduced by Aldous \cite{aldcrt1, aldcrt2, aldcrt3}; this model has been extended by Le Gall \& Le Jan: in \cite{legalllejan}, they define the \textit{height process} (further studied by D.~\& Le Gall \cite{Duquesne02}) that is the coding function of L\'evy trees. L\'evy trees appear as scaling limits of Galton-Watson trees and they are the genealogical structure of 
continuous state branching processes. Let us briefly recall here the definition of the height process and that of L\'evy trees. 

 The law of the height process is characterized by a function $\Psi \! :\!  \bbR_+\!  \rightarrow \! \bbR_+$ called \textit{branching mechanism}; we shall restrict our attention to the critical and subcritical cases, namely 
 when the branching mechanism $\Psi$ is of the following L\'evy-Khintchine form: 
\begin{equation}\label{eq: def-Psi}
\forall \lambda \in \bbR_+, \quad \Psi(\lam)=\alpha \lam+\beta \lam^2+\int_{(0, \infty)}\!\!\!\! \!\!\! \big( e^{-\lam r}\! -\! 1+\lam r \big) \, \pi(dr) \; , 
\end{equation}   
where $\alpha, \beta\! \in \! \bbR_+$ and where $\pi$ is the L\'evy measure on $(0, \infty)$ that satisfies $\int_{(0, \infty)}(r\! \wedge \! r^2)\, \pi(dr)\! <\! \infty$.  The height process is derived from a spectrally positive L\'evy process whose Laplace exponent is $\Psi$. It shall be convenient to work with the canonical process $X\! =\! (X_t)_{t\geq 0}$ on the space of c\`adl\`ag functions $\bD(\bbR_+, \bbR)$ equipped with the Skorohod topology. Let us denote by $\bbP$ the law of a spectrally positive L\'evy 
process starting from $0$ and whose Laplace exponent is $\Psi$. Namely, 
$$ \forall t, \lam \in \bbR_+, \quad \bbE\[\exp\(-\lam X_t\)\]=\exp\big(t \Psi(\lam)\big) \; . $$ 
Note that the form (\ref{eq: def-Psi}) ensures that $X$ under $\bbP$ does not drift to $\infty$: see for instance 
Bertoin \cite{Bebook}, Chapter VII for more details. 
Under the following assumption: 
\begin{equation}\label{eq: hyp}
\int_1^\infty \frac{d\lambda}{\Psi(\lambda)} <\infty, 
\end{equation}
Le Gall \& Le Jan \cite{legalllejan} (see also D.~\& Le Gall \cite{Duquesne02}) have proved that there exists a continuous process $H\! =\!  (H_t)_{t\geq 0}$ such that for all $t\! \in \! \bbR_+$, 
the following limit holds in $\bbP$-probability: 
\begin{equation}\label{eq: defH}
H_t=\lim_{\ep\to 0}\frac{1}{\ep}\int_0^t ds \,\indi_{\{I^s_t<X_s<I^s_t+\ep\}},
\end{equation}
where $I^s_t:=\inf_{s<r<t}X_r$. The process $H$ is called the \textit{$\Psi$-height process}. In the Brownian case, namely when $\Psi (\lambda) \! =\! \lam^2$, easy arguments show that $H$ is distributed as a reflected Brownian motion. Le Gall \& Le Jan \cite{legalllejan} have proved a Ray-Knight theorem for $H$, which shows that the height process $H$ 
codes the genealogy of continuous state branching processes (see also D.~\& Le Gall \cite{Duquesne02}, Theorem 1.4.1). Moreover, the $\Psi$-height process $H$ appears as the scaling limit of the discrete height process and the contour function of Galton-Watson discrete trees: see D.~\& Le Gall \cite{Duquesne02}, Chaper 2, for more details. 

 For all $x \! \in \! (0, \infty)$, we set $T_x \! =\!  \inf \{ t\! \in \! \bbR_+ : X_t= \! - x \}$, that is $\bbP$-a.s.~finite since $X$ under $\bbP$ does not drift to $\infty$. We next introduce the following law 
$\bP^x$ on $\bC(\bbR_+ , \bbR_+)$: 
\begin{equation} 
\label{xLevfo}
\textrm{$\bP^x$ is the law of $(H_{t\wedge T_x})_{t\geq 0}$ under $\bbP$.}
\end{equation}
The tree $\cT_H$ under $\bP^x (dH) $ is called the \textit{$\Psi$-L\'evy forest starting from a population of size $x$}. 
Then, the mass measure of $\cT_H$ under $\bP^x(dH)$ satisfies the following important properties: 
\begin{equation}
\label{conttreex}
\textrm{$\bP^x(dH)$-a.s.~$\bbm_H$ is diffuse and  $\bbm_H(\cT_H \backslash \mathtt{Lf}(\cT_H))= 0$,}
\end{equation}
where we recall from (\ref{brleafset}) that $\mathtt{Lf} (\cT_H)$ stands for the set of leaves of the tree $\cT_H$. 
The $\Psi$-L\'evy forest $(\cT_H, d_H , \rho_H, \bbm_H)$ is therefore a \textit{continuum tree} according to the definition of Aldous \cite{aldcrt3}. 

\medskip

Each excursion above $0$ of $H$ under $\bP^x$ corresponds to a tree of the L\'evy forest. Let us make this point precise by introducing a Poisson decomposition of $H$ into excursions above $0$. 
To that end, denote by $I$ the infimum process of $X$: 
$$\forall t \in \bbR_+, \quad  I_t= \inf_{0\leq r\leq t} X_r \; .$$
Observe that (\ref{eq: hyp}) entails that either 
\begin{equation}
\label{infvarhyp}
\beta \! >\! 0 \quad \textrm{or}  \quad \int_{(0, 1)} r\, \pi(dr)\! =\! \infty \; , 
\end{equation} 
which is equivalent for the L\'evy process $X$ to have unbounded variation sample paths; basic results of fluctuation theory (see for instance Bertoin \cite{Bebook}, Sections VI.1) entail that $X\! -\! I$ is a strong Markov process in $[0, \infty)$ and that $0$ is regular for 
$(0, \infty)$ and recurrent with respect to this Markov process. Moreover, $\! -I$ 
is a local time at $0$ for $X\! -\! I$ (see Bertoin \cite{Bebook}, Theorem VII.1). 
We denote by $\bN$ the corresponding excursion 
measure of $X\! -\! I$ above $0$.

  It is not difficult to derive from (\ref{eq: defH}) that $H_t$ only depends on the excursion of $X\! -\! I$ above $0$ which straddles $t$. Moreover, we get 
$\{ t\! \in \! \bbR_+ : H_t \! >\! 0 \}= \{ t \! \in \! \bbR_+ : X_t \! >\! I_t\}$ and if we denote by $(a_i, b_i)$, $i\! \in \! \cI$, the connected components of this set and if we set $H^{i}_s = H_{(a_i +s)\wedge b_i}$, $s\! \in \! \bbR_+$, then the point measure 
\begin{equation}\label{PoisdecH}
 \sum_{i\in \cI} \delta_{(-I_{a_i} ,\,  H^{i})}   
\end{equation} 
is a Poisson point measure on $\bbR_+ \! \times \! \bC(\bbR_+, \bbR_+)$ with intensity $dx \, \bN(dH)$, where, with a slight abuse of notation, $\bN (dH)$ stands for the 'distribution' of $H(X)$ under $\bN (dX)$.  
In the Brownian case, up to scaling, $\bN$ is It\^o positive excursion of Brownian motion and the decomposition (\ref{PoisdecH}) corresponds to the Poisson decomposition of a reflected 
Brownian motion above $0$. 

In what follows, we shall mostly work with the $\Psi$-height process 
$H$ under its excursion $\bN$ that is a sigma-finite measure on $\bC(\bbR_+, \bbR_+)$. 
We simply denote by $\zeta$ the \textit{lifetime} of $H$ under $\bN$ and we easily check that 
\begin{equation}\label{HNcoding}
 \textrm{$\bN$-a.e.} \quad \zeta \! <\!  \infty\, , \quad H_0\! =\! H_{\zeta}\! =\! 0 \quad \textrm{and} \quad H_t \! >\! 0 \;  \Longleftrightarrow \; t \! \in \! (0, \zeta) \; .
\end{equation} 
Also note that $X$ and $H$ under $\bN$ have the same lifetime $\zeta$ and basic results of fluctuation theory (see for instance Bertoin \cite{Bebook}, Chapter VII) also entail the following: 
\begin{equation}
\label{lifetimeexc}
\forall \lambda\in (0, \infty)\, ,  \quad \bN \big[ 1\! -\! e^{-\lambda \zeta}  \big]= \Psi^{-1} (\lambda ), 
\end{equation}
where $\Psi^{-1}$ stands for the inverse function of $\Psi$.

Note that (\ref{HNcoding}) shows that $H$ under $\bN$ is a coding function as defined above. D.~\& Le Gall \cite{Duquesne05} then define the \textit{$\Psi$-L\'evy tree} as the real tree coded by $H$ under $\bN$.

\begin{itemize}
\item[] \textit{Convention.} When there is no risk of confusion, we simply write 
$$ \big( \cT, d, \rho, \bbm, \leq, p, \Gamma, D \big) := \big( \cT_H, d_H, \rho_H, \bbm_H, \leq_H, p_H, \Gamma (H), D(H) \big) $$
when $H$ is considered under $\bN$, $\bP^x$ or under other measures on $\bC(\bbR_+, \bbR_+)$.  \cq 
\end{itemize}
Recall from (\ref{brleafset}) that $\mathtt{Lf} (\cT)$ stands for the set of leaves of $\cT$. 
Then the mass measure has the following properties: 
\begin{equation}
\label{conttree}
\textrm{$\bN$-a.e.~$\bbm$ is diffuse and $\bbm(\cT \backslash \mathtt{Lf}(\cT))= 0$.}
\end{equation}
The $\Psi$-L\'evy tree $(\cT, d , \rho, \bbm)$ is therefore a continuum tree according to the definition of Aldous \cite{aldcrt1}.

\paragraph{Diameter decomposition.} Recall from (\ref{heigdiam}) the definition of the total height $\Gamma$ and that of the diameter $D$. Let us first briefly recall results on the total height. 
One checks that the total height is $\bN$-a.s.~realized at a unique time   
(see D.~\& Le Gall \cite{Duquesne05} and also Abraham \& Delmas \cite{AbDe09}). Namely,   
\begin{equation}
\label{taudef} 
\textrm{$\bN$-a.e.~there exists a unique $\tau \! \in \! [0, \zeta]$ such that $H_{\tau} = \Gamma$}\; .
\end{equation}
Moreover, the distribution of the total height $\Gamma$ under $\bN$ is characterized as follows: 
\begin{equation}
\label{lawGamma}  
\forall t \in (0, \infty), \quad v(t):= \bN(\Gamma >t) \quad \textrm{satisfies} \quad \int_{v(t)}^\infty \! 
\frac{d\lambda}{\Psi (\lambda)}= t \; .
\end{equation}
Note that $v \! :\!  (0, \infty) \! \rightarrow \! (0, \infty)$ is a bijective decreasing $C^\infty$ function and (\ref{lawGamma}) implies that on $(0, \infty)$, $\bN( \Gamma \! \in \! dt) \! =\!  \Psi (v(t)) \, dt $.

Recall from (\ref{xLevfo}) that $\bP^x$ is the law of $(H_{t\wedge T_x})_{t\geq 0}$ under $\bbP$, where $T_{x}\! =\!  \inf \{ t \! \in \! \bbR_+: X_t\! =\!  \! -x \}$.  The Poisson decomposition (\ref{PoisdecH}) implies that $\sup_{t \in [0, T_x]} H_t \! = \! \max \{ \Gamma (H^{i} ) ; i \! \in \! \cI: -I_{a_i} \! \leq \! x\}$ and since $\Gamma$ under $\bN$ has a density, then (\ref{taudef}) and (\ref{lawGamma}) entail that  
\begin{equation}
\label{tauPx} 
\textrm{$\bP^x$-a.s.~there is a unique $\tau \! \in \! [0, \zeta]$ such that $H_{\tau} = \Gamma$} \quad \textrm{and} \quad \textrm{$\bP^x (\Gamma \leq t) \! = \! e^{-xv(t)}$, $t\! \in \! \bbR_+$. }
\end{equation}

In \cite{AbDe09}, Abraham \& Delmas generalize Williams' decomposition of the Brownian excursion to the excursion of the $\Psi$-height process: they first make sense of the conditioned law $\bN (\, \cdot \, | \, \Gamma \! = \! r)$. Namely they prove that $\bN (\, \cdot \, | \, \Gamma \! = \! r)$-a.s.~$\Gamma \! = \! r$, that  $r\mapsto \bN (\, \cdot \, | \, \Gamma \! = \! r)$ is weakly continuous on $\bC(\bbR_+, \bbR_+)$ and that 
\begin{equation}
\label{heightcondef}
 \bN = \int_0^\infty  \!\!\! \bN(\Gamma \! \in \! dr)\,  \bN (\, \cdot \, | \, \Gamma \! = \! r) \; .
 \end{equation}
Moreover they provide a Poisson decomposition along the total height of the process: see Section \ref{pfth1th2sec} where a more precise statement is recalled.

The first two results of our article provide a similar result for the diameter $D$ of the $\Psi$-L\'evy tree under $\bN$. 
Recall that $p: [0 , \zeta ] \!  \rightarrow \!  \cT$ stands for the canonical projection.  
\begin{thm}
\label{diamlaw} Let $\Psi$ be a branching mechanism of the form (\ref{eq: def-Psi}) that satisfies (\ref{eq: hyp}). 
Let $\cT$ be the $\Psi$-L\'evy tree that is coded by the $\Psi$-height process $H$ under the excursion measure $\bN$ as defined above. Then, the following holds true $\bN$-a.e.

\begin{itemize}
\item[(i)] There exists a unique pair $\tau_0, \tau_1 \! \in \! [0, \zeta]$ such that $\tau_0\! <\! \tau_1$ and $D= d(\tau_0, \tau_1)$. Moreover, either $H_{\tau_0}\! =\!  \Gamma$ or $H_{\tau_1} \! =\!  \Gamma$. Namely, either $\tau_0\! = \! \tau$ or  $\tau_1\! = \! \tau$, where $\tau$ is the unique time realizing the total height as defined by (\ref{taudef}).

\item[(ii)] Set $\gamma_0\! = \! p(\tau_0)$ and $\gamma_1\! = \! p(\tau_1)$. Then $\gamma_0$ and $\gamma_1$ are leaves of $\cT$. 
Let $\gamma_{\rm mid}$ be the mid-point of $\lgeo \gamma_0, \gamma_1 \rgeo$: namely, $\gamma_{{\rm mid}}$ is the unique point of  
$\lgeo\gamma_0, \gamma_1 \rgeo$ such that $d( \gamma_0 , \gamma_{\rm mid})= D/2$. 
Then, there are exactly two times $0 \! \leq \tau_{{\rm mid}}^-\! < \! \tau_{{\rm mid}}^+\! \leq \! \zeta$ such that $p( \tau_{{\rm mid}}^-)\! =\!  p(\tau_{{\rm mid}}^+) \! =\!  \gamma_{\rm mid}$,  
and $\gamma_{{\rm mid}}$ is a simple point of $\cT$: namely, it is neither a branching point nor a leaf of $\cT$.

\item[(iii)] For all $r \! \in \! (0, \infty)$, we get 
\begin{equation}
\label{diamrepa}
\bN \big( D \! > 2r  \big)= v(r) -\Psi \big( v(r)\big)^2 \!\! \int_{v(r)}^\infty \! \frac{d\lambda}{\Psi(\lambda)^2}\; .
\end{equation}
This implies that $\bN (D \! \in \! dr) \! = \! \varphi (r) \, \!  dr$ on $(0, \infty)$ where the density 
$\varphi\! : \! (0, \infty) \! \rightarrow \! (0, \infty)$ is given  by 
 \begin{equation}
\label{diamdens}
\forall r \in (0, \infty), \quad \varphi (2r)= \Psi (v(r))-\Psi(v(r))^2 \, \Psi^\prime(v(r)) \!\! \int_{v(r)}^\infty \! \frac{d\lambda}{\Psi(\lambda)^2} \; .
\end{equation}
\end{itemize}
\end{thm}

The second main result of our paper is a Poisson decomposition of the subtrees of $\cT$ grafted on the diameter $\lgeo \gamma_0, \gamma_1  \rgeo$. This result is stated in terms of coding functions and we first need to introduce the following notation: let $H, H^\prime \! \in \! \bC(\bbR_+, \bbR_+)$ be two coding functions as defined above; 
the \textit{concatenation} of $H$ and $H^\prime$ is the coding function denoted by $H\oplus H^\prime$ and given by 
\begin{equation}
\label{concadef}
\forall t\in \bbR_+, \qquad (H\oplus H^\prime)_t = H_t \quad \textrm{if $t\in [0, \zeta_H]$} \quad \textrm{and} \quad (H\oplus H^\prime)_t = H^\prime_{t-\zeta_H} \quad \textrm{if $t \geq \zeta_H$.}
\end{equation}
Moreover, to simplify notation we write the following: 
\begin{equation}
\label{htcondsim}
 \forall r \! \in \! (0, \infty), \quad \bN^{\Gamma}_{r}= \bN (\, \cdot \, | \, \Gamma \! = \! r) \; .
\end{equation}

\begin{thm}
\label{thm: decomp} 
Let $\Psi$ be a branching mechanism of the form (\ref{eq: def-Psi}) that satisfies (\ref{eq: hyp}). 
For all $r \! \in \! (0, \infty)$, we denote by $\bQ_r$ the law on $\bC(\bbR_+, \bbR_+)$ 
of $H \! \oplus \! H^\prime$ under $\bN^{\Gamma}_{r/2}( dH) \bN^{\Gamma}_{r/2}( dH^\prime)$, where 
$\bN^{\Gamma}_{r/2}$ is defined by (\ref{htcondsim}). 
Namely, for all measurable functions $F \! :\!  \bC (\bbR_+, \bbR_+) \! \rightarrow \! \bbR_+$,
\begin{equation}
\label{defQr}
\bQ_r \big[ F(H)\big]= \int \!\!\!\! \int_{\bC(\bbR_+, \bbR_+)^2} \!\!\!\!\!\!\!\! \!\!\!\!\!\!\!\! \!\!\!\! \bN^{\Gamma}_{r/2}( dH) \bN^{\Gamma}_{r/2}( dH^\prime) \; F\big( H \!  \oplus \! H^\prime \big)\;  \; .
\end{equation}
Then $\bQ_r$ satisfies the following properties. 
\begin{itemize}
\item[(i)] $\bQ_r$-a.s.~$D=r$ and there exists a unique pair of points $\tau_0, \tau_1\! \in \! [0, \zeta]$ such that $D= d(\tau_0, \tau_1)$.  
\item[(ii)] For all $r\! \in \! (0, \infty)$, $\bQ_r  [ \, \zeta \, ] \! =\!  2 \bN^{\Gamma}_{r/2} [\, \zeta \,  ] \! \in \! (0, \infty)$. Moreover, the application 
$r\! \mapsto \! \bQ_r$ is weakly continuous and for all measurable functions $F\! :\!  \bC (\bbR_+, \bbR_+) \! \rightarrow\!  \bbR_+$ and $f \! :\!  \bbR_+ \! \rightarrow \! \bbR_+ $, 
\begin{equation}
\label{desindiam}
 \bN \big[ f(D) F(H)\big]= \int_{0}^\infty  \! \frac{\bN (D\! \in \! dr)}{\bQ_r[\, \zeta\, ]}
\, f(r) \,  \bQ_r \Big[ \int_0^\zeta \!\!\! F\big(  H^{[t]} \big) \, dt \Big]    \; , 
\end{equation} 
where $H^{[t]}$ is defined by (\ref{rerootH}).

\item[(iii)] Recall the notation $\tau^-_{{\rm mid}}$ and 
$\tau^+_{{\rm mid}}$ from Theorem \ref{diamlaw} $(ii)$. Then, for all $r \! \in \! (0, \infty)$, 
\begin{equation}
\label{paraph}
 \bN \big[ F\big(H^{[\tau^-_{\rm mid}]} \big) \, \big| \, D\! = \! r\big]= \frac{1}{\bN^{\Gamma}_{r/2} [\, \zeta \, ] }  \int \!\!\!\! \int_{\bC(\bbR_+, \bbR_+)^2} \!\!\!\!\!\!\!\! \!\!\!\!\!\!\!\! \!\!\!\!\! \!   \bN^{\Gamma}_{r/2}( dH) \bN^{\Gamma}_{r/2}( dH^\prime)  \; 
 \,  \zeta_{H^\prime}   F\big( H \! \oplus \!  H^\prime \big) \; , 
\end{equation} 
where $\bN (\, \cdot \, \big| \, D\! = r\! \, )$ makes sense for all $r\! \in \! (0, \infty)$ thanks to (\ref{desindiam}). 
\item[(iv)] Recall from (\ref{xLevfo}) the notation $\bP^y$. To simplify notation, we write for all $y, b \! \in \! (0, \infty)$
\begin{equation}
\label{renota}
 \bN_b = \bN \big(\, \cdot \, \cap\,  \{ \Gamma \leq b\} \big) \quad \textrm{and} \quad \bP^y_b= \bP^y \big( \, \cdot \, \cap \,  \{ \Gamma \leq b\} \big) , 
\end{equation} 
Then, under $\bQ_r$, $\cM_{\tau_0, \tau_1}(da\, d\overleftarrow{H}\, d\overrightarrow{H})$, defined by (\ref{spint0t1}), is a Poisson point measure on $[0, r] \! \times \! \bC(\bbR_+, \bbR_+)^2$ whose intensity is 
\begin{align}
\label{decomQr}
 \beta \un_{[0, r]} (a)  da  & \, \Big( \delta_{\mathbf{0}} (d\overleftarrow{H}) \bN_{a \wedge (r\! - \! a)} ( d\overrightarrow{H}) + \bN_{a \wedge (r\! - \! a)}  ( d\overleftarrow{H}) \delta_{\mathbf{0}} (d\overrightarrow{H})  \Big) \nonumber \\
+ &\; \;  \un_{[0, r]} (a) da \int_{(0, \infty)}\!\!\!\!\!\!\!\! \!\! \pi (dz)\!\!  \int_0^z \!\! \!\! dx \; \bP^x_{a \wedge (r\! - \! a)} \big(d\overleftarrow{H} )\,  \bP^{z-x}_{a \wedge (r\! - \! a)} \big(  d\overrightarrow{H} ) , 
\end{align}
where $\beta$ and $\pi$ are defined in (\ref{eq: def-Psi}) and where $\mathbf{0}$ stands for the null function. 
\end{itemize}
\end{thm}
\begin{rem}
\label{mksense}
As already mentioned, the previous theorem makes sense of $\bN \big( \cdot   \big| D\! = r\! \, \big)$ and for all measurable functions $F \! :\!  \bC(\bbR_+, \bbR_+) \! \rightarrow \! \bbR_+$, we have 
\begin{equation}
\label{Hcondia}
\forall r \! \in \! (0, \infty), \qquad \bN \big[ \,F(H) \, \big| \, D\! = r\! \, \big]  =   \bQ_r \Big[ \int_0^\zeta \!\!\! F\big(  H^{[t]} \big) \, dt \Big]  \Big/ \bQ_r[\, \zeta\, ] \; , 
\end{equation}
Namely, Theorem \ref{thm: decomp} $(i)$ entails 
that $\bN (\, \cdot \, \big| \, D\! = r\! \, )$-a.s.~$D\! = \! r$. Then (\ref{defQr}) combined with the already mentioned 
continuity of $r\mapsto \bN (\, \cdot \, | \, \Gamma \! = \! r/2)$ easily implies that 
$r\mapsto \bN (\, \cdot \, \big| \, D\! = r\! \, )$ is weakly continuous on $\bC(\bbR_+, \bbR_+)$. 
Moreover, (\ref{desindiam}) can be rewritten as 
\begin{equation}
\label{condiame}
 \bN = \int_0^\infty  \!\!\! \bN(D\! \in \! dr)\,  \bN (\, \cdot \, | \, D \! = \! r) \; 
 \end{equation}
that is analogous to (\ref{heightcondef}). We mention that the proof of Theorem  \ref{thm: decomp} relies on the decomposition (\ref{heightcondef}) due to Abraham \& Delmas \cite{AbDe09}. \cq 
\end{rem}
\begin{rem}
\label{rereroot} It is easy to check from (\ref{rerootH}) that for all $t_0, t$, $(H^{[t]})^{[t_0]}= H^{[t+t_0]}$. Therefore,  
(\ref{desindiam}) implies that $H$ under $\bN$ is invariant under rerooting. Namely, for all measurable functions 
$F\! :\!  \bC(\bbR_+, \bbR_+)\! \rightarrow \! \bbR_+$, 
\begin{equation}
\label{renracH}
\forall t_0 \in \bbR_+, \quad \bN \big[ \un_{\{ \zeta \geq t_0\}} F\big( H^{[t_0]} \big)\big]= \bN \big[ \un_{\{ \zeta \geq t_0\}} F\big( H \big)\big] \; , 
\end{equation}
which is quite close to Proposition 2.1 in D.~\& Le Gall \cite{Duquesne09}, that is used in the proof of  Theorem  \ref{thm: decomp}. \cq 
\end{rem}
\begin{rem}
\label{Qinterest} As shown by (\ref{Hcondia}), $\bN \big( \cdot  \big|  D\! = r\! \, \big)$ is derived 
from $\bQ_r$ by a uniform rerooting. This property suggests that the law of  the compact real tree $(\cT, d)$ coded by $H$ under $\bQ_r$, \textit{without its root}, is the scaling limit of natural models of labeled \textit{unrooted} trees conditioned by their diameter. \cq 
\end{rem}

\begin{rem}
\label{Qinterestbis}
Another reason for introducing the law $\bQ_r$ is the following: we deduce from (\ref{Hcondia}) that  for all measurable functions $F\! :\!  \bC(\bbR_+, \bbR_+) \! \rightarrow \! \bbR_+$,
\begin{equation}
\label{rootdia}
\bN \big[ F( H^{[\tau_0]} ) \, \big| \, D\! = \! r \big]=   \bQ_r \big[\zeta  F( H^{[\tau_0]} ) \big] \big/  \bQ_r[\, \zeta\, ] \; , 
\end{equation} 
where $\tau_0$ is as in Theorem \ref{diamlaw}. As shown by Theorem \ref{thm: decomp} $(iv)$, 
$H$ under $\bQ_r$ enjoys a Poisson decomposition along its diameter, which is 
not the case of $H$ under $\bN(\, \cdot \, | \, D\! = \! r)$ by (\ref{rootdia}).  \cq 
\end{rem}
\paragraph{The law of the height and of the diameter of stable L\'evy trees conditioned by their total mass.}
In application of Theorem \ref{thm: decomp}, we compute the law of $\Gamma$ and $D$ under $\bN(\, \cdot \, | \, \zeta\! = \! 1)$ in the cases where $\Psi$ is a stable branching mechanism. Namely,  
we fix $\gamma \! \in \! (1, 2]$ and 
$$\Psi (\lambda)= \lambda^\gamma, \quad \lam \! \in \! \bbR_+\; , $$
that is called 
the \textit{$\gamma$-stable branching mechanism}. We first recall the definition of the law $\bN(\, \cdot \, | \, \zeta\! = \! 1)$ for such a branching mechanism. 

When $\Psi$ is $\gamma$-stable, the L\'evy process $X$ under $\bbP$ satisfies the following scaling property: for all  $r \! \in \! (0, \infty)$, 
$(r^{-\igam}\! X_{rt} )_{ t\geq 0}$ 
has the same law as 
$X$, which easily entails by (\ref{eq: defH}) that under $\bbP$, 
$( r^{-\frac{\gamma-1}{\gamma}} \! H_{rt})_{t\geq 0} $  
has the same law as $H$ and the 
Poisson decomposition (\ref{PoisdecH}) implies the following: 
\begin{equation}
\label{HscalingN}
\big( r^{-\frac{\gamma-1}{\gamma}}\! H_{rt} \big)_{t\geq 0} \quad  \textrm{under} \quad  
r^{\frac{1}{\gamma}} \,  \bN \quad  
\overset{\textrm{(law)}}{=} \quad H  \quad  \textrm{under}  \quad \bN \; . 
\end{equation}
We then easily derive from (\ref{lifetimeexc}) that 
\begin{equation}\label{eq: defp}
\bN(\zeta \! \in \! dr) \! =\! p_\gamma(r) \, dr\; , \quad \textrm{where} \quad p_\gamma(r)=c_\gamma r^{-1-\igam} \quad \textrm{with} \quad 1/c_\gamma = \gamma \Gamma_{\! e} \big(\frac{_{\gamma -1}}{^\gamma} \big) \; .
\end{equation}
Here $\Gamma_{\! e}$ stands for Euler's Gamma function. By 
(\ref{HscalingN}), there exists a family of laws on $\bC(\bbR_+, \bbR_+)$ denoted by  
$\bN(\, \cdot \, | \, \zeta \! = \! r)$, $r \! \in \! (0, \infty)$, such that 
$r \mapsto \bN(\, \cdot \, | \, \zeta \! =\!  r)$ is weakly continous on $\bC(\bbR_+, \bbR_+)$, 
such that $\bN(\, \cdot \, | \, \zeta \! = \! r)$-a.s.~$\zeta \! = \! r$ and such that 
\begin{equation}
\label{zetades}
\bN = \int_0^\infty \!\!\! \bN(\, \cdot \, | \, \zeta \! =\!  r) \, \bN (\zeta \! \in \! d r) \; .
\end{equation}
Moreover, by (\ref{HscalingN}), $\big( r^{-\frac{\gamma-1}{\gamma}} \! H_{rt} \big)_{t\geq 0}$ under 
$\bN(\, \cdot \, | \, \zeta \! = \! r)$ has the same law as $H$ under 
$\bN(\, \cdot \,| \, \zeta \! =\!  1 )$. We call $\bN(\, \cdot \,| \, \zeta \! =\!  1 )$ the \textit{normalized law of the $\gamma$-stable height process} and to simplify notation we set 
\begin{equation}
\label{Nnorm}
\bN_{\! {\rm nr}}:= \bN(\, \cdot \,| \, \zeta \! =\!  1 )
\end{equation}
Thus, for all measurable functions $F\! :\! \bC(\bbR_+, \bbR_+)\! \rightarrow \! \bbR_+$, 
\begin{equation}
\label{echtscal}
\bN \big[ F(H)\big]=  c_\gamma \int_0^\infty \!\!\! \! dr  \, r^{-1-\igam}\, \bN_{{\rm nr}} \Big[ F\Big( \big( r^{\frac{\gamma-1}{\gamma}}\! H_{t/r} \big)_{\! t\geq 0}\Big) \Big] \; . 
\end{equation}
When $\gamma \! = \! 2$, $\bN_{\! {\rm nr}}$ is, up to scaling, the normalized Brownian excursion that is, as shown by Aldous \cite{aldcrt3}, the scaling limit of the contour process of the uniform (ordered rooted) tree with $n$ vertices as $n\rightarrow \infty$;  Aldous \cite{aldcrt3} 
also extends this limit theorem to Galton-Watson trees conditioned to have $n$ vertices and 
whose offspring distribution has a
second moment. This result has been extended by D.~ \cite{Duquesne03} to Galton-Watson trees conditioned to have $n$ vertices and whose offspring distribution is in the domain of attraction of a 
$\gamma$-stable law, the limiting process being in this case the normalized excursion of the 
$\gamma$-stable height process. See also Kortchemski \cite{Kor12} 
for scaling limits of Galton-Watson tree conditioned to have $n$ leaves.  

\medskip

\begin{figure}[ht]
\begin{minipage}[l]{.45\linewidth}
      \includegraphics[scale=0.46]{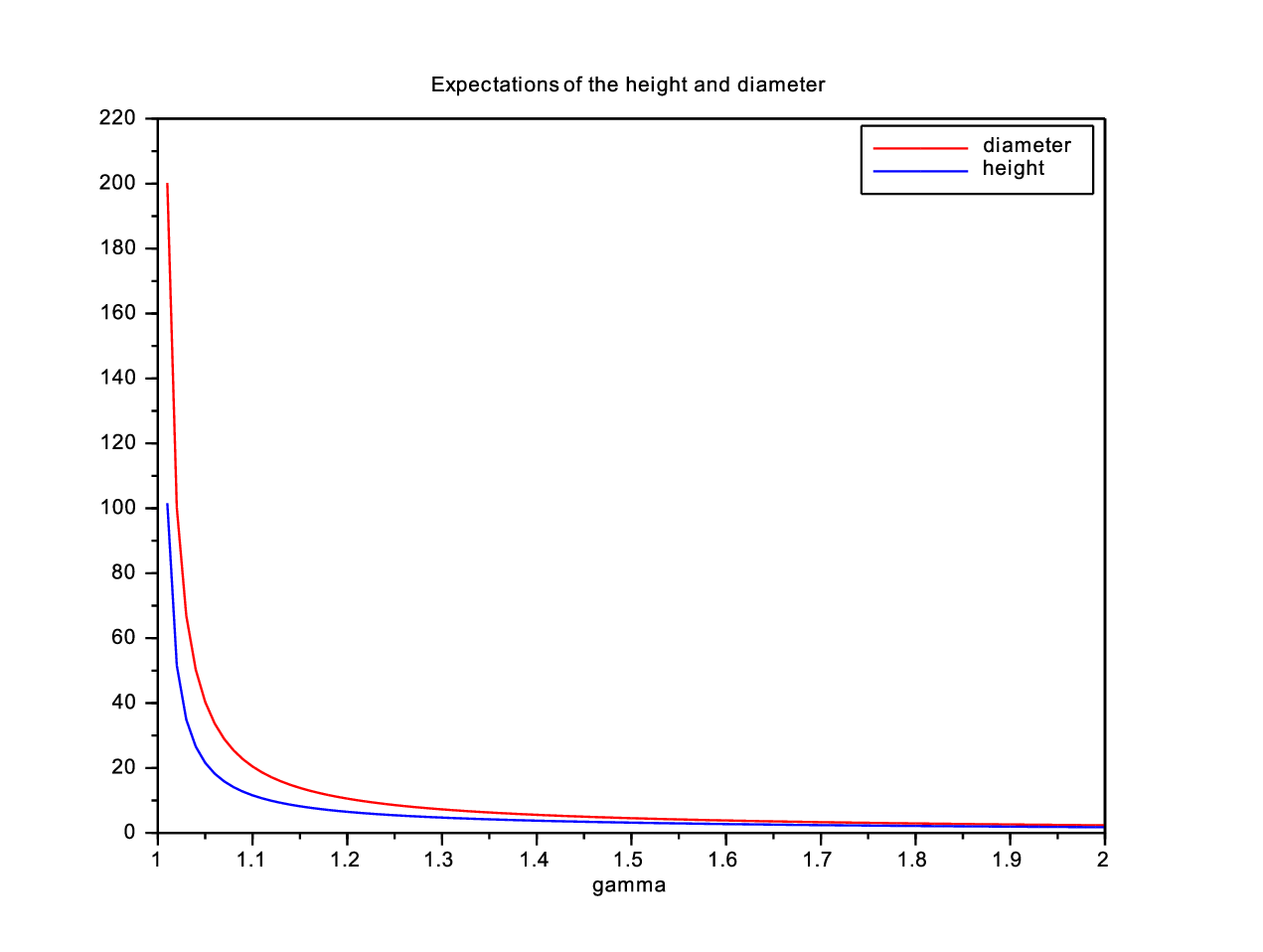}
   \end{minipage} \hfill
   \begin{minipage}[l]{.45\linewidth}
      \includegraphics[scale=0.46]{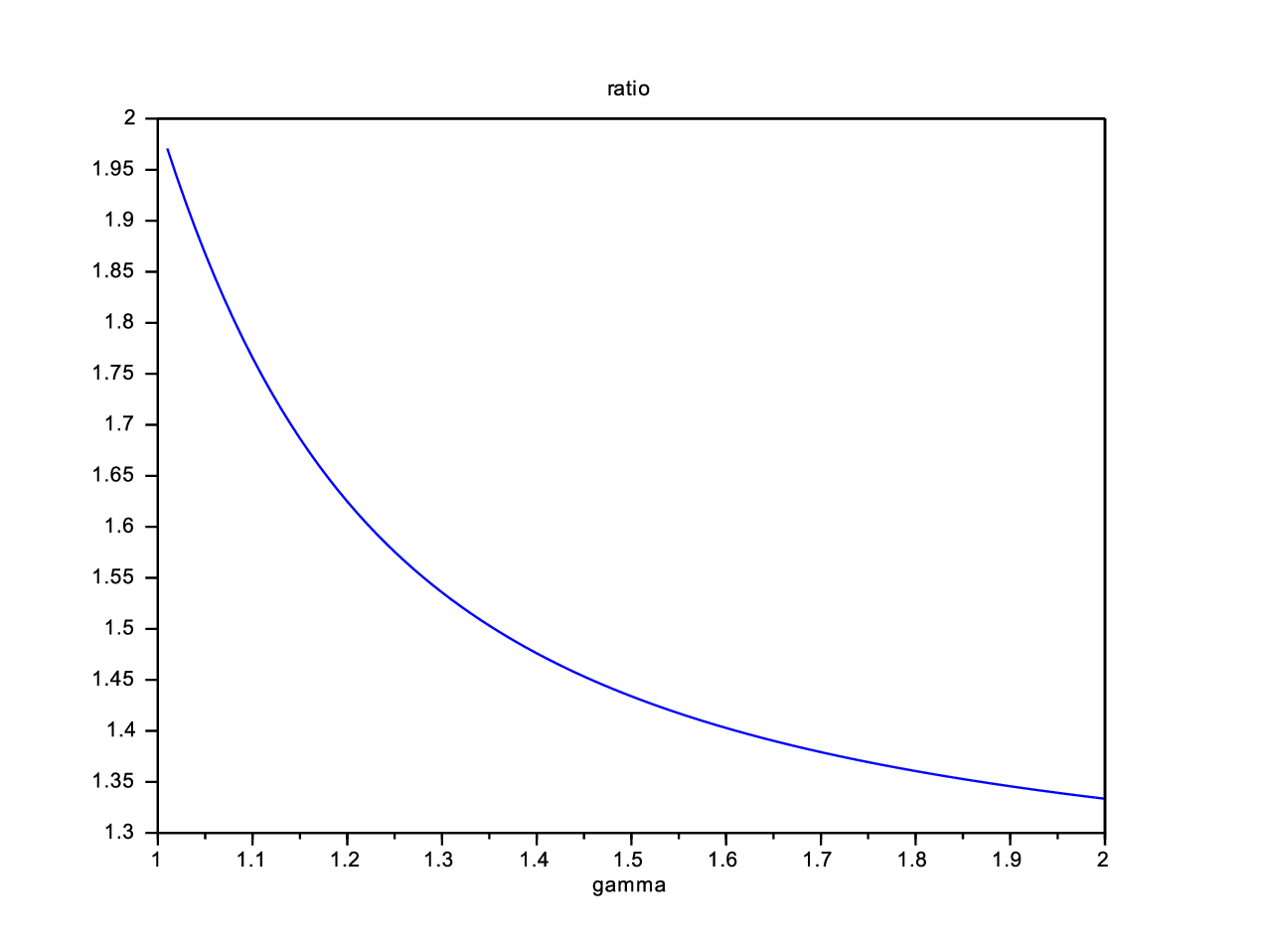}
   \end{minipage}
\caption{\label{simu}{\small \textit{\hspace{-1mm}numerical evaluations of $\nr[\Gamma]$ and $\nr[D]$ for $\gam\ino (1, 2]$. 
On the left hand side, the graphs of $\gamma \! \mapsto \! \nr[D]$ (in red) and $\gamma\mapsto \nr [\Gamma]$ (in blue). 
On the right hand side, the graph of $\gam \! \mapsto \! \nr [D]/\nr[\Gamma]$.}} }
\end{figure}

We next introduce $w\! :\!  (0, \infty) \! \rightarrow \! (1, \infty)$ that is the unique $C^\infty$ 
decreasing bijection that satisfies the following integral equation: 
\begin{equation}
\label{wdef1}
\forall y \in (0, \infty) , \qquad \int_{w(y)}^\infty \frac{du}{u^\gamma-1}=y \; .
\end{equation}
We refer to Section \ref{preprelisec} for a probabilistic interpretation of $w$ and further properties. 
The following proposition characterizes the joint law of $\Gamma$ and $D$ under 
$\bN_{\! {\rm nr}}$ by  means of Laplace transform.
\begin{prop}
\label{prop: lp}
Fix $\gam \! \in \! (1,2]$ and $\Psi (\lam)\! =\! \lam^\gam$, $\lam \! \in \! \bbR_+$. Recall from 
(\ref{Nnorm}) the definition of the law $\bN_{\! {\rm nr}}$ of 
the normalized excursion of the $\gamma$-stable height process. We then set  
\begin{equation}\label{eq: def_bL}
\forall \lam, y, z \in (0, \infty) , \quad {\rm L}_{\lambda} (y,z):= c_\gamma \int_0^\infty\!\!\!  e^{-\lambda r} 
r^{-1-\frac{1}{\gam}} \, \bN_{\! {\rm nr}} 
\big(\, r^{\frac{\gamma -1}{\gamma}} D\! >\! 2y \, ;\,  r^{\frac{\gamma -1}{\gamma}} \Gamma \! >\! z \big) \, dr \; , 
\end{equation}
where we recall from (\ref{eq: defp}) 
that $1/c_\gamma = \gamma \Gamma_{\! e} \big(\frac{{\gamma -1}}{\gamma} \big)$, $\Gamma_{\! e}$ standing for Euler's Gamma function. Note that 
\begin{equation}\label{scalebL}
\forall \lam, y, z \in (0, \infty) , \quad {\rm L}_1(y, z)= \lambda^{-\frac{1}{\gamma}} {\rm L}_\lambda \big( \lambda^{-\frac{\gamma-1}{\gamma}} y\, ,\,  \lambda^{-\frac{\gamma-1}{\gamma}} z \big)   \; .
\end{equation}
Recall from (\ref{wdef1}) the definition of $w$. Then, 
\begin{equation}
\label{Lreecrit}
{\rm L}_1(y, z)
\! = \!  w ( y \!  \vee \!   z) -1 
-\frac{_1}{^\gamma} \un_{\{ z < 2y\}}\!  \big( w(y)^\gamma \! -\! 1 \big)^{\! 2} \! \! \left( \! \frac{w \big( y \! \wedge \! (2y\! -\! z )\big)}{w \big( y\!  \wedge \! (2y\! -\! z ) \big)^\gamma \! -\! 1} -(\gamma\! -\! 1) \big(y \! \wedge \! (2y\! -\! z ) \big) \right) . 
\end{equation} 
In particular, for all $y, z\! \in \! (0, \infty)$,  
\begin{equation}
\label{htdianor}
{\rm L}_1 (0, z)= w(z)-1 \quad \textrm{and} \quad {\rm L}_1(y, 0) = w(y)-1 -\frac{_1}{^\gamma} \big( w(y)^\gamma \! -\! 1  \big) \Big( w(y) \! -\! (\gamma \! -\! 1) y \big( w(y)^\gamma \! -\! 1\big)  \Big) . 
\end{equation}
\end{prop}
\begin{rem}
\label{expechtdm} Proposition \ref{prop: lp} allows explicit computations of $\bN_{{\rm nr}} [\Gamma]$ and $\bN_{{\rm nr}} [D]$ in terms of $\gamma$: we refer to Proposition \ref{expcht} and Proposition 
\ref{expcdm} in Section \ref{Pfexpc} for precise results. In the Brownian case $\gamma \! = \! 2$, we recover that $\nr [\Gamma]\! =\! \sqrt\pi$ and $\nr [D]\! =\! \frac{4}{3}\sqrt\pi$, therefore $\nr[D]/\nr [\Gam]\! =\! \frac{4}{3}$. This ratio between the height and diameter of the Brownian tree is first observed in \cite{Sz83} and later Aldous gives an explanation of this fact in \cite{aldcrt2}. In the non-Brownian stable cases this explanation breaks down: as a consequence of Proposition \ref{expcht} and Proposition 
\ref{expcdm}, as $\gamma \! \to \! 1+$, we prove that 
\begin{equation}
\label{siuccd}
\nr [\Gamma]\! =\! \frac{1}{\gamma \! -\! 1} + \gam_{{\rm e}} \! +\! 1+  \mathcal{O} \big(\gam \! -\! 1\big) 
\quad \textrm{and} \quad  \nr [D]\! =\! \frac{2}{\gamma \! -\! 1} +2\gam_{{\rm e}} \! -\! 1+ \mathcal{O} \big( \gam \! -\! 1\big)\; , 
\end{equation}
where $\gam_{{\rm e}}$ stands for the Euler-Mascheroni constant. Thus, $\lim_{\gam \to 1+}  \nr [D] / \nr [\Gamma]\! = \! 2$. 
See Figure \ref{simu}.  
We refer to Section \ref{Pfexpc} for more details. \cq 
\end{rem}

\medskip

Proposition \ref{prop: lp} is known in the 
Brownian case, where $w(y)\! = \! \coth (y)$: see W.~\cite{Wang14} for the joint law. In the Brownian case, standard computations derived from (\ref{htdianor}) imply 
the following power expansions that hold true for all $r\! \in \! (0, \infty)$: 
\begin{equation}
\label{htBronor}
\bN_{\! {\rm nr}} \big(\Gamma \! >\! r \big)=2\sum_{n\geq 1}\(2n^2r^2 \! -\! 1\)e^{-n^2r^2} 
\end{equation}
and 
\begin{equation}
\label{diaBronor}
\bN_{\! {\rm nr}} \big( D \! >\! r\big) =  \sum_{n\geq 2} (n^2-1)\big(\frac{_1}{^6}n^4r^4-2n^2r^2+2\big)e^{-n^2r^2/4}   \; .
\end{equation}
These results can be derived from expressions in Szekeres \cite{Sz83} (see also W.~\cite{Wang14} for more details). 

  We next provide similar asymptotic expansions in the non-Brownian stable cases. 
To that end, we introduce $s_\gamma \!  :\!  (0, \infty) \! \to \! (0, \infty)$ as the continuous version of the density of the spectrally positive $\frac{\gamma-1}{\gamma}$-stable distribution; more precisely, $s_\gamma$ is characterized by the following:  
\begin{equation}
\label{denstable}
\forall \lambda \in \bbR_+, \qquad \int_0^\infty \!\!\! e^{-\lambda x}\! s_\gamma (x) \, dx=\exp(-\gamma\lambda^{\frac{\gamma -1}{\gamma}}) \; .
\end{equation}
The following asymptotic expansion of $s_\gamma$ at $0$ is due to Zolotarev (see Theorem 2.5.2 \cite{Zol}): for all integer $N \! \geq \! 1$, 
\begin{equation}\label{eq: est}
\big(2\pi (1\! -\! \frac{_1}{^\gamma} )\big)^{\frac{1}{2}}  \, x^{\frac{\gamma+1}{2}} e^{1/x^{\gamma-1}} \! 
s_\gamma \big((\gamma\! -\! 1)x \big)= 1+\! \! \sum_{1\leq n< N} \!\!\! S_n \, x^{n(\gamma -1)} +  \cO_{N, \gamma} \big( x^{N(\gamma-1)} \big),  
\quad \textrm{as $x\to 0$.}
\end{equation}
Here $ \cO_{N, \gamma}$ means that the expansion depends on $N$ and $\gamma$. Note that $S_n$ depends on $n$ and $\gamma$ but we skip the dependence in $\gamma$ to simplify notation. 
\begin{rem} 
\label{BrocoefSn} In the Brownian case where $\gamma \! = \! 2$, it is well-known that 
 $$ s_2 (x)= \pi^{-\frac{1}{2}}  x^{-\frac{3}{2}} e^{-1/x} , \quad x\in \bbR_+ \; .$$
Then, $S_0\! = \! 1$ and $S_n \! = \! 0$, for all $n\! \geq \! 1$. \cq 
\end{rem}
For generic $\gamma \! \in \! (1, 2)$,  this asymptotic expansion does not yield a converging power expansion (although it is the case if $\gamma \! = \! 2$). See Section \ref{preliasy} for more details on $s_\gamma$. 
To state our result we first need to introduce an auxiliary function derived from $s_\gamma$ as follows. 
\begin{prop}
\label{thetproo} Let $\gamma \! \in \! (1, 2]$. Recall from (\ref{denstable}) the definition of $s_\gamma$. We introduce the following function:
 \begin{equation}
\label{ethetadef}
\forall x \in \bbR_+, \quad \theta(x) \!: =\! (\gamma \! -\! 1) \, x^{-1} \! s_\gamma (x)
- \frac{_{\gamma  - 1}}{^\gamma} x^{-1-\frac{1}{\gamma}}\!\!  \int_0^x \!\!\! dy \, y^{\frac{1}{\gamma} -1} 
s_\gamma (y) \; .
\end{equation}
Then, the following holds true.  
\begin{itemize}
\item[(i)] $\theta$ is well-defined, continuous, 
\begin{equation}
\label{eLapltheta}
\int_0^\infty \!\!\! dx \,  |\theta (x)|  <  \infty \quad \textrm{and} \quad \int_0^\infty \!\!\! dx \, e^{-\lambda x} \theta(x) = \lambda^{\frac{1}{\gamma}} e^{-\gamma \lambda^{\frac{\gamma -1}{\gamma}}}, 
\quad \lambda \! \in \! \bbR_+.   
\end{equation}
\item[(ii)] Recall from (\ref{eq: est}) the definition of the sequence $(S_n)_{n\geq 0}$, with $S_0\! = \! 1$.  
Let $(V_n)_{n\geq 0}$ be a sequence of real numbers recursively defined by $V_0\! = \! 1 $ and  
\begin{equation}
\label{erecVp}
\forall n \! \in \! \bbN, \qquad V_{n+1}=  S_{n+1} +  \big(n \!-\! \frac{_1}{^2} \! -\! \frac{_1}{^{\gamma -1}}  \big) S_n - 
\big( n  \!-\! \frac{_1}{^2} \! -\! \frac{_1}{^{\gamma }}  \big) V_n \; . 
\end{equation}
Then, for all integers $N \! \geq \! 1$, 
\begin{equation}\label{eq: esth}
\big(2\pi (1\! - \frac{_1}{^\gamma} )\big)^{\frac{1}{2}}  \, x^{\frac{\gamma+3}{2}} e^{1/x^{\gamma-1}}  
\theta \big((\gamma\! -\! 1)x \big)= 1+\! \! \!\! \sum_{1\leq n< N} \!\!\! V_n \, x^{n(\gamma -1)} \; +  \cO_{N, \gamma} \big( x^{N(\gamma-1)} \big),  
\end{equation}
as $x\to 0$.
\end{itemize} 
\end{prop}
We use $\theta$ to get the asymptotic expansion of the law of the total height of the normalized $\gamma$-stable tree as follows.
\begin{thm}\label{thm: ht_asy}
Let $\gamma\in (1,2]$. We introduce the following function:
\begin{equation}
\label{edefxi}
\forall r \in \bbR_+, \quad \xi(r):= r^{-\frac{\gamma+1}{\gamma-1}} \theta \big( r^{-\frac{\gamma }{\gamma-1}} \big) \; . 
\end{equation}
where $\theta$ is defined in (\ref{ethetadef}). Then, there exists a real valued sequence $(\beta_n)_{n\geq 1}$ and $x_1 \! \in \! (0, \infty)$ such that 
\begin{equation}  
\label{propbetan}
\sum_{n\geq 1} |\beta_n| x_1^n \; < \infty \quad \textrm{and} \quad  \forall r \in  (0, \infty) , \quad \sum_{n\geq 1} |\beta_n| \sup_{s\in [r, \infty)} |\xi(ns)| \; < \infty  \; , 
\end{equation}
and such that 
\begin{equation}
\label{exactoht}
\forall r \in (0, \infty), \quad  c_\gamma \bN_{\! {\rm nr}} \big( \Gamma \! >\! r\big)  = \sum_{n\geq 1} \beta_n \, \xi (nr) \; , 
\end{equation}
where we recall from (\ref{eq: defp}) that $1/c_\gamma\! = \! 
\gamma \Gamma_{\! e} \big(\frac{{\gamma -1}}{\gamma} \big)$, $\Gamma_{\! e}$ standing for Euler's gamma function.  
Moreover, for all integers $N \! \geq \! 1$, as $r\to \infty$, 
\begin{equation}\label{eq: ht_asy}
\frac{1}{C_1}\, r^{-1-\frac{\gamma}{2}} \, e^{r^\gamma} \, \bN_{\! {\rm nr}} \Big( \Gamma \! >\! r (\gamma \! -\! 1)^{-\frac{\gamma -1}{\gamma}} \Big)\,  = \, 1+ \!\! \sum_{1\leq n<N} \!\!\! V_n \, r^{-n \gamma} \; + \cO_{\! N, \gamma } \big( r^{-N\gamma }\big)\; , 
\end{equation}
where $C_1:= (2\pi)^{-\frac{1}{2}} (\gamma -1)^{\frac{1}{2}+ \frac{1}{\gamma}} 
\gamma^{\frac{3}{2}}\,  \Gamma_{\! e} ( \frac{\gamma -1}{\gamma})  \exp (C_0) $, where 
\begin{equation}\label{eq: Cgam}
C_0:=\gamma \int_1^\infty \frac{du}{(u+1)^\gam -1}- \int_0^1 \frac{du}{u}\,  \frac{(u+1)^\gam-1-\gam u}{(u+1)^\gam-1}\; , 
\end{equation}
and where the sequence $(V_n)_{n\geq 1}$ is recursively defined by (\ref{erecVp}) in Proposition \ref{thetproo}.  
\end{thm}
\begin{rem}
\label{betarapid}
The convergence in (\ref{exactoht}) is rapid. Indeed, by (\ref{eq: esth}), we see that $\xi (nr)$ is of order 
$$(nr)^{1+\frac{\gamma}{2}} \exp (-n^\gamma (\gamma \! -\! 1)^{\gamma -1} r^\gamma ) \; .$$ 
Then, the asymptotic expansion (\ref{eq: ht_asy}) is that of the first term of (\ref{exactoht}) that is $c_\gamma^{-1} \! 
\beta_1\,  \xi(r)$. \cq
\end{rem}
\begin{rem}
\label{betanou}
The definition of the sequence $(\beta_n)_{n\geq 0}$ is involved: see Lemma \ref{vraiexphi} and its proof for a precise definition. 
However, in the Brownian case, everything can be explicitly computed: for all $n\! \geq 1\!$, $\beta_n\! = \! 2$, 
$ \xi(r)= (4\pi)^{-\frac{1}{2}} (2r^2 -1)e^{-r^2}$, $c_2\! = \! (4\pi)^{-\frac{1}{2}} $, and we recover (\ref{htBronor}) from (\ref{exactoht}); moreover, $C_0\! = \! \log 2$, $C_1 \! = \! 4$, $V_0\! = \! 1$, $V_1\! = \! -\frac{1}{2}$ and $V_n\! = \! 0$, for all $n\! \geq \! 2$. \cq 
 \end{rem}
To state the result concerning the diameter, we need precise results on the derivative of the $\frac{\gamma-1}{\gamma}$-stable density. 
\begin{prop}
\label{derivss} 
Let $\gamma \! \in \! (1, 2]$. Recall from (\ref{denstable}) the definition of the density $s_\gamma$. 
Then $s_\gamma$ is $C^1$ on $\bbR_+$, 
\begin{equation}
\label{eLaplsprim}
\int_0^\infty\!\!\! dx \,  |s^\prime_\gamma (x)| < \infty  \qquad \textrm{and} \qquad \int_0^\infty \!\!\! dx \, e^{-\lambda x} s_\gamma^\prime (x)= \lambda e^{-\gamma \lambda^{\frac{\gamma -1}{\gamma}}}, 
\quad \lambda \! \in \! \bbR_+.   
\end{equation}
Moreover, $s^\prime_\gamma$ has the following asymptotic expansion: recall from (\ref{eq: est}) the definition of the sequence $(S_n)_{n\geq 0}$, with $S_0\! = \! 1$;  
let $(T_n)_{n\geq 0}$ be a sequence of real numbers recursively defined by $T_0\! = \! 1 $ and 
\begin{equation}
\label{eTdedS}
\forall n\in \bbN, \qquad T_{n+1}:=  S_{n+1}+ \big( n-\frac{_1}{^2} -\frac{_1}{^{\gamma-1}}\big) S_{n} \; .
\end{equation} 
Then, for all positive integers $N$, we have 
\begin{equation}
\label{asprim}
\big(2\pi (1\! - \frac{_1}{^\gamma} )\big)^{\frac{1}{2}}  \, x^{\frac{3\gamma+1}{2}} e^{1/x^{\gamma-1}}  
s_\gamma^\prime \big((\gamma\! -\! 1)x \big)= 1+\! \! \!\! \sum_{1\leq n< N} \!\!\! T_n \, x^{n(\gamma -1)} \; +  \cO_{N, \gamma} \big( x^{N(\gamma-1)} \big), 
\end{equation}
 as $x\rightarrow 0$. 
 \end{prop}
The asymptotic expansion of the law of the diameter of the normalized $\gamma$-stable tree 
is then given in the following theorem. 
\begin{thm}
\label{thm: dm_asy}
Let $\gamma \! \in \! (1,2]$. Recall from (\ref{edefxi}) the definition of the function $\xi$. We also introduce the following function: 
\begin{equation}
\label{edefxibar}
\forall r \in \bbR_+, \quad \overline{\xi}(r):= r^{-\frac{\gamma+1}{\gamma-1}} s^\prime_\gamma \big( r^{-\frac{\gamma }{\gamma-1}} \big) \; , 
\end{equation}
where $s^\prime_\gamma$ is the derivative of the density $s_\gamma$ defined in (\ref{denstable}). Then there exist 
two real valued sequences $(\gamma_n)_{n\geq 2} $ and $(\delta_n)_{n\geq 2}$ and $x_2 \! \in \! (0, \infty)$ such that  
\begin{equation}  
\label{propgadeln}
\sum_{n\geq 2} \, (|\gamma_n| +|\delta_n|) x_2^n  <\!  \infty \quad \textrm{and} \quad  \forall r \in (0, \infty) , \quad 
\sum_{n\geq 2}  |\gamma_n| \!\!\! \sup_{\; s\in [r, \infty)} \!\!\!  | \overline{\xi}(ns)|  + |\delta_n | \!\!\! 
\sup_{\; s\in [r, \infty)} \!\! | \xi(ns)|  < \infty \; ,  
\end{equation}
and such that 
\begin{equation}
\label{exactodiam}
\forall r \in (0, \infty), \quad  c_\gamma \bN_{\! {\rm nr}} \big( D \! >\! 2r\big)  = \sum_{n\geq 2} 
\gamma_n \overline{\xi} (nr)+   \delta_n \xi (nr)  \; , 
\end{equation}
where we recall from (\ref{eq: defp}) that $1/c_\gamma\! = \! 
\gamma \Gamma_{\! e} \big(\frac{{\gamma -1}}{\gamma} \big)$, $\Gamma_{\! e}$ standing for Euler's gamma function.  
Moreover, for all integers $N \! \geq \! 1$, as $r\to \infty$, 
\begin{equation}\label{eq: dm_asy} 
\frac{1}{C_2}\, r^{-1-\frac{3\gamma}{2}} \, e^{r^\gamma} \, \bN_{\! {\rm nr}} \Big(D \! >\! r (\gamma \! -\! 1)^{-\frac{\gamma -1}{\gamma}} \Big)\,  = \, 1+ \!\! \sum_{1\leq n<N} \!\!\! U_n \, r^{-n \gamma} \; + \cO_{ N, \gamma } \big( r^{-N\gamma }\big)\; ,   
\end{equation}
where $C_2\! := \! (8\pi)^{-\frac{1}{2}} (\gamma \! -\! 1)^{\frac{3}{2} + \frac{1}{\gamma}} \gamma^{\frac{5}{2}} \,  \Gamma_{\! e} (\frac{\gamma-1}{\gamma})  \exp (2C_0) $, where $C_0$ is defined by (\ref{eq: Cgam}) and where 
the sequence $(U_n)_{ n\geq 1}$ is recursively defined by $U_0\! = \! 1$ and 
\begin{equation}
\label{Undefrec}
\forall n\geq 1 \,,  \qquad U_n= T_n-\frac{_{\gamma+1}}{^{\gamma (\gamma -1)}} V_{n-1}\; . 
\end{equation}
Here $(T_n)_{n\geq 0}$ is defined by (\ref{eTdedS}) and $(V_n)_{n\geq 0}$ is defined by (\ref{erecVp}).  
\end{thm}

\begin{rem}
\label{gadelnou}
The convergence in (\ref{exactodiam}) is rapid. Indeed, by (\ref{asprim}) and (\ref{eq: esth}) we see that $\overline{\xi} (nr/2)$ and $\xi (nr/2)$ are of respective 
order 
$$ (nr)^{1+\frac{3\gamma}{2}} \exp (-n^\gamma 2^{-\gamma}(\gamma \! -\! 1)^{\gamma -1} r^\gamma ) \quad \textrm{and} \quad (nr)^{1+\frac{\gamma}{2}} \exp (-n^\gamma 2^{-\gamma}(\gamma \! -\! 1)^{{\gamma -1}} r^\gamma ) \; .$$
Then the asymptotic expansion (\ref{eq: dm_asy}) is that of $c_\gamma^{-1} \gamma_2 \,  \overline{\xi}(r)+c_\gamma^{-1} \delta_2 \,  \xi(r) $. \cq
\end{rem}

\begin{rem}
\label{gadelnou}
The definitions of the sequences $(\gamma_n)_{n\geq 0}$ and $(\delta_n)_{n\geq 0}$ are involved: see the proof of Lemma \ref{asfundia} for a precise definition. However, in the Brownian case, everything can be computed explicitly:  
$$ \forall n\geq 2, \quad \gamma_n= \frac{_4}{^3} (n^2-1), \quad \delta_n= -2 (n^2-1) \quad \textrm{and} \quad \overline{\xi} (r)= \pi^{-\frac{1}{2}} r^2 \big( r^2-\frac{_3}{^2}\big) e^{-r^2}\; , $$
which allows to recover (\ref{diaBronor}) from (\ref{exactodiam}). Moreover, $C_2 \! =\!  8$, $U_0\! = \! 1$, $U_1\! = \! -3$, $U_2\! = \! -\frac{3}{4}$ and $U_n \! = \! 0$, for all 
$n\! \geq \! 3$. \cq
\end{rem}

\paragraph{The tail at $0+$ of the law of the total height and of the diameter of the normalised stable tree.} In the Brownian case $\gamma \! = \! 2$, 
it is not straightforward to derive from (\ref{htBronor}) and (\ref{diaBronor}) an asymptotic expansion of $\bN_{\! {\rm nr}} \big( \Gamma \! \leq \! r\big) $ and $ \bN_{\! {\rm nr}} \big( D \! \leq \! r\big)$ when $r\! \rightarrow \! 0$. To that end, we use Jacobi's identity on theta functions and we get 
\begin{equation}
\label{bhtBronor} 
\bN_{\! {\rm nr}} \big(\Gamma \! \leq \! r \big)=
\frac{4\pi^{5/2}}{r^{3}}
\sum_{n \geq 1}n^2e^{-n^2\pi^2/r^2} \! \underset{r\rightarrow 0}{\sim} 4\pi^{\frac{5}{2}}r^{-3} e^{-\frac{\pi^2}{r^2}}
\end{equation}
and 
\begin{equation}
\label{bdiaBronor} \bN_{\! {\rm nr}} \big( D \! \leq \! r\big) \! = \! 
\frac{\sqrt{\pi}}{3} \sum_{n\geq 1} \Big(  \frac{8}{r^3} \big( 24a_{n,r} \! -\! 36 a_{n,r}^2 +8a_{n,r}^3 \big) +\frac{16}{r} a_{n,r}^2 \Big) e^{-a_{n,r}}\!\! \underset{r\rightarrow 0}{\sim}\! \tfrac{1}{3}\, 2^{12}\pi^{\frac{13}{2}}r^{-9} \! e^{-\frac{4\pi^2}{r^2}}  , 
\end{equation}
where we have set $a_{n,r}= 4 (\pi n/r)^2$ for all $r \! \in \! (0, \infty)$ and for all $n\! \geq \! 1$. See Szekeres \cite{Sz83} and Aldous \cite{aldcrt2} for more detail and see W.~\cite{Wang14} for the joint of $D$ and $\Gamma$ in the Brownian case.

  In the non-Brownian stable cases, when $\gamma  \! \in \!  (1, 2)$, the asymptotic expansions (\ref{exactoht}) in Theorem \ref{thm: ht_asy} and (\ref{propgadeln}) in Theorem \ref{thm: dm_asy} are useless to get asymptotics of $\bN_{\! {\rm nr}} \big( \Gamma \! \leq \! r\big) $ and $ \bN_{\! {\rm nr}} \big( D \! \leq \! r\big)$ when $r\! \rightarrow \! 0$. In these cases we only prove the following result. 
\begin{thm} 
\label{htdm0} We fix $\gamma \ino (1, 2)$ (in particular, $\gamma\! \neq \! 2$). Then, as $r\! \rightarrow \! 0+$, 
\begin{align}
\label{dtcxtqcd}
& \nr(\Gamma \! \le \! r)  \sim C r^{\gam+2+\frac{1}{\gam-1}}\exp\big(\! -\! \llcr \, r ^{-\frac{\gam}{\gam-1}}\big) \\
\label{frqvutgc}
 \textrm{and}  \quad &  \nr(D \! \le \! 2r) \sim C^\prime r^{\gam+1}\exp\big(\! -\! \llcr \, r ^{-\frac{\gam}{\gam-1}}\big)\; ,  \\ 
 \label{lamstar}
\textrm{where} \quad  \llcr \! := \! & \Big( \frac{{\pi/\gamma}}{{\sin (\pi / \gamma)}}\Big)^{\frac{\gamma}{\gamma-1}}
\, , \quad 
 C\! := \frac{(\gam\! -\! 1)^{\gam+2}\Gamma_{\! e} (1\! -\! \frac{1}{\gam})}{\gam^{\gam-1}\llcr\Gamma_{\! e} (2\! -\! \gam)} \quad \textrm{and} \quad C^\prime \! := \! 
 2\llcr C \; .
\end{align}
\end{thm}
In table \ref{table}, we summarize the exponents of the tail probabilities for the total height and the diameter in the different asymptotic regimes. 
We make two remarks. 
\begin{table}[ht]
\begin{center}
\begin{tabular}[c]{|c| |c | c | c| }
\hline
&   & $\gamma\ino (1, 2) $ & $ \gamma \! =\! 2 $ \\
\hline      
$r \! \to \! \infty$ & $  \! -\! \log \nr \Big({\color{black}\Gamma} \! >\! r \Big) \sim $& $ (\gamma\! -\! 1)^{\gamma-1}r^\gamma$  &$ r^2$  \\
\hline
$r \! \to \! \infty$ & $  \! -\! \log \nr \Big({\color{black}D} \! >\! r \Big) \sim $& $ (\gamma\! -\! 1)^{\gamma-1}r^\gamma$  &$ r^2$  \\   
\hline
$r \! \to \! 0+$ & $ \!  - \! \log \nr \Big( {\color{black}\Gamma} \! \le\!  r \Big) \sim $ & $ \big(\frac{\gamma\sin(\pi/\gamma)}{\pi} \, r\big)^{_{-\frac{\gamma}{\gamma-1}}}$ & $\pi^2 \! /r^2$\\
\hline
$r \! \to \! 0+$ & $  \!  - \! \log \nr \Big ({\color{black}D}\! \le\!  r \Big) \sim $ & $ \big(\frac{\gamma\sin(\pi/\gamma)}{2\pi}\,  r \big)^{_{-\frac{\gamma}{\gamma-1}}}$ & $ 4\pi^2 \! /  r^2$\\
\hline
\end{tabular}
\end{center}
\caption{
\label{table} 
{\small \textit{\hspace{-1mm}asymptotic exponents for the height and the diameter of stable trees.}} }
\end{table}
\begin{rem}
\label{fljviu}
First note that $\! - \! \log \nr(\Gamma \! >\! r) \! \sim \! -\! \log \nr(D\! >\! r)$ as $r\! \to \! \infty$, 
while $ \! -\! \log \nr(\Gamma \! \le \! r) \! \sim \! -\! \log \nr(D\! \le \! 2r)$ as $r \! \to \! 0$. This can be explained informally as follows: roughly speaking, 
theorem \ref{thm: decomp} asserts that a stable tree conditioned by its total diameter 
$D$ is obtained by glueing at their roots two independent trees conditioned to have height $D/2$, the root is uniformly chosen according to the mass measure in the resulting tree and 
the height is the distance of the root from the most distant extremity of the diameter. When $r$ is large, one of the two trees has a much larger mass that is concentrated near its height, thus the root is close to one of the extremities of the diameter and $\Gamma$ is comparable to $D$. When $r$ is small, both trees have a comparable mass that is concentrated near their root (corresponding to the midpoint of the diameter). 
So the root of the tree conditioned by its diameter is close to the midpoint of the diameter and 
$\Gamma$ is comparable to $D/2$.  It is possible to make these observations rigorous by an argument based on Proposition \ref{prop: lp}. In the Brownian case, they are easily derived from the expressions for the joint law of $\Gamma$ and $D$ given in W.~\cite{Wang14}. \cq 
\end{rem}
\begin{rem}
\label{tvpgnxk} In the asymptotic regime $r \! \to \! 0+$, there is a discontinuity of the exponents as $\gamma\! \to \! 2$. 
This comes from the fact that $-\llcr$, as defined by (\ref{lamstar}) is a singular point of 
the continuation extension of $\lambda\!  \mapsto \! \mathrm L_\lam( 0, 1)$ when $\gamma \ino (1, 2)$, which is not the case when $\gamma\! = \! 2$: for more details, we refer to the proof of theorem \ref{htdm0} and Remark \ref{fufluns}. \cq 
\end{rem}

The paper is organized as follows. Section \ref{pfsecdiamdec} is devoted to the proof of Theorem \ref{diamlaw} and of Theorem \ref{thm: decomp}: in Section \ref{rtreesec}, we discuss an important geometric property of the diameter of real trees (Lemma \ref{lem: dm}) and we explain the spinal decomposition according to the total height, the result of Abraham \& Delmas \cite{AbDe09} being recalled in Section \ref{pfth1th2sec} where the proofs of  Theorem \ref{diamlaw} and Theorem \ref{thm: decomp} are actually given. 
Proposition \ref{prop: lp}, that characterizes the joint law of the total height and the diameter of normalized stable trees, is proved in Section \ref{sec: stable}. Theorem \ref{thm: ht_asy} and Theorem \ref{thm: dm_asy} are proved in Section \ref{pfasysec}. Theorem \ref{htdm0} is proved in Section \ref{Pfhtdm0}. There is  an appendix in two parts: the first part is devoted to the proof of a technical lemma (Lemma \ref{mesdec2}); the second part briefly recalls various results in complex analysis that are used in the proof of Theorems \ref{thm: ht_asy}, Theorem \ref{thm: dm_asy} and Theorem \ref{htdm0}.

\medskip

\noi
\textbf{Acknowledgment.} We are grateful to P.~Flajolet for questioning us about a possible asymptotic expansion of the right tail of the height of normalized non-Brownian stable trees. 
We want to thank N.~Curien and I.~Kortchemski for questions 
concerning the total height of stable trees when $\gamma$ goes to $1$. We warmly thank Z.~Shi for mentioning to us 
the modification of the Ikehara-Ingham Tauberian Theorem used in his paper \cite{HuShi97} written with Y.~Hu (this version of Ikehara-Ingham Tauberian Theorem is a key ingredient of the proof of Theorem \ref{htdm0}). 
We are also grateful to S.~Janson and R.~Abraham for valuable comments that improved 
a first version of the manuscript.

\section{Proof of the diameter decomposition.}
\label{pfsecdiamdec}
\subsection{Geometric properties of the diameter of real trees; height decomposition.}
\label{rtreesec}
In this section we gather deterministic results on real trees and their coding functions: we first prove a key lemma on the diameter of real trees; we next discuss how to reconstruct 
the coding function $H$ from a spinal decomposition $\cM_{0, t} (H)$, under a specific assumption on the mass measure $\bbm_H$ on $\cT_H$; then we discuss a decomposition related to the total height.     

\paragraph{Total height and diameter of compact rooted real trees.} The following result connects the total 
height and the diameter of a compact rooted real tree. 
\begin{lem}
\label{lem: dm} Let $(T, d, \rho)$ be a compact rooted real tree. We denote by $\Gamma$ and 
$D$ resp.~its total height and its diameter: $\Gamma\! :=\! \sup_{\sigma \in T} d(\rho ,\sigma)$ and $D\! = \! \sup_{\sigma, \sigma^\prime \in T} d(\sigma, \sigma^\prime)$. Then, the following holds true. 
\begin{itemize}
\item[(i)] There exist $\sigma, \sigma_0, \sigma_1 \! \in \! T$, such that 
$\Gamma \! = \! d(\rho, \sigma)$ and $D\! = \! d(\sigma_0, \sigma_1 )$. This entails 
\begin{equation}
\label{eq: bdD}
\Gamma \leq D \leq 2 \Gamma \; .
\end{equation}
\item[(ii)] Let $\sigma_0, \sigma_1 \! \in \! T$ be such that $D\! = \! d(\sigma_0, \sigma_1)$. Then, $\max \! \big( d(\rho, \sigma_0); d(\rho, \sigma_1) \big)\! =\!  \Gamma$. 
\end{itemize}
\end{lem}

\noi
\textbf{Proof.} First note that $\gamma \! \in \! T \mapsto d(\rho, \gamma)$ and $(\gamma, \gamma^\prime) \! \in \! T^2 \mapsto d(\gamma, \gamma^\prime)$ 
are real valued continuous functions defined on compact spaces; basic topological arguments entail the existence of $\sigma, \sigma_0, \sigma_1 \! \in \! T$ as in $(i)$. The inequality $\Gamma \! \leq \! D$ is an immediate consequence of the definitions of $\Gamma$ and $D$. The triangle inequality next entails that $D\! \leq \! d(\sigma_0, \rho)+ d(\rho, \sigma_1) \! \leq \! 2\Gamma$, which completes the proof of (\ref{eq: bdD}) and of ($i$). 

Let $\sigma, \sigma_0, \sigma_1 \! \in \! T$ be as in $(i)$. By the four-point condition (\ref{fourpoint}) and basic inequalities, we get 
\begin{eqnarray*}
\Gamma+D= d(\rho, \sigma) + d(\sigma_0, \sigma_1)& \leq&  \max \big( d(\rho, \sigma_0)+ d(\sigma, \sigma_1) \, ; \,  d(\rho, \sigma_1)+ d(\sigma, \sigma_0) \big) \\
& \leq & \max \big( d(\rho, \sigma_0); d(\rho, \sigma_1) \big) + \max \big( d(\sigma, \sigma_1); d(\sigma, \sigma_0) \big)\; . 
\end{eqnarray*}
If $  \max \! \big( d(\rho, \sigma_0); d(\rho, \sigma_1) \big)\! <\!  \Gamma $, then the previous inequality implies that $D \! < \! \max\!  \big( d(\sigma, \sigma_1); d(\sigma, \sigma_0) \big)$, which is absurd. Thus, 
$\max \! \big( d(\rho, \sigma_0); d(\rho, \sigma_1) \big)\! =\!  \Gamma$. \cqfd

\paragraph{Coding functions and their spinal decompositions.} 
Recall that $\mathbf{0}$ stands for the null function of $\bC(\bbR_+, \bbR_+)$. We denote by 
$\bC_c(\bbR_+, \bbR_+)$ the functions of $\bC(\bbR_+, \bbR_+)$ with compact support. 

\begin{defi}
\label{defExcc}
We introduce the set of coding functions: 
\begin{equation}
\label{defExc}
\mathtt{Exc}= \big\{ H \! \in \! \bC_c(\bbR_+, \bbR_+) : \; 
\textrm{$H_0\! = \! 0$, ~$H \! \neq \mathbf{0}$, $\bbm_H$ is diffuse and ~$\bbm_H (\cT_H \backslash \mathtt{Lf} (\cT_H))= 0 $}    \big\} \; ,
\end{equation}
where we recall from (\ref{codef}) the definition of the real tree $\cT_H$ coded by $H$, where we recall from 
(\ref{brleafset}) that $\mathtt{Lf} (\cT_H)$ stands for the set of leaves of $\cT_H$ and where we 
recall from (\ref{masmea}) that $\bbm_H$ stands for the mass measure of $\cT_H$. Then, we set 
\begin{equation}
\label{defcH}
\cH= \big\{ B \cap \mathtt{Exc} \, ; \, \textrm{$B$ Borel subset of $\bC(\bbR_+, \bbR_+)$}\big\} \; .
\end{equation}
that is the trace sigma field on $\mathtt{Exc}$  of the Borel sigma field of $\bC(\bbR_+, \bbR_+)$. \cq 
\end{defi}
\begin{rem}
\label{lateruse1}
Let $H \! \in \! \mathtt{Exc}$ and  
let $s_0, s_1 \! \in \! (0, \zeta_H)$ be such that $s_0 \! < \! s_1$ and $d_H(s_0, s_1)\! = \! 0$. 
Then, we easily check that $ H^{[s_0]}_{\cdot \,  \wedge (s_1-s_0)} \! \! \in  \mathtt{Exc}$.  \cq 
\end{rem}
\begin{rem}
\label{lateruse}
Recall from (\ref{conttreex}) and from (\ref{conttree}) that $\bP^x$ and $\bN$ are supported by $\mathtt{Exc}$.   \cq 
\end{rem}

\begin{defi}
\label{defME} 
We introduce the following  subset of $\bbR_+ \! \times \! \bC(\bbR_+, \bbR_+)^2$: 
\begin{equation}
\label{defEspa}
E := \bbR_+ \! \times \! \big( \mathtt{Exc} \! \times \! (\mathtt{Exc}\!  \cup \! \{ \mathbf{0} \}  )\cup  (\mathtt{Exc} \! \cup\!  \{ \mathbf{0} \}  ) \! \times \!  \mathtt{Exc}  \big) \;  
\end{equation}
and we denote by $\ccM_{{\rm pt}}  (E)$ the set of point measures 
$$M (da\, d\overleftarrow{H} \, d \overrightarrow{H})\! =\!  \sum_{a\in \cJ} \delta_{(a, \overleftarrow{H}^a, \overrightarrow{H}^a )} $$
on $E$ that satisfy the following conditions: 
\begin{align}
\label{hypME}
 \exists\,  r\!  \in \! \bbR_+  \;  \textrm{such that the closure} & \; \textrm{of the countable set}  \;  \cJ \; \textrm{is $[0, r]$ and} \nonumber  \\
\forall \varepsilon , \eta \!  \in \! (0, \infty) & , \quad \#   \big\{ a\!  \in\!  \cJ:  \;  \Gamma (\overleftarrow{H}^a) \! \vee \!  \Gamma (\overrightarrow{H}^a) \! >\! \eta \quad \textrm{or} \quad  \zeta_{\overleftarrow{H}^a} \!  \vee \! \zeta_{\overrightarrow{H}^a} \! > \! \varepsilon   \big\} < \infty \; . 
\end{align}
We then equip $\ccM_{{\rm pt}}  (E)$ with the sigma field $\cG$ generated by the applications $M\! \in \! \ccM_{{\rm pt}} (E) \mapsto M(A)$, where $A$ ranges among the Borel subsets of $\bbR_+ \! \times \! \bC(\bbR_+, \bbR_+)^2$. \cq 
\end{defi}

The following 
lemma, whose proof is postponed in Appendix, asserts that $H$ can be recovered in a measurable way from the spinal decomposition $\cM_{0, t} (H)$, as defined in (\ref{spindec0t}). 
\begin{lem} 
\label{mesdec2} Recall from above the definition of the measurable spaces $(\mathtt{Exc}, \cH)$ and $(\ccM_{{\rm pt}} (E), \cG)$. Then, the following holds true. 
\begin{itemize}
\item[(i)] For all $t \! \in \! (0, \infty)$, we set $\{ \zeta \! >\! t\} \! :=\!  \{ H \! \in \! \mathtt{Exc}: \zeta_H \! >\! t \}$. Then, $\{ \zeta \! >\! t\} \! \in\! \cH $ and 
$$ H \in \{ \zeta \! >\! t\} \longmapsto \cM_{0, t} (H) \in \ccM_{{\rm pt}} (E) \quad \textrm{is measurable.} $$
\item[(ii)] There exists a measurable function $\Phi\! : \! \ccM_{{\rm pt}} (E)\!  \rightarrow\!  \bbR_+ \!\! \times \! \mathtt{Exc}$ such that 
$$ \forall H \in \mathtt{Exc}, \; \forall t \in (0, \zeta_H) , \quad 
\Phi( \cM_{0, t} (H))= (t, H) \; .$$
\end{itemize}
\end{lem}

\noi
\textbf{Proof.} See Appendix \ref{pfseclem}. \cqfd

\begin{figure}[htp]
    \centering
    \includegraphics[scale=.7]{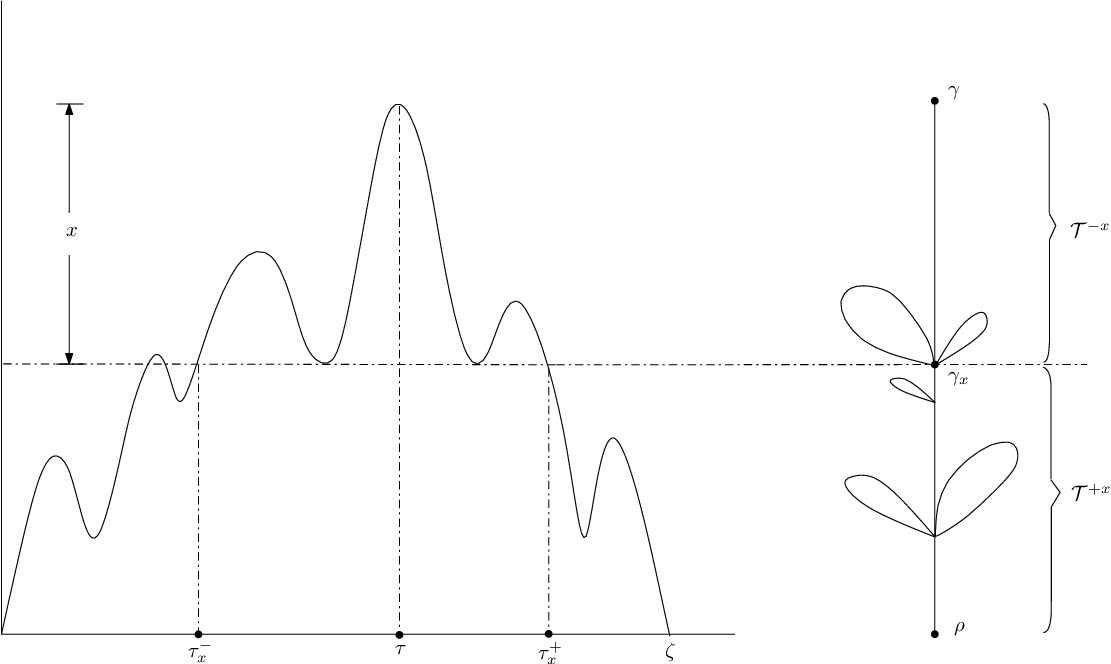}
    \caption{\label{fig: decomp_ht} {\small \textit{\hspace{-3mm}the left hand side figure illustrates the decomposition of $H$ into $H^{\ominus x}$ and $H^{\oplus x}$; the right hand side figure represents this decomposition in terms of the tree coded by $H$. }}}
\end{figure}

\paragraph{Decomposition according to the total height.} Let us fix $H \! \in \! \mathtt{Exc}$. Recall from (\ref{heigdiam}) the definition of 
$\Gamma (H)$, the height of $H$. We introduce the first time that realises the total height: 
\begin{equation}
\label{firsttau}
\tau (H)= \inf \{ t \! \in \! \bbR_+: H_t = \Gamma (H) \} \; .
\end{equation}
For all $x \! \in \! (0, \Gamma (H))$ we also introduce the following times: 
\begin{equation}
\label{tauplms}
\tau^-_x(H):= \sup \big\{ t \! <\! \tau (H) : H_t \! < \! \Gamma (H)\!-\! x  \big\}\quad \textrm{and} \quad 
\tau^+_x(H):= \inf \big\{ t \! >\! \tau (H) : H_t \! < \! \Gamma (H)\!-\! x  \big\}.
\end{equation}
Recall from (\ref{rerootH}) the definition of $H^{[s]}$. We then set 
\begin{equation}
\label{Hplmsx}
\forall \,  t \! \in \! \bbR_+, \qquad H^{\ominus x}_t = H^{[\tau_x^-]}_{t \,  \wedge (\tau^+_x -\tau^-_x)} \quad \textrm{and} \quad 
H^{\oplus x}_t = H^{[\tau_x^+]}_{t \, \wedge (\zeta -(\tau^+_x -\tau^-_x))}
\end{equation} 
where we denote $\tau^-_x\!: =\!  \tau^-_x (H)$, $\tau^+_x\! :=\!  \tau^+_x (H)$ and $\zeta\! := \! \zeta_H$ to simplify notation. See Figure \ref{fig: decomp_ht}.  

\medskip

Let us interpret $H^{\ominus x}$ and $H^{\oplus x}$ in terms of $\cT_H$. To that end, we recall that $p_H\! :\!  [0, \zeta] \! \rightarrow \! \cT_H$ stands for the canonical projection and we set 
$\gamma\! : = \! p_H (\tau(H))$. We first note that $d_H( \tau^-_x , \tau^+_x )\! = \! 0$. Then we set $\gamma (x)\! : =\! p_H(\tau^-_x)\! =\! p_H(\tau^+_x)$ that is the unique point of $\lgeo \rho, \gamma \rgeo$ such that $x\! = \! d(\gamma, \gamma (x))$ and thus, $d(\rho, \gamma (x))\! = \! \Gamma (H)\! -\! x$.   
We denote by $\cT^{o}$ the connected component of $\cT_H \backslash \{ \gamma (x)\}$ that contains the root $\rho$ and we set 
$$ \cT^{- x}= \cT_H \backslash \cT^{o} \quad \textrm{and} \quad \cT^{+x} = \{ \gamma (x)\} \cup \cT^{o} \; .$$
Thus $(\cT^{-x}\! ,  d, \gamma (x))$ is coded by  $H^{\ominus x}$ and  $(\cT^{+x}\! , d, \gamma (x))$ is coded by  $H^{\oplus x}$. See Figure \ref{fig: decomp_ht}.

\medskip

Recall from (\ref{rerootH}) the spinal decomposition of $H$ at a time $t$. We shall use the following notation: 
$$ \cM_{0, \tau(H)} (H)= \sum_{a \in \cJ_{0, \tau (H)}} \delta_{ (a, \overleftarrow{H}^a, \overrightarrow{H}^a) } \; . $$
This is a measure on $[0, \Gamma (H)] \! \times \! \mathtt{Exc}$ provides the spinal decomposition along the geodesic realising the total height. 
Let us first make the following remark. 

\begin{rem}
\label{simplept}
Let  $x \! \in \! (0, \Gamma (H))$ and recall the notation $\gamma (x)\! =\!  p_H(\tau_x^-(H)) \! =\!  p_H (\tau^+_x (H))$. 
Observe that if $x \notin \cJ_{0, \tau (H)}$, then $H_t \! >\!  \Gamma (H) \!-\! x $, for all $t \! \in \! (\tau^-_x (H), \tau^+_x (H))$ and thus, $\tau^-_x (H)$, $\tau^+_x (H)$ are the only time $t \! \in \! [0, \zeta_{H}]$ such that $p_H(t)\! =\! \gamma (x)$, which implies that $\gamma (x)$ is not a branching point of $\cT_H$: since it is not a leaf, it has to be a simple point of $\cT_H$. \cq
\end{rem}

For all  $x \! \in \! (0, \Gamma (H))$, we next introduce the following restriction of 
$ \cM_{0, \tau(H)} (H)$: 
\begin{equation}
\label{cMxplms}
\cM^{-x}_{0, \tau(H)} (H)=  \!\!\!\!\!\! \sum_{a \in \cJ_{0, \tau (H)} \cap [0, x]} \!\!\!\!\!\! 
\delta_{ (a, \overleftarrow{H}^a,  \overrightarrow{H}^a) }  \quad \textrm{and} \quad \cM^{+x}_{0, \tau(H)} (H)  = \!\!\!\!\!\! \sum_{a \in \cJ_{0, \tau (H)} \cap ( x, \Gamma (H)]}  \!\!\!\!\!\! 
\delta_{ (a, \overleftarrow{H}^a,  \overrightarrow{H}^a) } \; , 
\end{equation}
so that $\cM_{0, \tau(H)} (H)\! =\!  \cM^{-x}_{0, \tau(H)} (H)+\cM^{+x}_{0, \tau(H)} (H)$. Observe that 
\begin{equation}
\label{nesthtdec}
\tau (H)= \tau_x^- (H)+ \tau (H^{\ominus x} )  \quad \textrm{and} \quad \cM_{0,  \tau (H^{\ominus x} )  } (H^{\ominus x}) = \cM^{-x}_{0, \tau(H)} (H) \; .
\end{equation}

 For all $H^\prime\! \in \! \mathtt{Exc}$, we next denote by $\Lambda (H^\prime) \! := \! (H^\prime_{(\zeta_{H^\prime}-t)_+})_{t\geq 0}$ the function that reverses $H^\prime$ at its lifetime. We easily check that $\Lambda\! : \! \mathtt{Exc} \! \rightarrow \! \mathtt{Exc}$ is measurable; with a slight abuse of notation, we also set: 
$$ \Lambda \big( \cM_{0, \tau (H)}^{+x} (H) \big) =\!\!\! \!\!\!  \sum_{a \in \cJ_{0, \tau (H)} \cap ( x, \Gamma (H)]}  \!\!\!\!\!\! 
\delta_{ (\Gamma (H) -a \, ,\,  \Lambda (\overrightarrow{H}^a)\, ,\,    \Lambda (\overleftarrow{H}^a)) } \; .$$
It is easy to check first 
that  $\Lambda \big( \cM_{0, \tau (H)}^{+x} (H)\big) $ is a measurable function of $\cM_{0, \tau (H)}^{+x} (H)$ and next that 
\begin{equation}
\label{antinest} 
\cM_{0 , \zeta_H- \tau^+_x (H)} (H^{\oplus x} ) =  \Lambda \big( \cM_{0, \tau (H)}^{+x} (H) \big) \; .
\end{equation}
This combined with (\ref{nesthtdec}) and Lemma \ref{mesdec2} immediately implies the following lemma. 
\begin{lem}
\label{reconplms} There 
are two measurable functions $\Phi, \Phi^\prime: \ccM_{{\rm pt}}(E) \! \rightarrow \bbR_+ \! \times \! \mathtt{Exc}$ such that 
\begin{align*}
 \forall H \! \in \!  \mathtt{Exc}, \; \forall x \! \in \! (0, \Gamma (H)),  & \quad 
 \Phi \big(\cM^{-x}_{0, \tau(H)} (H) \big)= \big( \tau(H)\! -\!  \tau^-_x (H)\, , \,  H^{\ominus x} \big) \\
 \textrm{and} &  \quad \Phi^\prime\!  \big(\cM^{+x}_{0, \tau(H)} (H) \big)= 
\big( \zeta_H\! -\!  \tau^+_x (H)\, , \,  H^{\oplus x} \big)\; , 
\end{align*}
where $\tau(H)$ is defined by (\ref{firsttau}), $\tau^-_x (H)$ and $\tau^+_x (H)$ by (\ref{tauplms}), $H^{\ominus x}$ and $H^{\oplus x}$ by (\ref{Hplmsx}) and  $\cM^{-x}_{0, \tau(H)} (H)$ and 
$\cM^{+x}_{0, \tau(H)} (H)$ by (\ref{cMxplms}). 
\end{lem}

\subsection{Proofs of Theorem \ref{diamlaw} and of Theorem \ref{thm: decomp}.}
\label{pfth1th2sec}
As already mentioned, Abraham \& Delmas in \cite{AbDe09} make sense of the conditioned law $\bN (\, \cdot \, | \, \Gamma \! = \! r)$: namely they prove that $\bN (\, \cdot \, | \, \Gamma \! = \! r)$-a.s.~$\Gamma \! = \! r$, that  $r\mapsto \bN (\, \cdot \, | \, \Gamma \! = \! r)$ is weakly continuous on $\bC(\bbR_+, \bbR_+)$ and that (\ref{heightcondef}) holds true. Recall from (\ref{htcondsim}) and (\ref{renota}) the short-hand notations  
\begin{equation}
\label{rerenota}
\forall r, b, y \! \in \! (0, \infty), \quad \bN^{\Gamma}_{r}\! = \! \bN (\, \cdot \, | \, \Gamma \! = \! r), \quad  \bN_b \! = \! \bN \big(\, \cdot \, \cap\,  \{ \Gamma \leq b\} \big) \quad \textrm{and} \quad \bP^y_b\! =\!  \bP^y \big( \, \cdot \, \cap \,  \{ \Gamma \leq b\} \big), 
\end{equation}
where we recall from (\ref{xLevfo}) the notation $\bP^y$.  
Also recall from (\ref{taudef}) that $\bN^{\Gamma}_r$-a.s.~there exists a unique $\tau \! \in \! [0, \zeta]$ such that $H_{\tau} = \Gamma$. Recall from (\ref{spindec0t}) that $\cM_{0, \tau} (H)$ gives the excursions coding the subtrees grafted on $\lgeo \rho, p(\tau)\rgeo$ listed according to their 
distance of their grafting point from $p(\tau)$ (here $p\! : \! [0, \zeta] \! \rightarrow \! \cT$ stands for the canonical projection). In the following lemma, we recall from Abraham \& Delmas \cite{AbDe09} 
the following Poisson decomposition of $H$ under $\bN^\Gamma_r$ at its maximum, which extends Williams' decomposition that corresponds to the Brownian case.  
\begin{lem}
[Abraham \& Delmas \cite{AbDe09}]
\label{recAbDe} Let $\Psi$ be a branching mechanism of the form (\ref{eq: def-Psi}) that satisfies (\ref{eq: hyp}). We keep the previous notation. Let $r \! \in \! (0, \infty)$. Then, under $\bN^{\Gamma}_r$, 
\begin{equation}
\label{notaMtau}
\cM_{0, \tau} (da \, d\overleftarrow{H}\, d\overrightarrow{H})= \sum_{j\in \cJ_{0, \tau}} \delta_{(a, \overleftarrow{H}^a, \overrightarrow{H}^a)} 
\end{equation}
is Poisson point process on $[0, r] \! \times \! \bC(\bbR_+, \bbR_+)^2$ whose intensity is 
\begin{align}
\label{decomGr}
\mathbf{n}_r (da \,d\overleftarrow{H}\, d\overrightarrow{H}) := \beta \un_{[0, r]} (a)  da  & \, \Big( \delta_{\mathbf{0}} (d\overleftarrow{H}) \bN_{a} ( d\overrightarrow{H}) + \bN_{a}  ( d\overleftarrow{H}) \delta_{\mathbf{0}} (d\overrightarrow{H})  \Big) \nonumber \\
+ &\; \;  \un_{[0, r]} (a) da \int_{(0, \infty)}\!\!\!\!\!\!\!\! \!\! \pi (dz)\!\!  \int_0^z \!\! \!\! dx \;  \bP^x_{a } \big(d\overleftarrow{H} )\,  \bP^{z-x}_{a } \big(  d\overrightarrow{H} ) , 
\end{align}
where $\beta$ and $\pi$ are defined in (\ref{eq: def-Psi}) and where 
$\mathbf{0}$ stands for the null function. 
\end{lem}

  We first discuss several consequences of Lemma \ref{recAbDe}. To that end, we set 
$$ \bnu_{r, a} (d\overleftarrow{H}\, d\overrightarrow{H})= \beta \delta_{\mathbf{0}} (d\overleftarrow{H}) \bN_{a} ( d\overrightarrow{H}) + \beta\bN_{a}  ( d\overleftarrow{H}) \delta_{\mathbf{0}} (d\overrightarrow{H}) + \int_{(0, \infty)}\!\!\!\!\!\!\!\! \!\! \pi (dz)\!\!  \int_0^z \!\! \!\! dx \; \bP^x_{a } \big(d\overleftarrow{H} )\,  \bP^{z-x}_{a } \big(  d\overrightarrow{H} ) ,  $$ 
so that $\mathbf{n}_r (da \,d\overleftarrow{H}\, d\overrightarrow{H}) = \un_{[0, r]} (a) da  \, \bnu_{r, a} (d\overleftarrow{H}\, d\overrightarrow{H})$. Denote by $\langle \bnu_{r,a} \rangle$ the total mass of $\bnu_{r, a}$. We claim that 
$\langle \bnu_{r,a} \rangle\! = \! \infty$. 
\textit{Indeed}, first recall that $\bN$ is an infinite measure. Since $\bN ( \Gamma \! > \! a ) \! < \! \infty$ (by (\ref{lawGamma})), 
$\bN_a $ is also an infinite measure. Thus, if $\beta \! >\! 0$,  
$\langle \bnu_{r,a} \rangle\! = \! \infty $. Suppose now that $\beta\! = \! 0$. Then by (\ref{tauPx}), we get 
$ \langle \bnu_{r,a} \rangle\! = \!  \! \int_{(0, \infty)}\!  \pi (dz) z e^{-z v(a)}\! = \! \infty $, since 
$\int_{(0, \infty)} \! z \, \pi(dz) \! = \! \infty$, by (\ref{infvarhyp}). 

Therefore, standard results on Poisson point measures entail 
that $\bN^{\Gamma}_r$-a.s.~the closure of $\cJ_{0, \tau}$ is $[0, r]$. This point combined with the fact that $H$ is $\bN^{\Gamma}_r$-a.s.~continuous with compact support implies that 
$\bN^\Gamma_r$-a.s.~$\cM_{0, \tau} \! \in \! \ccM_{{\rm pt}} (E)$, 
where the set of point measures $\ccM_{{\rm pt}} (E)$ is defined in Definition \ref{defME}. 

Recall from (\ref{defExc}) the definition 
of $\mathtt{Exc}$ and recall from (\ref{conttreex}) and from (\ref{conttree}) that $\bP^x$ and $\bN$ are supported by $\mathtt{Exc}$. We easily derive from (\ref{heightcondef}) 
that $\bN^{\Gamma}_r$-a.s.~$H\ino \mathtt{Exc}$. 

Next recall that $\Lambda \! :\!  \mathtt{Exc}\! \rightarrow \!  \mathtt{Exc}$, its the functional that reverses excursions at their lifetime: namely 
for all $H\! \in \! \mathtt{Exc}$, we denote by $\Lambda (H)= (H_{(\zeta_{H}-t)_+})_{t\geq 0}$. Then, Corollary 3.1.6 \cite{Duquesne02} asserts that $H$ and $\Lambda (H)$ have the same distribution under $\bN$. This also implies that 
$H$ and $\Lambda (H)$ have the same law under $\bP^x$ and by (\ref{heightcondef}) we easily see that $H$ and $\Lambda (H)$ have the same law under $\bN^{\Gamma}_r$. 

  We thus have proved the following. 
 \begin{equation}
\label{revHr}   
\textrm{$H$ and $\Lambda (H)$ have the same law under $\bN^\Gamma_r$ and $\bN^\Gamma_r$-a.s.} \quad H \in \mathtt{Exc} \; \;  \textrm{and} \quad \cM_{0, \tau} \! \in \! \ccM_{{\rm pt}} (E) \; .
\end{equation} 

Recall from (\ref{firsttau}) the definition of $\tau(H)$, from (\ref{tauplms}) that of 
$\tau^-_x (H)$ and $\tau^+_x (H)$, from (\ref{Hplmsx}) that of $H^{\ominus x}$ and $H^{\oplus x}$, and from 
(\ref{cMxplms}) that of $\cM^{-x}_{0, \tau(H)} (H)$ and 
$\cM^{+x}_{0, \tau(H)} (H)$. To simplify notation we simply write $\tau$, $\tau^-_x$, $\tau^+_x$, 
$\cM^{-x}_{0, \tau}$ and $\cM^{+x}_{0, \tau}$. 
We then prove the following lemma. 
\begin{lem} 
\label{nestlem} We keep the same assumptions as in Lemma \ref{recAbDe} and the notation therein. 
Let $x \! \in \! (0, r)$. Then, the following holds true. 
\begin{itemize} 
\item[(i)] Under $\bN^{\Gamma}_r$, $\cM^{-x}_{0, \tau}$ and $\cM^{+x}_{0, \tau}$ are independent Poisson point measures. 
\item[(ii)] $\bN^{\Gamma}_r$-a.s.~$x\notin \cJ_{0, \tau}$. 
\item[(iii)] $\cM^{-x}_{0, \tau}$ under $\bN^{\Gamma}_r$ has the same law as $\cM_{0, \tau}$ under 
$\bN^{\Gamma}_x$. Thus the law of $H^{\ominus x}$ under $ \bN^{\Gamma}_r$ is $\bN^{\Gamma}_x$. 
\end{itemize}
\end{lem}
\noi
\textbf{Proof.} Point $(i)$ is a consequence of Lemma \ref{recAbDe} and of 
basic results on Poisson point measures. Moreover, $\cM^{-x}_{0, \tau}$ 
under $\bN^{\Gamma}_r$ has intensity $\un_{[0, x]} (a) da  \bnu_{r, a} (d\overleftarrow{H}\, d\overrightarrow{H})$ which is equal to $\mathbf{n}_x$. This implies that
$\cM^{-x}_{0, \tau}$ under $\bN^{\Gamma}_r$ has the same law as $\cM_{0, \tau}$ under 
$\bN^{\Gamma}_x$. By Lemma \ref{mesdec2} and Lemma \ref{reconplms}, it implies that 
$$ ( \tau \! -\! \tau^-_x, H^{\ominus x} ) = \Phi \big( \cM_{0, \tau}^{-x}\big) \quad \textrm{under} \quad \bN^{\Gamma}_r \quad \overset{\textrm{law}}{=} \quad (\tau, H)= \Phi \big(\cM_{0, \tau} \big)  \quad \textrm{under} \quad 
\bN^{\Gamma}_x\; , $$
which entails $(iii)$. Since the intensity measure $\mathbf{n}_r (da\, d\overleftarrow{H}\, d\overrightarrow{H})$ is diffuse in the variable $a$, standard results on Poisson point measures entail $(ii)$. \cqfd 

\paragraph{Proof of Theorem \ref{diamlaw} $(i)$.} We keep the previous notation and we set 
\begin{equation}
\label{notadelt}
\forall b \in (0, \infty), \; \forall  \overleftarrow{H}, \overrightarrow{H} \in \mathtt{Exc}, \qquad \Delta_{b, \overleftarrow{H}, \overrightarrow{H} }= b+ \Gamma (\overleftarrow{H}) \! \vee \! \Gamma (\overrightarrow{H}) \; .
\end{equation}
Recall from (\ref{lawGamma}) and (\ref{tauPx}) that the distributions of $\Gamma$ under $\bN$ and under $\bP^x$ are diffuse. Thus, for all $a\! \in (0, \infty)$, the distributions of $\Gamma$ under $\bN_a$ and under $\bP^x_a$ are also diffuse. Recall the notation (\ref{notaMtau}) for $\cM_{0, \tau}$. 
Then, Lemma \ref{recAbDe} combined with Lemma \ref{lem: dm} implies that $\bN^{\Gamma}_r$-a.s.~there exists a unique $Y \! \in \! (0, r) \cap \cJ_{0, \tau}$ such that 
\begin{equation}
\label{laYdef}
D= Y+ \Gamma (\overleftarrow{H}^Y) \! \vee \! \Gamma (\overrightarrow{H}^Y)= \Delta_{Y, \overleftarrow{H}^Y, \overrightarrow{H}^Y }\; >  \!\! \!\!\! \sup_{ \; \; \; a \in \cJ_{0, \tau} \backslash \{ Y\}} \!\! \!\!\!  \Delta_{a, \overleftarrow{H}^a, \overrightarrow{H}^a } \; .
\end{equation} 
Then either $ \Gamma (\overleftarrow{H}^Y) \! <\! \Gamma (\overrightarrow{H}^Y) $ or 
$ \Gamma (\overleftarrow{H}^Y) \! >\! \Gamma (\overrightarrow{H}^Y) $. Let us us consider these two cases.

\smallskip

\noi
$\bullet$ If $ \Gamma (\overleftarrow{H}^Y) \! <\! \Gamma (\overrightarrow{H}^Y) $ then by (\ref{taudef}) and (\ref{tauPx}) there exists a unique point 
$t_*$ such that 
$\overrightarrow{H}^Y_{\! t_*}\! \! =\!  \Gamma (\overrightarrow{H}^Y)$. 
This entails Theorem \ref{diamlaw} $(i)$ in this case under $\bN^{\Gamma}_r$ and we have 
$ \tau_0\! = \! \tau$ and 
$$ \tau_1= \tau + t_*+\!\!\!\!\!\!\!  \sum_{\; \; a\in \cJ_{0, \tau} \cap [0, Y)} \!\!\!\!\!\!\! \zeta_{\overrightarrow{H}^a} \; .$$

\noi
$\bullet$ If $ \Gamma (\overleftarrow{H}^Y) \! >\! \Gamma (\overrightarrow{H}^Y) $ then by (\ref{taudef}) and (\ref{tauPx}) there exists a unique point 
$t_*$ such that 
$\overleftarrow{H}^Y_{\! t_*}\! \! =\!  \Gamma (\overleftarrow{H}^Y)$. 
This entails Theorem \ref{diamlaw} $(i)$ in this case under $\bN^{\Gamma}_r$ and we have 
$ \tau_1\! = \! \tau$ and 
$$   \tau_0= t_* + \!\!\!\!\!\!\!  \sum_{\; \; a\in \cJ_{0, \tau} \cap (Y, r]}  \!\!\!\!\!\!\! \zeta_{\overleftarrow{H}^a}  \; .$$
Theorem \ref{diamlaw} $(i)$ is then proved under $\bN_{r}^\Gamma$, for all $r \! \in \! (0, \infty)$, which implies Theorem \ref{diamlaw} $(i)$ (under $\bN$) by (\ref{heightcondef}).  \cq

\paragraph{Proof Theorem \ref{diamlaw} $(ii)$.} Recall from (\ref{cMxplms}) the notation 
$\cM^{-x}_{0, \tau}$ and 
$\cM^{+x}_{0, \tau}$. We shall use the following lemma. 
\begin{lem}
\label{deviss1} We keep the same assumptions as in Lemma \ref{recAbDe} and the notation therein. 
Recall from Definition \ref{defME} the notation $ \ccM_{{\rm pt}} (E)$. Then, for all $r \! \in \! (0, \infty)$ and for all measurable functions $G_1, G_2 \! : \! \ccM_{{\rm pt}} (E) \! \rightarrow \! \bbR_+$, 
$$\bN^\Gamma_r \Big[ \un_{\{ \tau= \tau_0\}} G_1\big( \cM_{0, \tau}^{- \frac{_1}{^2} \! D}\big) G_2\big( \cM_{0, \tau}^{+ \frac{_1}{^2} \! D}\big)\Big]=  \bN^\Gamma_r \Big[ \un_{\{ \tau= \tau_0\}} \bN^{\Gamma}_{\! \frac{_1}{^2}\!   D} [ G_1 (\cM_{0, \tau})]  G_2\big(\cM_{0, \tau}^{+ \frac{_1}{^2} \! D} \big)\Big] \; , $$
with a similar statements where $\tau_0$ is replaced by $\tau_1$. Moreover, 
by (\ref{heightcondef}) a similar statement holds true under $\bN$.  
\end{lem}
Before proving this lemma, we first complete the proof of Theorem \ref{diamlaw}.  Recall from the notation (\ref{notaMtau}) and from (\ref{cMxplms}) that 
$$ \cM_{0, \tau} = \sum_{j\in \cJ_{0, \tau}} \delta_{(a, \overleftarrow{H}^a, \overrightarrow{H}^a)} \quad \textrm{and} \quad 
\cM_{0, \tau}^{- \frac{_1}{2} \! D}= \sum_{j\in \cJ_{0, \tau}\cap [0, \frac{_1}{2} \! D] } \delta_{(a, \overleftarrow{H}^a, \overrightarrow{H}^a)}\; .$$
We the next see the event $\{  \frac{_1}{^2} D \! \in \! \cJ_{0, \tau}\}$ the the event that $\cM_{0, \tau}^{- \frac{_1}{^2} \! D}$ has an atom "at" $\frac{_1}{^2} \! D$. 
By Lemma \ref{deviss1} with $G_2\equiv 1$ we then get  
$$  \bN \big( \frac{_{_1}}{^{^2}} D \! \in \! \cJ_{0, \tau} \big) \! = \! \int_0^\infty \!\! \!\!\!\!\! N(D\ino dr) \,  \bN^{\Gamma}_{ \! \frac{_{1}}{^{2}} \!  r} \big( \frac{_{_1}}{^{^2}} r \ino \cJ_{0, \tau}\big) = 0 $$
because for any $b \! \in \! (0, \infty)$, Lemma \ref{recAbDe} asserts that under 
$\bN^{\Gamma}_b$, $\cM_{0, \tau}$ is a Poisson point measure 
with intensity $\mathbf{n}_b$, which implies that  $\bN^{\Gamma}_b$-a.s.~$b\! \notin\! \cJ_{0, \tau}$. We next use Remark \ref{simplept} with $x\! =\!  \frac{_1}{^2} D$ that asserts that 
\begin{equation}
\label{notatasse}
 \tau^-_{{\rm mid}}\! :=\!  \tau_{\frac{1}{2} D}^- \quad \textrm{and} \quad   
\tau^+_{{\rm mid}}\! :=\!  \tau_{\frac{1}{2} D}^+
\end{equation}
are the only times $t\! \in \! [0, \zeta]$, such that $d( p(\tau_1), p(t))\! = \! \frac{_1}{^2} D$,  
which completes the proof of Theorem \ref{diamlaw} $(ii)$. \cq

\paragraph{Proof Theorem \ref{diamlaw} $(iii)$.} Let $r, y \! \in \! (0, \infty)$ be such that $\frac{_1}{^2}y\! < \! r \! < \!  y$. 
We first work under $\bN^{\Gamma}_r$. Recall from (\ref{notaMtau}) the notation for $\cM_{0, \tau}$ and recall notation (\ref{notadelt}). 
Then (\ref{laYdef}) combined with Lemma \ref{recAbDe} that asserts that under $\bN^{\Gamma}_r$, $\cM_{0, \tau}$ is a Poisson point measure with intensity $\mathbf{n}_r$, we get 
\begin{equation}
\label{NrGaDy}
\bN^\Gamma_r \big( D \! \leq \! y\big)\! = \! \bN^\Gamma_r \big( \sup\!  \big\{ \Delta_{a, \overleftarrow{H}^a ,  \overrightarrow{H}^a} \, ; \, a \in \cJ_{0, \tau} \big\} 
\leq  y\big) \! =\!  \exp \Big( \! - \!\! \int \! \mathbf{n}_r (da\, d\overleftarrow{H} \,  d\overrightarrow{H})  \un_{\{ \Delta_{a, \overleftarrow{H}^a ,  \overrightarrow{H}^a} >y \}} \Big),  
\end{equation}
where $\mathbf{n}_r$ is given by (\ref{decomGr}). Recall from (\ref{lawGamma}) that $\bN(\Gamma \! >\! t)\! = \! v(t)$ and from (\ref{tauPx}) that $\bP^x (\Gamma \! \leq \! t)\! =\! e^{-x v(t)}$.  Thus, 
\begin{align*}
\int \! \mathbf{n}_r (da\, d\overleftarrow{H} \,  d\overrightarrow{H})  \un_{\{ \Delta_{a, \overleftarrow{H}^a ,  \overrightarrow{H}^a} >y \}}= 2 \beta\int_0^r \!\! \! da\, & \bN(y \!-\! a\! <\! \Gamma \! \leq  \! a  ) \\
+ \int_0^r \!\! \!  da\!  \int_{(0, \infty)}\!\!\!\!\!\!\!\! \!\! \pi (dz)\!\!  \int_0^z \!\! \!\! dx \!\! \int\!  \bP^x_{a } \! \big(d\overleftarrow{H} ) & \! \int\! \bP^{z-x}_{a } \! \big(  d\overrightarrow{H} ) \big( 1-
\un_{\{ \Gamma(\overleftarrow{H})\leq y-a \}  } \un_{\{ \Gamma(\overrightarrow{H})\leq y-a \}  }\big). 
\end{align*}
If $a\! <\! \frac{_1}{^2}y$, then $\bN(y \!-\! a\! <\! \Gamma \! \leq  \! a  )\!= \! 0$ and if 
 $a\! >\! \frac{_1}{^2}y$, then  
$\bN(y \!-\! a\! <\! \Gamma \! \leq  \! a  )\! = \! v(y\!-\! a) \! -\! v(a)$. Recall that the total mass 
of $\bP^x_a$ is $\bP^x (\Gamma \! \leq \! a)\! = \! \exp(-xv(a))$ and observe that 
$\bP^x_a (\Gamma \! \leq \! y\! - \! a)=\bP^x (\Gamma \! \leq \! a \! \wedge \! (y\! - \! a))= \exp (-xv(a \! \wedge \! (y\! - \! a))$. Thus   
$$\int \! \bP^x_{a } \! \big(d\overleftarrow{H} )\! \int \!   \bP^{z-x}_{a } \! \big(  d\overrightarrow{H} ) \big( 1-
\un_{\{ \Gamma(\overleftarrow{H})\leq y-a \}  } \un_{\{ \Gamma(\overrightarrow{H})\leq y-a \}  }\big)= e^{-zv(a)} -e^{-zv(a \wedge (y -  a))}, $$
which is null if  $a\! <\! \frac{_1}{^2}y$. Note that this expression does not depend on $x$. Consequently,  
\begin{eqnarray*}
\int \! \! \mathbf{n}_r (da\, d\overleftarrow{H} \,  d\overrightarrow{H})  \un_{\{ \Delta_{a, \overleftarrow{H}^a ,  \overrightarrow{H}^a} >y \}}\!\!\!\!\!\!  &=& \!\!\!\!\!\!
\int_{\frac{1}{2}y}^r \!\!\! \!\! da \, 2\beta \big(  v(y\!-\! a) \! -\! v(a)  \big) + \int_{\frac{1}{2}y}^r \!\!\!\!\!  da\!\!  \int_{(0, \infty)}\!\!\!\!\!\!\!\! \!\! \pi (dz)\, z \big( e^{-zv(a)} -e^{-zv(y -  a)} \big) \\
& =& \!\!\!\!\!\! \int_{\frac{1}{2}y}^r \!\!\! \!\! da \big(\Psi^\prime (v(y\!-\! a)) \! -\!  \Psi^\prime (v(a)) \big)\! = \!
\int^{\frac{1}{2}y}_{y-r}\!\!\!\! db \, \Psi^\prime (v(b)) \! -\!\! \int_{\frac{1}{2}y}^r \!\!\! \!\! db\,    \Psi^\prime (v(b))
\end{eqnarray*}
by (\ref{eq: def-Psi}). Recall that $v$ satisfies $\int_{v(b)}^\infty d\lambda / \Psi (\lambda) \! = \! b$.The change of 
variable $\lambda\! = \! v(b)$ entails 
\begin{eqnarray*}
\int  \mathbf{n}_r (da\, d\overleftarrow{H} \,  d\overrightarrow{H})  \un_{\{ \Delta_{a, \overleftarrow{H}^a ,  \overrightarrow{H}^a} >y \}}  &=&   \int_{v(\frac{1}{2}y)}^{v(y-r)} \!\!\!\!\! d\lambda \, 
\frac{\Psi^\prime (\lambda)}{\Psi (\lambda)}\;  -\int_{v(r)}^{v(\frac{1}{2}y)} \!\!\!\!\! d\lambda \, 
\frac{\Psi^\prime (\lambda)}{\Psi (\lambda)} \\
 &=& \log \frac{\Psi(v(y\! -\! r))}{\Psi(v(\frac{_1}{^2}y))} \; -\; \log  \frac{\Psi(v(\frac{_1}{^2}y))}{\Psi(v(r))} . 
 \end{eqnarray*}
By (\ref{NrGaDy}), we get 
\begin{equation}
\label{equation}
\forall \, r \! \in \! (0, \infty) ,\;  \forall \, y \! \in \! (r, 2r), \quad \bN^\Gamma_r \big( D \! \leq \! y\big)= 
\frac{\Psi(v(\frac{_1}{^2}y))^2 }{\Psi(v(r))\Psi(v(y\! -\! r))}\; .
\end{equation}
Now observe that $\bN^\Gamma_r \big( D \! > \! y\big)\! = \! 0$, if $y\! \geq \! 2r$ and that 
$\bN^\Gamma_r \big( D \! \geq \! y\big)\! = \! 1$, if $y \! \leq \! r$. Thus by (\ref{heightcondef}), 
\begin{eqnarray*}
\bN(D \! > \! y) \!\! &=& \!\!\! \int_0^\infty \!\!\!\! \bN(\Gamma \! \in \! dr) \, \bN^\Gamma_r(D \! > \! y)= 
\bN(\Gamma \! >\! y)+\! \int_{\frac{1}{2}y}^{y} \!\!\!\! dr \,\Psi (v(r)) \Big(1-\frac{\Psi(v(\frac{_1}{^2}y))^2 }{\Psi(v(r))\Psi(v(y\! -\! r))} \Big) \\
& =& v(\frac{_{_1}}{^{^2}} y) -\Psi(v(\frac{_{_1}}{^{^2}}y))^2  \! \int_{\frac{1}{2}y}^{y} \!\! 
\frac{dr}{\Psi(v(y\! -\! r))}=  v(\frac{_{_1}}{^{^2}} y) -\Psi(v(\frac{_{_1}}{^{^2}}y))^2  \! \int_{v(\frac{1}{2}y)}^\infty 
\!\! \frac{d\lambda}{\Psi (\lambda)^2} , 
 \end{eqnarray*}
where we use the change of variable $\lambda\! = \! v(y\!-\!r)$ in the last equality. This proves (\ref{diamrepa}) that 
easily entails (\ref{diamdens}), which completes the proof of Theorem \ref{diamlaw} $(iii)$. \cq

\paragraph{Proof of Lemma \ref{deviss1}.} To completes the proof of Theorem \ref{diamlaw}, it remains to prove 
Lemma \ref{deviss1} that is also the key argument to prove Theorem \ref{thm: decomp}. We first work under $\bN^\Gamma_r$. Recall the notation (\ref{notaMtau}) for $\cM_{0, \tau}$ and $\cJ_{0, \tau}$ and recall from 
(\ref{cMxplms}) the following definitions (with $x\! = \! \frac{_1}{^2} D$),
$$ \cM_{0, \tau} \! =\!\!\!\!  \sum_{j\in \cJ_{0, \tau}} \!\! \delta_{(a, \overleftarrow{H}^a, \overrightarrow{H}^a)} \; , \quad 
\cM_{0, \tau}^{- \frac{_1}{2} \! D}\!= \!\!\!\! \!\!\!\!  \sum_{j\in \cJ_{0, \tau}\cap [0, \frac{_1}{2} \! D] }\!\!\!\! \!\! \!\!  \delta_{(a, \overleftarrow{H}^a, \overrightarrow{H}^a)} \quad \textrm{and} \quad  \cM_{0, \tau}^{+ \frac{_1}{2} \! D}\!= \!\!\!\! \!\! \!\! 
\sum_{j\in \cJ_{0, \tau}\cap (\frac{_1}{2} \! D , r]}\!\!\!\! \!\!  \delta_{(a, \overleftarrow{H}^a, \overrightarrow{H}^a)} . $$
Recall from (\ref{laYdef}) the definition of the random variable $Y$: since $ \Gamma (\overleftarrow{H}^Y) \! \vee \! \Gamma (\overrightarrow{H}^Y)\! <\! Y$, 
we get $Y\! >\! \frac{_1}{^2}D$ and $(Y, \overleftarrow{H}^Y, \overrightarrow{H}^Y)$ is an atom of $\cM_{0, \tau}^{+ \frac{_1}{2} \! D}$. This argument, combined with (\ref{laYdef}) and the Palm formula for Poisson point measures, implies 
\begin{align}
\label{Palm1}
 \bN^\Gamma_r \!  \Big[ &\un_{\{ \tau= \tau_0\}}  F \big( Y, \overleftarrow{H}^Y, \overrightarrow{H}^Y \big) 
 G_1\big( \cM_{0, \tau}^{- \frac{_1}{^2} \! D}\big) 
G_2\big( \cM_{0, \tau}^{+ \frac{_1}{^2} \! D}\big)\Big]
= \nonumber \\
&  \int\!  \mathbf{n}_r (dy \,dH^\prime \, dH^{\prime \prime}) 
\un_{\{ \Gamma (H^{\prime  \prime})  >
\Gamma (H^{\prime}) \}} F( y,H^\prime ,H^{\prime \prime})  \\
\times & \; \bN^\Gamma_r \Big[ G_1 \big(  \cM_{0, \tau}^{- \frac{_1}{^2} \! \Delta_{y,H^{\prime }\!\! , 
H^{\prime \prime}\!}}  \big)
G_2\Big( \cM_{0, \tau}^{+ \frac{_1}{^2} \! \Delta_{y,H^{\prime }\!\! , 
H^{\prime \prime}\!}}\! + \delta_{(y,H^{\prime}, H^{\prime \prime}) } \Big) \un_{\big\{ \Delta_{y,H^{\prime }\!\! , 
H^{\prime \prime}\!} \, > \, \sup \{\Delta_{a,\overleftarrow{H}^a, \overrightarrow{H}^a} \, ; \, a\in \cJ_{0, \tau} 
\}    \big\}}  \Big]    \nonumber 
\end{align}
where we recall that $\tau_0\! = \! \tau$ iff $\Gamma (\overrightarrow{H}^Y) \! >\!
\Gamma (\overleftarrow{H}^Y)$. Then observe that $\mathbf{n}_r \! \otimes \! \bN^\Gamma_r$-a.e.~for all $a\!  \in \! 
\cJ_{0, \tau} \! \cap\!  [0, 
\frac{_{_1}}{^{^2}} \Delta_{y,H^{\prime }\!\! , 
H^{\prime \prime}\!} ]$, we have $\Delta_{a,\overleftarrow{H}^a, \overrightarrow{H}^a} \! < \! 2a  \leq \Delta_{y,H^{\prime }\!\! , 
H^{\prime \prime}\!}$. Thus, $\mathbf{n}_r \!  \otimes \! \bN^\Gamma_r$-a.e.
$$\un_{\big\{ \Delta_{y,H^{\prime }\!\! , 
H^{\prime \prime}\!} \, > \, \sup \{\Delta_{a,\overleftarrow{H}^a, \overrightarrow{H}^a} \, ; \, a\in \cJ_{0, \tau} 
\}    \big\}}=\un_{\big\{ \Delta_{y,H^{\prime }\!\! , 
H^{\prime \prime}\!} \, > \, \sup \{\Delta_{a,\overleftarrow{H}^a, \overrightarrow{H}^a} \, ; \, a\in \cJ_{0, \tau} \cap (\frac{_{_1}}{^{^2}}  \Delta_{y,H^{\prime }, H^{\prime \prime}},  r]
\}    \big\}} $$
that only depends on $y,H^{\prime}\!\!, H^{\prime \prime}\!$ and of $\cM_{0, \tau}^{+ \frac{_1}{^2} \! \Delta_{y,H^{\prime }\!\! , 
H^{\prime \prime}\!}}\!\!\! $. By (\ref{Palm1}) with $F\equiv 1$ and by Lemma \ref{nestlem} $(i)$ and $(iii)$ with $x\! = \!  \frac{_{_1}}{^{^2}} \! \Delta_{y,H^{\prime }\!\! , 
H^{\prime \prime}\!}$, we get 
\begin{align*}
\bN^\Gamma_r \!  \Big[ &\un_{\{ \tau= \tau_0\}} 
G_1\big( \cM_{0, \tau}^{- \frac{_1}{^2} \! D}\big) 
G_2\big( \cM_{0, \tau}^{+ \frac{_1}{^2} \! D}\big)\Big]
= \nonumber \\
&  \int\!  \mathbf{n}_r (dy \,dH^\prime \, dH^{\prime \prime}) 
\un_{\{ \Gamma (H^{\prime  \prime})  >
\Gamma (H^{\prime}) \}} \, \bN^\Gamma_{\frac{_1}{^2} \! \Delta_{y,H^{\prime }\!\! , 
H^{\prime \prime}\!}}  \big[ G_1 (\cM_{0, \tau} ) \big]  \\
\times & \; \bN^\Gamma_r \Big[  
G_2\Big( \cM_{0, \tau}^{+ \frac{_1}{^2} \! \Delta_{y,H^{\prime }\!\! , H^{\prime \prime}\!\!}}\! + \delta_{(y,H^{\prime}, H^{\prime \prime}) } \Big) \un_{\big\{ \Delta_{y,H^{\prime }\!\! , 
H^{\prime \prime}\!} \, > \, \sup \{\Delta_{a,\overleftarrow{H}^a, \overrightarrow{H}^a} \, ; \, a\in \cJ_{0, \tau} 
\}    \big\}}  \Big]  .  \nonumber \\
& =  \bN^{\Gamma}_r \Big[ \un_{\{ \tau= \tau_0\}} \bN^{\Gamma}_{\! \frac{_1}{^2}\!   D} [ G_1 (\cM_{0, \tau})]  G_2\big(\cM_{0, \tau}^{+ \frac{_1}{^2} \! D} \big)\Big] , 
\end{align*}
which completes the proof of Lemma \ref{deviss1} when $\tau\! = \! \tau_0$ under $\bN^\Gamma_r$. When $\tau\! = \! \tau_1$, the proof is quite similar. Then, (\ref{heightcondef}) immediately entails the same result under $\bN$.  \cq

\paragraph{Proof of Theorem \ref{thm: decomp} $(iii)$.} Recall from (\ref{Hplmsx}) the definition of $H^{\ominus x}$ and $H^{\oplus x}$. Then, 
Lemma \ref{deviss1} under $\bN$ and Lemma \ref{reconplms} imply that 
for all measurable functions 
$F_1, F_2 \! :\!  \bC(\bbR_+, \bbR_+) \! \rightarrow \! \bbR_+$, $f: \bbR_+ \rightarrow \bbR_+$, 
\begin{equation}
\label{quartcas}
\bN \Big[ \un_{\{ \tau= \tau_0\}} f(D) F_1\big( H^{\ominus \frac{1}{2} \! D}\big) F_2\big( H^{\oplus \frac{1}{2} \! D}\big)\Big]=  \bN \Big[ \un_{\{ \tau= \tau_0\}}f(D)  \bN^{\Gamma}_{\frac{1}{2}\!   D} [ F_1 (H)] 
 F_2\big( H^{\oplus \frac{1}{2}\! D}\big)\Big] \; , 
 \end{equation}
with a similar statement with $\tau\! = \! \tau_1$. To simplify notation, we next set 
\begin{equation}
\label{dgcxtswa}
 H^{\ominus}:= H^{\ominus \frac{1}{2} \! D} 
\quad \textrm{and} \quad H^{\oplus}:= H^{\oplus \frac{1}{2} \! D} \; .
\end{equation}
By adding (\ref{quartcas}) with the analogous equality with $\tau\! = \! \tau_1$, we get   
\begin{equation}
\label{demicas}
\bN \Big[ f(D) F_1\big( H^{\ominus }\big) F_2\big( H^{\oplus }\big)\Big]=  \bN \Big[ f(D)  \, \bN^{\Gamma}_{\! \frac{1}{2}\!   D} [ F_1 (H)] 
\,  F_2\big( H^{\oplus }\big)\Big] \; . 
 \end{equation}
Recall from (\ref{notatasse}) that $ \tau^-_{{\rm mid}}\! =\!  \tau_{\frac{1}{2} D}^-$ and $ 
\tau^+_{{\rm mid}}\! =\!  \tau_{\frac{1}{2} D}^+$; rewriting (\ref{Hplmsx}) with $x\! = \! \frac{1}{2} D$ yields
\begin{equation}
\label{idnotapff}
H^{\ominus}= H^{[\tau^-_{{\rm mid}}]}_{ \cdot \; \wedge \, (\tau^+_{{\rm mid}}-\tau^-_{{\rm mid}})}, \quad 
H^{\oplus}= H^{[\tau^+_{{\rm mid}}]}_{ \cdot \; \wedge \; (\zeta-(\tau^+_{{\rm mid}}-\tau^-_{{\rm mid}}))} \quad \textrm{and thus} \quad H^{[\tau^-_{{\rm mid}}]}= H^{\ominus} \oplus H^{\oplus}, 
\end{equation}
where we recall from (\ref{concadef}) that $H^\prime \! \oplus \!  H^{\prime \prime}$ stands for the concatenation of the functions 
$H^\prime$ and $H^{\prime \prime}$. 

\medskip

Let us briefly interpret $H^{\ominus}$ and $H^{\oplus}$ in terms of the tree $\cT$. 
To that end, first recall that $\gamma \! = \! p(\tau)$,  $\gamma_0 \! = \! p(\tau_0)$ and  $\gamma_1 \! = \! p(\tau_1)$, where $p\! :
\! [0, \zeta] \! \rightarrow \cT$ stands for the canonical projection. Recall that $\gamma_{{\rm mid}}$ is the mid point of the diameter $\lgeo \gamma_0, \gamma_1 \rgeo$: namely $d(\gamma_0, \gamma_{{\rm mid}})\! =\! d(\gamma_1, \gamma_{{\rm mid}})\! =\! \frac{_1}{^2}D$. Recall from Theorem \ref{diamlaw} ($ii$) 
that $\tau^-_{{\rm mid}}$ and $\tau^+_{{\rm mid}}$ are the only times $t\! \in \! [0, \zeta]$ such that $p(t)\! = \! \gamma_{{\rm mid}}$; thus, $\gamma_{{\rm mid}}$ is a simple point of $\cT$; namely, 
$\cT\backslash \{ \gamma_{{\rm mid}}\}$ has only two connected components. Denote by 
$\cT^{o}$ the connected component containing $\gamma$: it does not contain the root;  if we set $\cT^-\! =\!  \{ \gamma_{{\rm mid}}\}\cup \cT^{o}$ and $\cT^+\! = \! \cT \backslash \cT^{o}$, then $H^\ominus$ codes $(\cT^-, d,\gamma_{{\rm mid}})$ and $H^{\oplus}$ codes $(\cT^+, d,\gamma_{{\rm mid}})$.

\medskip

In the following lemma we recall 
Proposition 2.1 from D.~\& Le Gall \cite{Duquesne09} that asserts that $H$ is invariant under uniform re-rooting. Recall from (\ref{rerootH}) the definition of $H^{[t]}$.  
\begin{lem} 
[ D.~\& Le Gall \cite{Duquesne09}]
\label{recDuLG} For all measurable functions $F\! : \! \bbR_+ \! \times \! \bC(\bbR_+, \bbR_+)\! \rightarrow \! \bbR_+$ and $g: \bbR_+ \rightarrow \bbR_+$, 
$$\bN \Big[ g(\zeta) \! \int_0^\zeta \!\!\! dt\,  F\big( t, H^{[t]} \big)\Big]= \bN \Big[ g(\zeta)\!  \int_0^\zeta \!\!\! dt\,  F ( t, H) \Big] \; .$$
\end{lem}
By applying this property we first get 
\begin{equation}
\label{reroospec}
 \bN \big[ \zeta F_1 (H^{\ominus}) F_2 (H^{\oplus})\big] \! = \! \bN \Big[ \int_0^\zeta \!\!\! dt \, 
F_1 (H^{\ominus}) F_2 (H^{\oplus})\Big] \! = \! \bN \Big[ \int_0^\zeta \!\!\! dt \, 
F_1 \big( (H^{[t]})^{\ominus}) F_2 ((H^{[t]})^{\oplus})\Big] \; .
\end{equation}
Next observe the following: if $t \! \in \! (\tau^-_{{\rm mid}}, \tau^+_{{\rm mid}})$, then $(H^{[t]})^{\ominus}\! =\!  H^{\oplus}$ and $(H^{[t]})^{\oplus}\! =\!  H^{\ominus}$, and if 
$t \! \in \! (0, \tau^-_{{\rm mid}}) \cup (\tau^+_{{\rm mid}} , \zeta)$, then 
$(H^{[t]})^{\ominus}\! =\!  H^{\ominus}$ and $(H^{[t]})^{\oplus}\! =\!  H^{\oplus}$. Thus, 
\begin{eqnarray*}
\int_0^\zeta \!\!\! dt \, 
F_1 \big( (H^{[t]})^{\ominus}) F_2 ((H^{[t]})^{\oplus}) \!\!  & = & \!\!  \big(\tau^+_{{\rm mid}}\!\! -\! \tau^-_{{\rm mid}} \big) F_1 (H^{\oplus}) F_2(H^{\ominus}) +  \big( \zeta \! -\! \tau^+_{{\rm mid}}\! +\tau^-_{{\rm mid}} \big) F_1 (H^{\ominus}) F_2 (H^{\oplus}) \\
& =& \!\!  \zeta_{H^{\ominus}} F_1 (H^{\oplus}) F_2(H^{\ominus}) + \zeta_{H^{\oplus}} F_1 (H^{\ominus}) F_2 (H^{\oplus}). 
\end{eqnarray*}
This equality, (\ref{reroospec}) and (\ref{demicas}) with $f\equiv 1$ imply the following: 
\begin{eqnarray}
\label{dedemi}
 \bN \big[ \zeta F_1 (H^{\ominus}) F_2 (H^{\oplus})\big] \!\! &=&  \! \!  \bN\big[ \zeta_{H^{\ominus}} F_1 (H^{\oplus}) F_2(H^{\ominus})\big]  + \bN \big[ \zeta_{H^{\oplus}} F_1 (H^{\ominus}) F_2 (H^{\oplus}) \big]  \nonumber \\
 &= &  \!\! \bN \Big[ \,  \bN^{\Gamma}_{\! \frac{1}{2}\!   D} \big[ \zeta F_2 (H) \big] \;  F_1 (H^{\oplus}) \Big] + 
 \bN \Big[\,   \bN^{\Gamma}_{\! \frac{1}{2}\!   D} [ F_1 (H)] 
\;   \zeta_{H^{\oplus}}   F_2\big( H^{\oplus }\big)\Big] \; .
\end{eqnarray}
Next observe that $\zeta_{H^{\ominus}}+  \zeta_{H^{\oplus}}\! = \! \zeta$. Thus, by (\ref{demicas}) we also get 
\begin{eqnarray}
\label{demimi}
 \bN \big[ \zeta F_1 (H^{\ominus}) F_2 (H^{\oplus})\big] \!\! &=&  \! \!  \bN\big[ \zeta_{H^{\ominus}} F_1 (H^{\ominus}) F_2(H^{\oplus})\big]  + \bN \big[ \zeta_{H^{\oplus}} F_1 (H^{\ominus}) F_2 (H^{\oplus}) \big]  \nonumber \\
 &= &  \!\! \bN \Big[ \,  \bN^{\Gamma}_{\! \frac{1}{2}\!   D} \big[ \zeta F_1 (H) \big] \;  F_2 (H^{\oplus}) \Big] + 
 \bN \Big[\,   \bN^{\Gamma}_{\! \frac{1}{2}\!   D} [ F_1 (H)] 
\;   \zeta_{H^{\oplus}}   F_2\big( H^{\oplus }\big)\Big] \; .
\end{eqnarray}
Then by (\ref{dedemi}) and (\ref{demimi}), we get 
$ \bN \big[ \,  \bN^{\Gamma}_{\! \frac{1}{2}\!   D} \big[ \zeta F_1 (H) \big] \;  F_2 (H^{\oplus}) \big]  \! =\! \bN \big[ \,  \bN^{\Gamma}_{\! \frac{1}{2}\!   D} \big[ \zeta F_2 (H) \big] \;  F_1 (H^{\oplus}) \big]$. Since the total height of $H^{\ominus}$ and $H^{\oplus}$ is $\frac{_1}{^2} D$, for all measurable functions $F_1, F_2 \! :\!  \bC(\bbR_+, \bbR_+) \! \rightarrow \! \bbR_+$, $f: \bbR_+ \rightarrow \bbR_+$, we get 
\begin{equation}
\label{dyssim}
\bN \Big[  f(D) \,  \bN^{\Gamma}_{\! \frac{1}{2}\!   D} \big[ \zeta F_1 (H) \big] \;  F_2 (H^{\oplus}) \Big]=  \bN \Big[  f(D) \,  \bN^{\Gamma}_{\! \frac{1}{2}\!   D} \big[ \zeta F_2 (H) \big] \;  F_1 (H^{\oplus}) \Big]  . 
\end{equation}  
By taking in (\ref{dyssim}) 
$F_1 \equiv 1$ and by substituting $f (D)$ with $f(D) / \bN^\Gamma_{\frac{1}{2}\!   D} [\,  \zeta \, ]$, we get 
$$  \bN \big[  f(D) \, F_2 (H^{\oplus}) \big]= \bN \Big[ f(D) \, 
\bN^{\Gamma}_{\! \frac{1}{2}\!   D} \big[ \zeta F_2 (H) \big]  \big/ \bN^\Gamma_{\frac{1}{2}\!   D} [\,  \zeta \, ] \Big] ,  $$
and by (\ref{demicas}), it entails 
\begin{equation} 
 \label{fullcas}
\bN \big[  f(D) \, F_1(H^{\ominus}) F_2 (H^{\oplus}) \big]= \bN \Big[ f(D) \, 
\bN^{\Gamma}_{\! \frac{1}{2}\!   D} \big[ F_1 (H) \big] \, 
\bN^{\Gamma}_{\! \frac{1}{2}\!   D} \big[ \zeta F_2 (H) \big]  \big/ \bN^\Gamma_{\frac{1}{2}\!   D} [\,  \zeta \, ] \Big]  \; .
\end{equation} 
Recall from (\ref{idnotapff}) that $H^{[\tau^-_{{\rm mid}}]}= H^{\ominus} \oplus H^{\oplus}$. Then, (\ref{fullcas}) implies for all measurable functions $F \! :\!  \bC(\bbR_+, \bbR_+) \! \rightarrow \! \bbR_+$, $f: \bbR_+ \rightarrow \bbR_+$, that 
\begin{align}
\label{devidys}
 \bN \big[ f(D)& F\big(H^{[\tau^-_{\rm mid}]} \big)\big] = \\ 
&  \int_0^\infty \!\!\!  \bN(D \! \in \! dr)
 \frac{f(r)}{\bN^{\Gamma}_{r/2} [\, \zeta \, ] }   \int \!\!\!\! \int_{\bC(\bbR_+, \bbR_+)^2} \!\!\!\!\!\!\!\! \!\!\!\!\!\!\!\! \!\!\!\!\! \!   \bN^{\Gamma}_{r/2}( dH) \bN^{\Gamma}_{r/2}( dH^\prime)  \; 
 \,  \zeta_{H^\prime}   F\big( H \! \oplus \!  H^\prime \big) \; , \nonumber
\end{align} 
which implies Theorem \ref{thm: decomp} $(iii)$ as soon as one makes sense of $\bN(\, \cdot \, | \, D\! = \! r)$. \cq 

\paragraph{Proof of Theorem \ref{thm: decomp} $(ii)$.} Recall that $\Lambda \! :\!  \mathtt{Exc}\! \rightarrow \!  \mathtt{Exc}$ is the functional that reverses excursions at their lifetime: namely 
for all $H\! \in \! \mathtt{Exc}$, $\Lambda (H)= (H_{(\zeta_{H}-t)_+})_{t\geq 0}$. Recall from (\ref{revHr}) that for all $r\! \in \! (0, \infty)$, $H$ and $\Lambda(H)$ have the same law under $\bN^\Gamma_r$, which entails the following by (\ref{fullcas}): 
\begin{equation}
\label{revHplms}
\textrm{$\big( \Lambda( H^{\ominus} ), \Lambda(H^{\oplus}) \big)$ and $(H^{\ominus}, H^{\oplus})$ have the same distribution under $\bN$.}
\end{equation}
Next, observe that $D(\Lambda (H))\! = \! D$, 
$\tau (\Lambda(H))\! = \! \zeta \!-\!\tau$, $\tau_0 (\Lambda (H))\! = \! \zeta \!-\!\tau_1$ and 
$\tau_1 (\Lambda(H))\! = \! \zeta \!-\!\tau_0$. Moreover, $(\Lambda (H))^{\ominus}= \Lambda( H^{\ominus})$ and   $(\Lambda (H))^{\oplus}= \Lambda( H^{\oplus})$. This combined with (\ref{revHplms}) and  (\ref{fullcas})  implies that 
\begin{eqnarray}
\label{inforac}
\frac{_1}{^2} \bN \big[  f(D)  F_1(H^{\ominus}) F_2 (H^{\oplus}) \big] & =&  \bN \big[  \un_{\{ \tau= \tau_0\} }f(D)  F_1(H^{\ominus}) F_2 (H^{\oplus}) \big] \\
& =&  \bN \big[ \un_{\{ \tau= \tau_1\} } 
f(D)  F_1(H^{\ominus}) F_2 (H^{\oplus}) \big]\; . \nonumber 
\end{eqnarray}
We then define 
$$ \textrm{$\tau^*\! :=\! \tau^-_{{\rm mid}}$ if $\tau\! = \! \tau_0$} \quad \textrm{and} \quad  
\textrm{$\tau^*\! :=\! \tau^+_{{\rm mid}}$ if $\tau\! = \! \tau_1$.}$$ 
By (\ref{idnotapff}), we get 
$$ H^{[\tau^*]}= H^{\ominus } \oplus H^{\oplus} \; \, \textrm{on $\{ \tau\! = \! \tau_0\}$} \quad \textrm{and} \quad  H^{[\tau^*]}= H^{\oplus } \oplus H^{\ominus} \; \, \textrm{on $\{ \tau\! = \! \tau_1\}$.} $$
This, combined with (\ref{inforac}) and (\ref{fullcas}) entails 
\begin{align}
\label{Htaustar}
 \bN \big[ f(D)& F\big(H^{[\tau^*]} \big)\big] = \\ 
&  \int_0^\infty \!\!\!  \bN(D \! \in \! dr)
 \frac{f(r)}{2\bN^{\Gamma}_{r/2} [\, \zeta \, ] }   \int \!\!\!\! \int_{\bC(\bbR_+, \bbR_+)^2} \!\!\!\!\!\!\!\! \!\!\!\!\!\!\!\! \!\!\!\!\! \!   \bN^{\Gamma}_{r/2}( dH)\,  \bN^{\Gamma}_{r/2}( dH^\prime)  \; 
 \,  (\zeta_H +\zeta_{H^\prime})   F\big( H \! \oplus \!  H^\prime \big) \; . \nonumber
\end{align} 
Recall from (\ref{defQr}) the definition of the law $\bQ_r$. Since $r\mapsto \bN^\Gamma_r$ is weakly continuous, 
it is easy to check that $r\mapsto \bQ_r $ is also weakly continuous. Then observe that $\bQ_r[\, \zeta \, ]= 2 \bN^\Gamma_{r/2} [\, \zeta \,]$. Therefore (\ref{Htaustar}) can be rewritten as 
\begin{equation}
\label{HtaustarQ}
 \bN \big[ f(D) F\big(H^{[\tau^*]} \big)\big] =  \int_0^\infty \!\!\!  \bN(D \! \in \! dr)\, f(r) \, 
 \bQ_r \big[ \zeta F(H) \big] \big/ \bQ_r[\, \zeta \, ] \; .
\end{equation} 
Next observe that for all $t \! \in \! [0, \zeta]$, $ (H^{[\tau^*]})^{[t]}= H^{[\tau^*+t]}$ and that $D(H^{[t]})\! = \! D$. Thus, (\ref{HtaustarQ}) implies  
\begin{eqnarray*}
\int_0^\infty \!\!\! \!\!  \bN(D \! \in \! dr)\, f(r)  
 \bQ_r \! \Big[ \zeta\!\! \int_0^\zeta \!\!\! dt \, F\big(H^{[t]} \big) \Big] \big/ \bQ_r[\, \zeta \, ] & =& \bN \Big[ f(D)\!\! \int_0^\zeta \!\!\! dt \, F\big(H^{[\tau^*+t]} \big)\Big]  \\
 & = &   \bN \Big[\! \int_0^\zeta \!\!\! dt \, f\big( D\big( H^{[t]} \big)\big)F\big(H^{[t]} \big)\Big] \\
 & =& \bN \big[ \zeta f(D) F(H)\big] \; , 
\end{eqnarray*}
where we have used Lemma \ref{recDuLG} in the last line. This proves (\ref{desindiam}) in Theorem \ref{thm: decomp} $(ii)$. 

\paragraph{Proof of Theorem \ref{thm: decomp} $(i)$ and $(iv)$.} The rest of the proof is now easy: we fix $r\! \in \! (0, \infty)$ and we denote by $\Pi_r (dH^\prime\! dH^{\prime \prime})$ the product 
law $\bN^\Gamma_{r/2} (dH^\prime) \bN^\Gamma_{r/2} (dH^{\prime \prime})$; we then set $H\! = \! 
H^\prime\! \oplus \! H^{\prime \prime}$. Thus, by definition, $H$ under $\Pi_r$ has law $\bQ_r$. Observe that if $t \! \neq\! \tau (H^\prime)$ (resp.~$t \! \neq\! \tau (H^{\prime\prime})$) then $H^\prime_{\! t}  \! < \! r/2$ 
(resp.~$H^{\prime \prime}_{\! t}  \! <\! r/2)$. Note that if $s \! \in \! [0, \zeta_{H^{\prime}}]$ and $t \! \in \! [\zeta_{H^{\prime}}, \zeta_{H^{\prime\prime}}]$, then $\inf_{[s, t]} H \! =\! 0 $ and $d_H(s, t)\! =\!  H^{\prime}_s + 
H^{\prime\prime}_{t-\zeta_{H^\prime}}$. This easily entails that $\Pi_r$-a.s.~$D(H)\! = \! r $ and that 
$\tau (H^{\prime})$ and $\zeta_{H^\prime}+ \tau (H^{\prime\prime})$ are the two only times $s\! <\! t$ such that $d_H(s, t) \! = \! D(H)$, which completes the proof of  \ref{thm: decomp} $(i)$. 

  The fact that $\bQ_r$-a.s.~$D\! = \! r$, combined with (\ref{desindiam}) and with the fact that $r \mapsto \bQ_r$ is weakly continuous, allows to make sense of $\bN(\, \cdot \, | \, D\! = \! r)$ that is a regular version of the conditional distribution of 
$\bN$ knowing that $D\! = \! r$. Moreover, (\ref{desindiam}) entails (\ref{Hcondia}) for all $r\! \in \! (0, \infty)$. Furthermore (\ref{devidys}) entails (\ref{paraph}) that was the last point to clear in the Theorem \ref{thm: decomp} $(iii)$, as already mentioned. 

It remains to prove Theorem \ref{thm: decomp} $(iv)$. We keep the previous notations and we introduce the following: 
$$ \cM_{0, \tau(H^{\prime})} (H^{\prime})= \sum_{a\in \cJ_{0, \tau^\prime}} \delta_{ (a, \overleftarrow{H}^{ a},  \overrightarrow{H}^{ a} )} \quad \textrm{and} \quad   \cM_{0, \tau(H^{\prime \prime})} (H^{\prime \prime})= \sum_{a\in \cJ_{0, \tau^{\prime \prime}}} \delta_{ (a, \overleftarrow{H}^{ a},  \overrightarrow{H}^{ a}) }\; , $$
that are under $\Pi_r$ independent Poisson point measures with the same intensity $\mathbf{n}_{r/2}$, by Lemma \ref{recAbDe}. 
We then set $\tau_0 (H)\! :=\!  \tau(H^{\prime})$ and $\tau_1(H)\! :=\! \zeta_{H^{\prime}} +\tau(H^{\prime \prime})$, that are the only pair of times realizing the diameter $D(H)$ under $\Pi_r$, as already shown. Observe that under $\Pi_r$, 
$$ \cM_{\tau_0(H), \tau_1(H)} (H)=  \sum_{a\in \cJ_{0, \tau^\prime}} \delta_{ (r-a, \Lambda (\overrightarrow{H}^{ a}),  \Lambda (\overleftarrow{H}^{ a}) )} +   \cM_{0, \tau(H^{\prime \prime})} (H^{\prime \prime}) \; ,$$
where we recall here that $\Lambda$ reverses excursions at their lifetime and that $\Lambda$ is invariant under $\bN_a$ and $\bP^x_a$. Thus, basic results on Poisson point measures and an 
easy calculation show that $ \cM_{\tau_0(H), \tau_1(H)} (H)$ is a Poisson point measure whose intensity is given by (\ref{decomQr}) in Theorem \ref{thm: decomp} $(iv)$, which completes the proof of \ref{thm: decomp} $(iv)$ because $H$ under $\Pi_r$ has law $\bQ_r$ and thus $\cM_{\tau_0(H), \tau_1(H)} (H)$ under $\Pi_r$ has the same law as $\cM_{\tau_0, \tau_1}$ under 
$\bQ_r$. This completes the proof of Theorem  \ref{thm: decomp}. \cqfd

\section{Total height and diameter of normalized stable trees.}
\label{sec: stable}
\subsection{Preliminary results.}
\label{preprelisec}
In this section, we gather general results that are used to prove Proposition \ref{prop: lp}. 
\textit{Unless the contrary is explicitly mentioned, $\Psi$ is a general branching mechanism of the form (\ref{eq: def-Psi}) that satisfies (\ref{eq: hyp})}.
We first introduce the following function 
\begin{equation}
\label{wlam}
\forall \, a, \lambda \in (0, \infty) , \quad w_\lambda (a):= \bN \big[ 1-\un_{\{\Gamma \leq a\}} e^{-\lambda \zeta} \big]\; .
\end{equation}
For all fixed $\lambda\! \in \! (0, \infty)$, note that $a\mapsto w_\lambda (a)$ is non-increasing, that $\lim_{a\rightarrow 0} w_\lambda (a)\! = \! \infty$ and by (\ref{lifetimeexc}) $\lim_{a\rightarrow \infty} w_\lambda (a)\! = \! \bN[1\! -\! e^{-\lambda \zeta}]= \Psi^{-1} (\lambda)$. As proved by Le Gall \cite{LG99}, Section II.3 (in the more general context of superprocesses) $w_\lambda (a)$ is the only solution of the following integral equation, 
\begin{equation}
\label{inteqw}   
\forall \, a, \lambda \in (0, \infty) , \quad \int_{w_\lambda(a)}^\infty \! \! \frac{du}{\Psi (u)\! -\! \lambda}= a \; ,
\end{equation}
that makes sense thanks to (\ref{eq: hyp}). 

Let us next consider $H$ under $\bbP$ and recall from (\ref{xLevfo}) that $\bP^x$ stands for the law of $H_{\cdot \wedge T_{x}}$ where $T_x= \inf \{ t\! \in \! \bbR_+: X_t\! =\!-x\}$. 
Recall from (\ref{PoisdecH}) that $\sum_{i\in \cI} \delta_{(-I_{a_i} , H^{i})}$ stands for the decomposition of $H$ into excursions above $0$; thus, 
the excursions of $H_{\cdot \wedge T_{x}}$ above $0$ are the $H^{i}$ where $i\! \in \! \cI$ is such that $-I_{a_i} \! \in \! [0, x]$. Elementary results on Poisson point processes then imply the following: 
\begin{eqnarray}
\label{sousPax}
\bE^x_a \big[ e^{-\lambda \zeta}\big] &= & \bE^x \big[ e^{-\lambda \zeta}\un_{\{ \Gamma \leq a \}} \big]  
=  \bbE \Big[\exp \Big( \! -\!\! \sum_{i\in \cI} \lambda \zeta_{H^{i}} \un_{[0, x]} (-I_{a_i})\Big)\un_{\big\{ \Gamma (H^{i}) \leq a \, , \, i \in \cI: -I_{a_i} \leq x \big\} } \Big] \nonumber  \\
& =& \exp \big(\! -\! x w_\lambda (a)  \big) \; .
\end{eqnarray}
We first prove the following lemma. 
\begin{lem}
\label{calculw} Let $\Psi$ be a branching mechanism of the form (\ref{eq: def-Psi}) that satisfies (\ref{eq: hyp}). Recall from (\ref{wlam}) the definition of $w_\lambda (a)$. First observe that for all $a, \lambda \! \in\!  (0, \infty)$, 
\begin{equation}
\label{deriv} \partial_a w_\lambda (a) = \lambda\! -\!  \Psi (w_\lambda (a) ) \quad \textrm{and} \quad \int_{w_\lambda(a)}^\infty \! \! \frac{du}{(\Psi (u)\! -\! \lambda)^2}= \frac{ \partial_\lambda w_\lambda (a)}{\Psi (w_\lambda (a) ) \! -\! \lambda } \; . 
\end{equation}  
Recall from (\ref{lawGamma}) the definition of the function $v$. Then, for all  $a, \lambda \! \in\!  (0, \infty)$,
\begin{equation}
\label{brinbrin}
\lim_{\lambda \rightarrow 0+} w_\lambda (a)= v(a)  \quad \textrm{and} \quad v(a) \leq w_\lambda(a) = v(a) + \bN_a \big[ 1\!-\! e^{-\lambda \zeta}\big] \leq v(a) + \Psi^{-1} (\lambda) \; , 
\end{equation} 
where we recall from (\ref{renota}) the notation $\bN_a$. Then, for all $r_1\! \geq \! r_0 \! >\! 0$, we get 
\begin{equation}
\label{intpsiw}
\int_{r_0}^{r_1} \!\!\! da \, \Psi^\prime (w_\lambda (a)) = \log \frac{\Psi (w_\lambda(r_0))\! -\! \lambda }{\Psi (w_\lambda (r_1))\! -\! \lambda } 
\quad \textrm{and} \quad 
\int_{r_0}^{r_1} \!\!\! da \, \Psi^\prime (v(a)) = \log \frac{\Psi (v(r_0))}{\Psi (v(r_1))} \; .
\end{equation}
\end{lem}
\noi
\textbf{Proof.} Note that (\ref{deriv}) and (\ref{brinbrin}) are easy consequences of resp.~(\ref{inteqw}) and the definition (\ref{wlam}). Let us first prove the first equality of (\ref{intpsiw}): to that end we use the change of variable $u= w_\lambda(a)$, $\lambda$ being fixed. Then, by (\ref{deriv}), $-du/ (\Psi (u)\!- \! \lambda)= da$, and we get 
$$ \int_{r_0}^{r_1} \!\!\! da \, \Psi^\prime (w_\lambda (a))= \int_{w_\lambda (r_1)}^{w_\lambda (r_0)} \!\!\!\!\! du \, \frac{\Psi^\prime (u)}{\Psi (u)\! -\! \lambda} = \log \frac{\Psi (w_\lambda(r_0))\! -\! \lambda}{\Psi (w_\lambda (r_1))\! -\! \lambda} \; ,  $$
which implies the second equality in (\ref{intpsiw}) as $\lambda \! \rightarrow \! 0$ by (\ref{brinbrin}). \cqfd 
\begin{prop}
\label{Ngamrlaws} Let $\Psi$ be a branching mechanism of the form (\ref{eq: def-Psi}) that satisfies (\ref{eq: hyp}). Let $r\! \in \! (0, \infty)$. Recall from (\ref{htcondsim}) the definition of $\bN^\Gamma_r$ and  
recall from (\ref{wlam}) the definition of $w_\lambda (a)$. Then for all $\lambda \! \in \! (0, \infty)$, we first get 
\begin{equation}
\label{Ngamass}
\bN^\Gamma_r \big[ e^{-\lambda \zeta}\big]= \exp \Big(\!\!  -\!\!  \int_{0}^{r} \!\!\! da \, \big( \Psi^\prime (w_\lambda (a))\! -\!  \Psi^\prime (v(a))\big)\Big) =   \frac{\Psi (w_\lambda (r)) \! -\! \lambda }{\Psi (v(r))}\; .
\end{equation}
We next set $q_\lambda(y,r)\! := \!  \bN^\Gamma_r \big[ e^{-\lambda \zeta}\un_{\{ D>2y\}}\big] $. 
Then for all $y\! \in \! (\frac{_1}{^2} r, r)$,  we have 
\begin{equation}
\label{Ngamdiam}
q_\lambda(y,r)= \frac{\Psi(w_\lambda (r)) \!-\!  \lambda}{\Psi(v(r))} \, \Big(1- \frac{ \big( \Psi (w_\lambda (y))  \!-\!  \lambda \big)^2}{\big( \Psi(w_\lambda (2y\!-\! r)) \!-\!  \lambda \big) \big( \Psi(w_\lambda (r)) \!-\!  \lambda \big)} \Big)\; .
\end{equation}
If $y\! \leq  \! \frac{_1}{^2} r$, then $q_\lambda(y,r)\! = \! \bN^\Gamma_r \big[ e^{-\lambda \zeta}\big]$ and if $y\! >\! r$, then $q_\lambda(y,r)\! = \! 0$. 
\end{prop}
\noi
\textbf{Proof.} Recall from (\ref{notaMtau}) the notation $\cM_{0, \tau}$ 
and recall from (\ref{notadelt}) the notation $\Delta_{b, \overleftarrow{H},  \overrightarrow{H} }$. Then, for all $r,y,\lambda \! \in (0, \infty)$, we get $\bN^\Gamma_r$-a.s.
$$e^{-\lambda \zeta}\un_{\{ D\leq 2y\}}= \exp \Big( \!\! -\! \lambda \!\! \sum_{a\in \cJ_{0, \tau}} \!\!\! 
\big( \zeta_{ \overleftarrow{H}^a}+ \zeta_{ \overrightarrow{H}^a} \big)\Big)\un_{\big\{ \forall a \in \cJ_{0, \tau}\, :\,  
\Delta_{a, \overleftarrow{H}^a,  \overrightarrow{H}^a } \leq 2y\big\}} \; .$$
Lemma \ref{recAbDe} asserts that under $\bN^\Gamma_r$, $\cM_{0, \tau}$ is a Poisson point measure with intensity $\mathbf{n}_r$ given by (\ref{decomGr}). Thus, 
elementary results on Poisson point measures imply that 
$$ 
\bN^\Gamma_r \big[ e^{-\lambda \zeta}\un_{\{ D\leq 2y\}}\big]= \exp \Big(\!\! -\!\! \underset{K}{\underbrace{\int \! \mathbf{n}_r (da d\overleftarrow{H}d\overrightarrow{H} )  \big( 1-\un_{\{ \Delta_{a, \overleftarrow{H},  \overrightarrow{H} }\leq 2y\}} e^{-\lambda\zeta_{ \overleftarrow{H}}-\lambda \zeta_{ \overrightarrow{H}}  }} \big)}\Big). 
$$
Recall that the total mass of $\bP^x_a$ is $e^{-xv(a)}$ and recall (\ref{sousPax}). Thus,  
$$ K\! = \! \int_0^r \!\!\!  da \, 2\beta\,  \bN_a \big[ 1\! -\! \un_{\{ \Gamma \leq  2y-a \}}e^{-\lambda \zeta} \big] + \int_0^r \!\!\!  da \! 
\int_{(0, \infty)} \!\!\!\!\!\! \!\!\! \! \pi(dz) \, z \big( e^{-zv(a)} -e^{-zw_\lambda (a\wedge (2y-a))}\big) . $$
Now observe that 
$$ \bN_a \big[ 1\! -\! \un_{\{ \Gamma \leq  2y-a \}}e^{-\lambda \zeta} \big] \! = \! \bN \big[ 1\! -\! \un_{\{ \Gamma \leq a\wedge ( 2y-a) \}}e^{-\lambda \zeta} \big] \! -\!  \bN[ \un_{\{ \Gamma >a \}}\big]= w_\lambda \big( a \! \wedge \! (2y\!-\! a)\big) \! -\! v(a) \; .$$
Consequently, 
\begin{equation} 
\label{PoNGadia}
\bN^\Gamma_r \big[ e^{-\lambda \zeta}\un_{\{ D\leq 2y\}}\big] = \exp \Big(\!\!  -\!\!  \int_{0}^{r} \!\!\! da \, \big( \Psi^\prime (w_\lambda ( a \! \wedge \! (2y\!-\! a)))\! -\!  \Psi^\prime (v( a))\big) \Big). 
\end{equation}
Then observe that if $y\! >\! r$, the $\bN^\Gamma_r \big[ e^{-\lambda \zeta}\un_{\{ D\leq 2y\}}\big] =\bN^\Gamma_r \big[ e^{-\lambda \zeta}\big]$ because $D\! \leq \! 2\Gamma$. This combined with (\ref{PoNGadia}) entails the first equality of (\ref{Ngamass}). Then, use (\ref{intpsiw}) in Lemma \ref{calculw} to get for any $\varepsilon \! \in \! (0, r)$, 
$$\int_{\varepsilon}^{r} \!\!\! da \, \big( \Psi^\prime (w_\lambda (a))\! -\!  \Psi^\prime (v(a))\big)= 
\log \frac{\Psi (v(r))}{\Psi (w_\lambda (r))\! -\! \lambda}-\log \frac{\Psi (v(\varepsilon))}{\Psi (w_\lambda (\varepsilon))\! -\! \lambda} .$$
This show that $\varepsilon \mapsto \Psi (v(\varepsilon))/(\Psi (w_\lambda (\varepsilon))\! -\! \lambda)$ is increasing and tends to a finite constant $C_\lambda\! \in \! (0, \infty)$ as $\varepsilon \! \rightarrow \! 0$. Then, 
$ C_\lambda^{-1}\Psi (v(r))\bN^\Gamma_r \big[ e^{-\lambda \zeta}\big] \! =\! \Psi (w_\lambda (r)) \! -\! \lambda$, which is equal to 
$-\partial_r w_\lambda (r)$ by (\ref{deriv}) in Lemma \ref{calculw}. Then recall from (\ref{lawGamma}) that 
$\bN(\Gamma \! \in \! dr) \! =\! \Psi (v(r))\, dr$; thus by (\ref{heightcondef}) and the fact that $w_\lambda(r)$ tends 
to $\Psi^{-1} (\lambda)$ as $r\! \rightarrow \! \infty$, we get for all $b\! \in \! (0, \infty)$, 
\begin{eqnarray*}
 w_\lambda(b)\! -\! \Psi^{-1} (\lambda)& = & 
\int_b^\infty \!\!\! dr \,  C_\lambda^{-1}\Psi (v(r))\bN^\Gamma_r \big[ e^{-\lambda \zeta}\big] = 
  C_\lambda^{-1} \bN \big[ e^{-\lambda \zeta}\un_{\{ \Gamma > b\}}\big]  \\
&= & C_\lambda^{-1}\big(  \bN\big[ 1\!-\! \un_{\{ \Gamma \leq b\}}e^{-\lambda \zeta} \big] \!-\! \bN\big[ 1\!-\!e^{-\lambda \zeta} \big] \big)=  C_\lambda^{-1}\big( w_\lambda(b)\! -\! \Psi^{-1} (\lambda) \big). 
 \end{eqnarray*}
This implies that $C_\lambda\! = \! 1$, which completes the proof of (\ref{Ngamass}). 

We next assume that $y\! \in \! (\frac{_1}{^2} r, r)$. Observe that $ a \! \wedge \! (2y\!-\! a)\! = \! a$ if $a \! \in \! (0, y)$ and that $a \! \wedge \! (2y\!-\! a)\! = \! 2y\!-\! a$ if $a \! \in \! (y, r)$. 
By (\ref{PoNGadia}) and (\ref{Ngamass}), we then get 
$$ q_{\lambda} (y,r) = \bN^\Gamma_r \big[ e^{-\lambda \zeta}\big]-\bN^\Gamma_r \big[ e^{-\lambda \zeta}\un_{\{ D\leq 2y\}}\big]=  \frac{\Psi (w_\lambda (r)) \! -\! \lambda }{\Psi (v(r))} \Big( 1- e^{-  \int_{y}^{r} \! da \, \big( \Psi^\prime (w_\lambda ( 2y\!-\! a))\! -\!  \Psi^\prime (w_\lambda (a) \big) } \Big) , $$
which easily implies (\ref{Ngamdiam}) by (\ref{intpsiw}) in Lemma \ref{calculw} since 
$$  \int_{y}^{r}\!\! \! \! da \Psi^\prime (w_\lambda ( 2y\!-\! a))\! =\!\!  \int_{2y-r}^y \!\!\!\!\!\!\!\!\! da  \Psi^\prime (w_\lambda (a))\! =\! \log  \frac{\Psi (w_{\lambda} (2y\! -\! r) ) \!-\! \lambda}{\Psi (w_{\lambda} (y) ) - \lambda} \; \textrm{and} \, \int_{y}^r \! \!\!\! da  \Psi^\prime (w_\lambda (a))\! =\! \log  \frac{\Psi (w_{\lambda} (y) ) \!-\! \lambda}{\Psi (w_{\lambda} (r) ) \!-\! \lambda}. $$
The other statements of the lemma follow immediately. \cqfd 

\begin{prop}
\label{genjoint} Let $\Psi$ be a branching mechanism of the form (\ref{eq: def-Psi}) that satisfies (\ref{eq: hyp}). For all $y,z,\lambda \! \in \! (0, \infty)$, we have 
\begin{align}
\label{jointlawc}
{\rm L}_\lambda  (y,z) :=  \bN \big[ e^{-\lambda \zeta}& \un_{\{D>2y \, ; \,  \Gamma >z\}}\big] \nonumber \\
= & w_\lambda (y\! \vee \! z) - \Psi^{-1} (\lambda)\, - \, \un_{\{ z\leq 2y\}} \big( \Psi (w_\lambda (y))  \!-\!  \lambda \big)^2 \!\! \int^{\infty}_{w_{\lambda} (y\wedge (2y-z))}
\frac{du}{\big(\Psi(u) \!-\!  \lambda\big)^2} \\
= &  w_\lambda (y\! \vee \! z) - \Psi^{-1} (\lambda)\, - \, \un_{\{ z\leq 2y\}} \big( \Psi (w_\lambda (y))  \!-\!  \lambda \big)^2 \frac{\partial_\lambda w_\lambda \big( y \! \wedge\! (2y\! -\! z) \big)}{
 \Psi \big(w_\lambda (y \! \wedge\! (2y\! -\! z))\big) \! -\! \lambda } \; . \nonumber 
\end{align}
\end{prop}
\noi
\textbf{Proof.} Recall notation $q_{\lambda} (y,r)$ from Proposition \ref{Ngamrlaws}, which asserts that 
$q_{\lambda} (y,r)\!= \! 0$ if $r\! <\! y$ and that $\Psi (v(r))q_{\lambda} (y,r)\! = \! - \partial_r w_\lambda (r)$, if $r \! \geq\! 2y$. 
Then, by (\ref{heightcondef}), we get 
\begin{equation}
\label{micmaczon}
 {\rm L}_{\lambda}(y,z)= \int_z^\infty \!\!\!\! dr \,  \Psi (v(r)) q_\lambda (y,r)\! =\! \un_{\{ z\leq 2y \}} \int_{z\vee y}^{2y} dr \,  \Psi (v(r)) q_\lambda (y,r) 
 - \int_{z\vee 2y }^\infty  \!\!\!\!\!\!\!\! dr \, \partial_r w_\lambda (r) . 
\end{equation} 
Since 
$\lim_{r\rightarrow \infty} w_\lambda (r)\! = \! \Psi^{-1} (\lambda)$, we get 
\begin{equation}
\label{trivzon}
-\int_{z\vee 2y }^\infty  \!\!\!\!\! dr \, \partial_r w_\lambda (r) = w_\lambda (z\vee 2y) -\Psi^{-1} (\lambda) \; . 
\end{equation}
We next assume that $z \! \in \! (y, 2y)$. By (\ref{Ngamdiam}) and since $\Psi (w_\lambda (r)) \! -\! \lambda\! = \! - \partial_r w_\lambda (r)$, we get 
\begin{eqnarray*}
\int_{z}^{2y} dr \,  \Psi (v(r)) q_\lambda (y,r) & = & \! -\! \int_z^{2y} \!\!\!\!  dr\,  \partial_r w_\lambda (r) -\big( \Psi (w_\lambda (y))  \!-\!  \lambda \big)^2\!\int_z^{2y}  \!\!\!\!  \frac{dr}{\Psi(w_\lambda (2y\!-\! r)) \!-\!  \lambda} \\
& =& w_\lambda (z) \! -\! w_\lambda (2y)- \big( \Psi (w_\lambda (y))  \!-\!  \lambda \big)^2\!\int_{0}^{2y-z} \!\!\!\!  \!\!\!\!  \frac{dr}{\Psi(w_\lambda (r)) \!-\!  \lambda} \\
& =& w_\lambda (z) \! -\! w_\lambda (2y)- \big( \Psi (w_\lambda (y))  \!-\!  \lambda \big)^2\!\int^{\infty}_{w_{\lambda} (2y-z)}
\frac{du}{\big(\Psi(u) \!-\!  \lambda\big)^2} \; , 
\end{eqnarray*}
with the change of variable $u\! =\! w_\lambda (r)$ in the last line. This combined with (\ref{trivzon}) easily 
entails the first equality in (\ref{jointlawc}). The second one follows from (\ref{deriv}) in Lemma \ref{calculw}.  \cqfd 

\subsection{Proof of Proposition \ref{prop: lp}.}
\label{pfproplp}
In this section,  we fix $\gamma \! \in \! (1, 2]$ and we take $\Psi (\lambda)\! = \! \lambda^\gamma$, $ \lambda \! \in \! \bbR_+$. 
Recall from (\ref{wlam}) the definition of $w_\lambda (a)$. We then set 
\begin{equation}
\label{wwlam}
\forall y \! \in \! (0, \infty) , \quad w(y):= w_1(y) \; .
\end{equation}
Note that $w$ satisfies (\ref{wdef1}) that is (\ref{inteqw}) with $\lambda \! = \! 1$. By an easy change of variable (\ref{inteqw}) implies that 
\begin{equation}
\label{scalew}
\forall a, \lambda \in (0, \infty), \quad w_\lambda (a)= \lambda^{\frac{1}{\gamma}} w \big( a\lambda^{\frac{\gamma-1}{\gamma}} \big)\; .
\end{equation}
Recall from Proposition \ref{genjoint} the definition of ${\rm L}_{\lambda} (y,z)$. Then observe that the scaling property (\ref{echtscal}) entails (\ref{eq: def_bL}). Moreover (\ref{scalebL}) follows from a simple change of variable. 
Next note from (\ref{scalew}) that 
$$ \partial_\lambda w_\lambda (a)= \frac{_1}{^\gamma} \lambda^{\frac{1}{\gamma}-1} 
w \big( a\lambda^{\frac{\gamma-1}{\gamma}} \big)+  \frac{_{\gamma-1}}{^\gamma} a w^\prime \big( a\lambda^{\frac{\gamma-1}{\gamma}} \big) . $$
This, combined with the fact that $-w^\prime (y)\! = \! -\partial_y w_1(y)\! = \! w(y)^\gamma \! -\! 1$, implies 
$$\frac{\partial_\lambda w_1 (y)}{w(y)^\gamma \! -\! 1}= \frac{_1}{^\gamma} \frac{w(y)}{w(y)^\gamma \! -\! 1} -  \frac{_{\gamma-1}}{^\gamma} y \; ,$$
which implies (\ref{Lreecrit}) thanks to the second equality in (\ref{jointlawc}) in Proposition \ref{genjoint}. This completes the proof of  Proposition \ref{prop: lp}. 

\subsection{Explicit computation of $\bN_{{\rm nr}} [\Gamma]$ and $\bN_{{\rm nr}} [D]$.}
\label{Pfexpc}
We can deduce from Proposition \ref{prop: lp} explicit expressions for the first moment of $\Gamma$ and $D$ under $\nr$. 
\begin{prop}
 \label{expcht} We fix $\gam \ino (1, 2]$ and to simplify notation we set $\delta=1\! -\! \frac{1}{\gam}$. Then we get: 
\begin{align}\label{EGamma} 
\nr[\Gamma]&=\frac{2^{-1+\frac{2}{\gam}}\sqrt{\pi} }{\Gamma_{\! {\rm e}}(\tfrac{3}{2}\! -\! \tfrac{1}{\gam})}\int_0^1 dv \,v^{-\frac{1}{\gam}}\frac{(1-v^{\frac{1}{\gam}})(1-v^{\frac{\gam-1}{\gam}})}{(1-v)^2} \\ 
\label{EGamma'}
& = \frac{\sqrt{\pi} \,2^{-2\delta} }{\Gamma_{\! {\rm e}}(\tfrac{1}{2}+\delta)}\Big( \frac{1}{\delta}- \frac{2\delta}{1+ \delta} + 2\delta (1\! -\! \delta)  \!  \sum_{n\geq 1} \frac{2n+1+2\delta}{(n+\delta)(n+1+\delta)(n+2\delta)}\, \Big) \; .
\end{align}
\end{prop}
\noi
\textbf{Proof.} The scaling property \eqref{echtscal} entails that for any $\lambda \ino (0, \infty)$,
\begin{align}\label{one-hd}
\bN\big[\zeta e^{-\lambda\zeta}\Gamma\big] \! =\! c_\gamma\int_0^\infty \!\!\!\!\! dr\, r^{-1-\frac{1}{\gamma}} re^{-\lambda r} r^{\frac{\gamma -1}{\gamma}} 
\nr [\Gamma]=\lam^{-2+\frac{2}{\gam}}c_\gamma \Gamma_{\! {\rm e}} \big(2\! -\! \tfrac{2}{\gamma}\big) \nr [\Gamma \big] \; .
\end{align}
Recall from Proposition \ref{genjoint} that $\mathrm{L}_\lam(0, z) \! =\! \bN [e^{-\lambda\zeta}\indi_{\{\Gamma>z\}}]$. Thus, 
\begin{align}\label{other-hd}
\bN \big[\zeta e^{-\lambda\zeta}\Gamma\big] \! =\! \int_0^\infty \!\!\!\! dz \,\bN \big[ \zeta e^{-\lambda\zeta}\indi_{\{\Gamma>z\}}\big] \! =\! \int_0^\infty \!\!\!\! \!\!  dz\,  
\big(\! -\! \partial_\lam \mathrm{L}_\lam(0, z)\big).
\end{align}
Recall from (\ref{wdef1}) the definition of the function $w$. By \eqref{scalebL} and \eqref{htdianor} in Proposition \ref{prop: lp}, we get  
$$ -\partial_\lam\mathrm{L}_\lam(0, z)= \tfrac{1}{\gam} \lam^{\frac{1}{\gam} -1} \big( 1 \! -\! w\big(  z\lam^{\frac{\gam-1}{\gam}} \big) \big) - \tfrac{\gam-1}{\gam} zw^\prime \big(  z\lam^{\frac{\gam-1}{\gam}} \big) \; .$$
%
Recall that $1/(\gamma c_\gam)\! =\! \Gamma_{\! {\rm e}}(1\! -\! \frac{1}{\gamma})$. The previous equality, combined with \eqref{other-hd} and \eqref{one-hd} with $\lam\! =\! 1$, implies 
\begin{equation}\label{EGamma-1}
\nr [\Gamma] \! =\! \frac{\Gamma_{\! {\rm e}} \big(1\! -\! \tfrac{1}{\gam} \big)}{\Gamma_{\! {\rm e}} \big(2\! -\! \tfrac{2}{\gam} \big)} \int_0^\infty \!\!\!\! \!\!   dz\,  \big(1\! -\! w(z)\! -\! (\gam\! -\! 1)zw'(z)\big) \! = \! \frac{2^{-1+\frac{2}{\gam}}\sqrt{\pi}}{\Gamma_{\! {\rm e}}(\tfrac{3}{2} \! -\! \tfrac{1}{\gam})} \int_0^\infty \!\!\!\! \!\!   dz\,  \big(1\! -\! w(z)\! -\! (\gam\! -\! 1)zw'(z)\big) \; , 
\end{equation}
by the duplication formula for the gamma function: $\Gamma_{\! {\rm e}} \big(1\! -\! \tfrac{1}{\gam} \big)/ \Gamma_{\! {\rm e}} \big(2\! -\! \tfrac{2}{\gam} \big)\! = \! 2^{-1+\frac{2}{\gam}}\sqrt{\pi}/ \Gamma_{\! {\rm e}} \big(\frac{3}{2}\! -\! \tfrac{1}{\gam} \big)$. 
Recall that $w$ satisfies the integral equation \eqref{wdef1}. By the change of variable $y\! := \! w(z)$, we easily get 
\begin{equation}\label{EG-2}
\int_0^\infty \!\!\!\! \!\!   dz\,  \big(1\! -\! w(z)\! -\! (\gam\! -\! 1)zw'(z)\big)\!= \! \int_1^\infty\!\!\!\!\!  dy \, \Big( \frac{1\! -\! y}{y^\gam\! -\! 1}+(\gam\! -\! 1) \! \int_y^\infty \!\!\!\!  \frac{du}{u^\gam \! -\! 1}\Big).
\end{equation}
Note that $(1\! -\! y)/(y^\gam \! -\! 1) \! =\! \int_y^\infty du \, ((1\! -\! \gam)u^\gam \! -\! 1+\gam u^{\gam-1})/ (u^\gam \! -\! 1)^2$. Then, \eqref{EG-2} equals
\begin{align*}
&\int_1^\infty \!\!\!\!\! \! dy \! \int_y^\infty \!\!\!\! \!\! du\, \Big( \frac{(1\! -\! \gam)u^\gam \! -\! 1+\gam u^{\gam-1}}{(u^\gam \! -\! 1)^2}+\frac{\gam\! -\! 1}{u^\gam \! -\! 1} \Big)\! =\! \gam \! \int_1^\infty \!\!\!\!\! \! dy \! \int_y^\infty \!\!\!\! \!\! du\,\frac{u^{\gam-1} \! -\! 1}{(u^\gam \! -\! 1)^2} \\
&=\gam \! \int_1^\infty \!\!\!\!\! du \, \frac{(u\! -\! 1)(u^{\gam-1}\! -\! 1)}{(u^\gam\! -\! 1)^2} \! =\! \int_0^1 \!\!\!\! dv \,v^{-\frac{1}{\gam}}\frac{(1\! -\! v^{\frac{1}{\gam}})(1\! -\! v^{\frac{\gam-1}{\gam}})}{(1\! -\! v)^2} \; , 
\end{align*}
where we have used Fubini in the second equality and the change of variable $v\! =\! u^{-\gam}$ in the last one. 
By \eqref{EG-2} and \eqref{EGamma-1}, we get \eqref{EGamma}. 
We then use the expansion $(1\! -\! v)^{-2}\! =\! \sum_{n\ge 0}(n+1)v^n$ in \eqref{EGamma} to get \eqref{EGamma'} by straightforward computations.   \cqfd

\bigskip

We also get an explicit formula for $\nr[D]$ in terms of $\delta\! := \! 1\! -\! \frac{1}{\gam}$. The method is the same as in Proposition \ref{expcht} but computations are much longer; we skip the proof and we just state the result. 
\begin{prop}
\label{expcdm}  We fix $\gam \ino (1, 2]$ and to simplify notation we set $\delta=1\! -\! \frac{1}{\gam}$. Recall from (\ref{wdef1}) the definition of the function $w$. Then, 
\begin{equation}
\label{Ediam}
\nr[D]=\frac{2^{\frac{2}{\gam}}\sqrt{\pi} }{\Gamma_e(\tfrac{3}{2}-\tfrac{1}{\gam})}\int_1^\infty \!\!\!\! dx \, W(x) =
\frac{\sqrt{\pi}\, 2^{-2\delta}}{\Gamma_{\!{\rm e} }( \, \frac{1}{2}+ \delta)} \Big(\frac{2}{\delta} - 3 +
 \delta (A_1(\delta) + A_2 (\delta) + A_3 (\delta)) \big)  \Big) 
 \end{equation}
%
where for all $x \ino (0, \infty)$,  
$$ W(x) \! :=\!  2(\gam \! -\! 1)^2 x^{\gam -1} (x^\gam \! -\! 1)\Big(\!  \int_x^\infty\!  \!\!\!\!  \frac{du}{u^\gam\! -\! 1}\Big)^2 - \tfrac{(\gam-1)(2\gam +1)}{\gam}  
(x^\gam \! -\! 1)\!  \int_x^\infty \!\! \!\!\!  \frac{du}{u^\gam\! -\! 1}- \frac{x\! -\! 1}{x^\gam \! -\! 1} + \tfrac{1}{\gam} x \; , $$
%
where 
$$ A_1(\delta) = \tfrac{4(1-\delta)}{(1+\delta)^2} + \tfrac{3}{2+ \delta} \; , \quad A_2(\delta)=\! \! \!  \sum_{\substack{m,n\geq 0, \\ m+n\ge 3}} \! \! \tfrac{8(1-\delta) \delta}{(m+n-2+2\delta)(m+\delta) (n+\delta)}-\sum_{\substack{m,n\geq 0, \\ m+n\ge 2}} \! \! \tfrac{8(1-\delta) \delta}{(m+n-1+2\delta)(m+\delta) (n+\delta)} \; .$$
$$ \textrm{and} \quad A_3 (\delta)= \sum_{n\geq 2} \tfrac{4(1-\delta)}{(n-1+2\delta)(n-1+ \delta)} - \sum_{n\geq 3} \tfrac{4(1-\delta)(3-\delta)}{(n+\delta)(n-1+2\delta)(n-2+ 2\delta)} \; . $$

Note that $A_1(\delta) + A_2 (\delta) + A_3 (\delta) \! = \! \mathcal{O} (1)$ as $\delta \! \to \! 0$ (namely as $\gam \! \to \! 1$). 
\end{prop}
In the special case $\gam\! =\! 2$, \eqref{EGamma} implies $\nr [\Gamma] \! =\! \sqrt\pi$ and  (\ref{Ediam}) 
implies $\nr [D] \! =\! \frac{4}{3}\sqrt\pi$, that are known results which can be found in Szekeres \cite{Sz83} or Aldous \cite{aldcrt2}. As $\gam \! \to \! 1+$ (namely as $\delta \! \to \! 0+$), 
we use (\ref{EGamma'}), (\ref{Ediam}) and the well-known Taylor expansion of the gamma function: 
$$ \Gamma_{\! {\rm e}} \big(\tfrac{1}{2} + \delta \big) \! =\!  \sqrt{\pi}  - \delta \sqrt{\pi} (2\log 2 + \gamma_{{\rm e}} ) + \mathcal{O}(\delta^2) \; , $$
where $\gamma_{{\rm e}}$ stands for the Euler-Mascheroni constant, to get (\ref{siuccd}) in Remark \ref{expechtdm}.

\section{Proof of Theorems \ref{thm: ht_asy} and \ref{thm: dm_asy}. }
\label{pfasysec}
\subsection{Preliminary results.}
\label{preliasy}
In this section we prove several estimates that are used in the proof of Theorems \ref{thm: ht_asy} and \ref{thm: dm_asy}. We fix $\gamma \! \in \! (1, 2]$ and we take $\Psi (\lambda)\! = \! \lambda^\gamma$, $\lambda \! \in \! \bbR_+$.

\paragraph{Laplace transform.} 
We next introduce the following notation for the Laplace transform of Lebesgue integrable functions: 
for all measurable functions $f\! :\!  \bbR_+ \! \rightarrow \! \bbR$
such that there exists $\lam_0 \! \in \! \bbR_+$ satisfying
\begin{equation*}  
\int_0^\infty\!\!\! dx \,  e^{-\lam_0 x}|f(x)| < \infty , \quad \textrm{we set} \quad \cL_\lambda (f):= \int_0^\infty \!\!\! dx \, e^{-\lambda x}f(x), \quad \lambda \! \in \! [\lam_0, \infty) \; ,    
\end{equation*}
which is well-defined. 
We shall need the following lemma. 
\begin{lem}
\label{identiLa} Let $f, g_n, h_n\! : \! \bbR_+ \! \rightarrow \! \bbR_+$, $n\! \in \! \bbN$, 
be continuous and nonnegative functions. We set $f_n\! := \! g_n\! -\! h_n$. Let $(q_n)_{n\geq 0}$ be a real valued sequence. We make the following assumptions. 
\begin{align*}
\tag{$a$} \exists \lam_0 \in \bbR_+:  \quad \int_0^\infty \!\!\! dx \, e^{-\lambda_0 x}f(x)\! < \! \infty \quad \textrm{and} \quad \sum_{n\geq 0} \! |q_n|\!  \int_0^\infty \!\!\! \! \! dx \, e^{-\lambda_0 x} \big(g_n(x) +h_n(x) \big) \,  <  \infty \; .
\end{align*}
This makes sense of the sum $\sum_{n\geq 0} q_n \cL_\lambda (f_n)$ for all $\lambda \! \in \! [\lambda_0, \infty)$ and we  assume that 
\begin{align*}
\tag{$b$} \forall \lambda \! \in \! [\lambda_0, \infty) , \quad \cL_\lambda (f)=\sum_{n\geq 0} q_n \cL_\lambda (f_n)\; .
\end{align*}
We furthermore assume 
\begin{align*}
\tag{$c$} \forall x\in \bbR_+ , \quad\sum_{n\geq 0} |q_n| \, \big( \! \sup_{\; y\in [0, x] } \!\! g_n (y) \, + \sup_{\; y\in [0, x] }\!\!  h_n (y)\big) < \infty \; . 
\end{align*}
Then, 
$$ \forall x \in \bbR_+ , \quad f(x)= \sum_{n\geq 0} q_n f_n (x) \; , $$
where the sum in the right member makes sense thanks to ({\rm c}). 
\end{lem}
\noi
\textbf{Proof.} We denote by $(\cdot)^+$ and $(\cdot)^-$ resp.~the positive and negative part functions. Assumption ($c$) ensures that the following functions are well-defined for all $x\! \in \! \bbR_+$, continuous on $\bbR_+$ and nonnegative: 
$$ G \! :=\! f +   \sum_{n\geq 0} (q_n)^-\! g_n +(q_n)^+h_n \quad \textrm{and} \quad H\! := \sum_{n\geq 0} 
(q_n)^+g_n +(q_n)^-h_n . $$
Since the functions are nonnegative, for all $\lambda \! \in \! [\lambda_0, \infty)$, we get 
$$ \cL_\lambda (G)= \cL_\lambda (f)+ \sum_{n\geq 0} (q_n)^- \cL_\lambda  (g_n) +(q_n)^+ \cL_\lambda (h_n) 
 \quad \textrm{and} \quad  \cL_\lambda (H)= \sum_{n\geq 0} (q_n)^+ \cL_\lambda  (g_n) +(q_n)^-  \cL_\lambda (h_n) . $$
By Assumption ($a$), $\cL_\lambda (G)$ and $\cL_\lambda (H) \! $ are finite quantities for all $\lambda \! \geq \! \lambda_0$. 
Assumption ($b$) then entails that $\cL_\lambda (G)\! = \! \cL_\lambda (H)$, for all for all $\lambda \! \geq \! \lambda_0$: this implies that the Laplace transform of the finite Borel measures $e^{-\lambda_0 x} G(x) dx$ and 
$e^{-\lambda_0 x} H(x) dx$  are equal. Consequently, these measures are equal. Thus $G\! = \! H$ Lebesgue-almost everywhere. Since $G$ and $H$ are continuous, $G\! = \! H$ everywhere, which implies the desired result. \cqfd

\paragraph{Estimates for stable distributions.} Let $(\Omega, \cF, \bP)$ be an auxiliary space. Let $S\! : \! \Omega \! \rightarrow \!\bbR_+$ be a spectrally positive 
$\frac{\gamma -1}{\gamma}$-stable random variable such that 
\begin{equation}
\label{rebsgam}
 \forall \lambda \in \bbR+, \quad \bE \big[ e^{-\lambda S}\big]= \int_0^\infty \!\!\!\!\! dx \, s_\gamma (x) \exp(-\lambda x) = 
\exp \big(-\gamma \lambda^{\frac{\gamma -1}{\gamma}} \big), 
\end{equation}
where we recall from (\ref{denstable}) that $s_\gamma : \bbR_+\rightarrow \bbR_+$ is the continuous version of the density of the $\frac{\gamma -1}{\gamma}$-stable 
distribution. We recall here from Ibragimov \& Chernin \cite{IbrChe59} (see also 
Chambers, Mallows \&  Stuck \cite{ChMaSt76} formula (2.1) p.~341 or Zolotarev \cite{Zol}) the following representation of such a $\frac{\gamma -1}{\gamma}$-stable law: to that end, we first set 
\begin{equation}
\label{mvdef}
\forall v \! \in \! (-\pi, \pi), \quad m_\gamma(v)= \left( \frac{ \gamma \sin \big( \frac{_{\gamma -1}}{^\gamma} v\big) }{\sin v}\right)^{\!\! \gamma -1} 
\frac{\gamma \sin \big( \frac{_1}{^\gamma } v\big)}{\sin v} . 
\end{equation} 
Let $V, W$ be two independent random variables defined on $(\Omega, \cF, \bP)$ such that $V$ is uniformly 
distributed on $[0, \pi]$ and such that $W$ is exponentially distributed with mean $1$. Then, 
$$ S\quad  \overset{\textrm{(law)}}{=}\quad \left( \frac{m_\gamma(V)}{W}\right)^{\frac{1}{\gamma -1}} \; , $$
which easily implies that 
\begin{equation}
\label{repsgam}
\forall x \in (0, \infty) , \quad s_\gamma (x) = \frac{\gamma \! -\! 1}{\pi} x^{-\gamma}\!\!  \int_0^{\pi} \! \! \!\!  dv \, m_\gamma (v) \exp \big(\! - \! x^{-(\gamma -1)} m_\gamma (v) \big) \; .
\end{equation} 
Observe that $m_\gamma(-v)\! =\! m_\gamma (v)$ and $m_\gamma (0)\! =\! (\gamma \! -\! 1)^{\gamma  - 1}$. 
Moreover, the function $m_\gamma$ is increasing on $[0, \pi)$ 
and $ m_\gamma (v) / m_\gamma (0)=1+ \frac{{\gamma -1}}{2\gamma} v^2 + \mathcal{O}_\gamma (v^4) $.

As proved in Theorem 2.5.2 in Zolotarev \cite{Zol}, an extension of Laplace's method (proved in Zolotarev \cite{Zol}, Lemma 2.5.1, p.~97) yields the asymptotic expansion (\ref{eq: est}) that can be rewritten as follows: recall from (\ref{eq: est}) the definition of the sequence $(S_n)_{n\geq 1}$; then set 
\begin{equation}
\label{bbddeeff}
\forall x \! \in \! (0, \infty) \quad b(x) \! = \! \Big( \frac{\gamma \! -\! 1}{x}\Big)^{\gamma-1}
\quad \textrm{and} \quad S^*_n\! := \! 
\big( 2\pi \big( 1\! -\! \frac{{_1}}{{^\gamma}}\big)\big)^{\!\! -\frac{1}{2}} \! (\gamma \! -\! 1)^{\frac{\gamma+1}{2}-n(\gamma-1) } S_n, \quad n\geq 0\; , 
\end{equation}
where recall that $S_0\! = \! 1$. Then, 
for all positive integers $ N$, as $x\rightarrow 0$, we have 
\begin{equation}
\label{expasyreb}
s_\gamma (x)=\sum_{0\leq n<N} \! S^*_n \,\,  x^{n(\gamma-1) -\frac{\gamma+1}{2}} \, e^{-b(x)}  + \mathcal{O}_{\! N, \gamma} \big(x^{N(\gamma-1) -\frac{\gamma+1}{2}} e^{-b(x)} \big) \; .
\end{equation}
For all $a \! \in \! \bbR$, we next set 
\begin{equation}
\label{Jadef}
\forall x \in \bbR_+ \quad J_a (x):= \int_0^x \!\! dy \, y^a e^{-b(y)} \; .
\end{equation}
An integration by parts entails 
\begin{equation}
\label{partinteg}
 \forall a \in \bbR\backslash \{ -\gamma \}, \; \forall x \in \bbR_+, \quad J_a (x) =  (\gamma \! -\! 1)^{-\gamma}  x^{a+\gamma} e^{-b(x)} -(\gamma \! -\! 1)^{-\gamma}  ( a+\gamma)  J_{a+\gamma -1} (x) \; , 
\end{equation}
which proves that $J_a (x) = \mathcal{O}_\gamma (x^{a+\gamma} e^{-b(x)} )$ as $x \rightarrow 0$. This also entails the following lemma. 
\begin{lem} 
\label{asymJa} 
Let $\gamma \! \in \! (1, 2]$. Let $a\in \bbR$. We assume that $-(a\! +\! 1)/ (\gamma \!-\!1)$ is not a positive integer. 
Recall from (\ref{bbddeeff}) the definition of the function $b$ and from (\ref{Jadef}) the definition of the function $J_a$.  
Then, we set 
\begin{equation}
\label{cnagamdef}
\forall q\in \bbN\backslash \{ 0\}, \quad c_q (a, \gamma):=  (-1)^{q} (\gamma \! -\! 1)^{-(q+1)\gamma} \!\! \prod_{1\leq k\leq q} \!\!\! \big(a\! +\! 1\! +\!  k(\gamma \! -\! 1) \big) \; , 
\end{equation}
with the convention that $c_0(a, \gamma)\! =\!  (\gamma \! -\! 1)^{-\gamma}$. Then, for all positive integers $p$, 
\begin{equation}
\label{partintn}
J_a (x)= \!\! \! \sum_{0\leq q< p}\!\!\!  c_q(a, \gamma)\;   x^{a+\gamma + q(\gamma -1)} \, e^{-b(x)} \; \, + (\gamma \! -\! 1)^{\gamma}c_{ p} (a, \gamma) \, J_{a+ p(\gamma-1)} (x) \; .
\end{equation}
This implies that for all positive integers $p$, as $x\rightarrow 0$, 
\begin{equation}
\label{asycsq}
x^{-a-\gamma} e^{b(x)} J_a (x)=  \!\!  \sum_{0\leq q<p} \!\! c_q(a, \gamma) \; x^{q(\gamma -1)}  \; \, + \mathcal{O}_{p, a,\gamma} 
\big(x^{ p(\gamma -1)}   \big), 
\end{equation}
where $\mathcal{O}_{p, a,\gamma}$ depends on $p, a$ and $\gamma$. 
\end{lem}
\noi
\textbf{Proof.} (\ref{partintn}) follows from (\ref{partinteg}), by induction. Since $J_{a+ p(\gamma-1)} (x)= \mathcal{O}_\gamma \big(x^{a+ p(\gamma-1)+ \gamma} e^{-b(x)}\big)$, 
(\ref{asycsq}) is an immediate consequence of (\ref{partintn}). \cqfd

\bigskip

We next prove the following lemma. 
\begin{lem}
\label{erderivss} Let $\gamma \! \in \! (1, 2]$. Recall from (\ref{denstable}) (or from (\ref{rebsgam})) the definition of the density $s_\gamma$. Recall from (\ref{mvdef}) the definition of $m_\gamma$. We set for all $x\! \in \! \bbR_+$, 
\begin{align}
\label{sigplms}
& \sigma^+ (x)\! := \! \frac{_{(\gamma-1)^2}}{^\pi} x^{-2\gamma} \!\! \! \int_0^\pi \!\!\! dv\,  m_\gamma (v)^2  
e^{- x^{-(\gamma -1)} m_\gamma (v)}  \\
\quad \textrm{and} \qquad & \sigma^-(x)\! : =  \! \gamma x^{-1} s_\gamma (x) =   \frac{_{\gamma (\gamma  - 1)}}{^\pi}
x^{-\gamma -1} \!\! \! \int_0^\pi  \!\!\!\! dv \, m_\gamma (v) 
e^{- x^{-(\gamma -1)} m_\gamma (v)} \; . \nonumber
\end{align}
Then, the following holds true. 
\begin{itemize}
\item[(i)] $\sigma^+$ and $\sigma^-$ are well-defined on $\bbR_+$, the function $s_\gamma$ is differentiable on $\bbR_+$ and $s^\prime_\gamma= \sigma^+-\sigma^-$. Moreover, $\sigma^+$, $\sigma^+$ are continuous, nonnegative, 
Lebesgue integrable and for all $\lambda \! \in \! \bbR_+$, 
\begin{equation}
\label{Laplsplms}
\cL_\lambda (\sigma^+)= \lambda e^{-\gamma \lambda^{\frac{\gamma -1}{\gamma}}} \! +\gamma\!  \int_\lambda^\infty 
 \!\!\! \!\! d\mu \, e^{-\gamma \mu^{\frac{\gamma-1}{\gamma}}}  \quad \textrm{and} \quad 
 \cL_\lambda (\sigma^-)= \gamma \!  \int_\lambda^\infty 
 \!\!\! \!\! d\mu \, e^{-\gamma \mu^{\frac{\gamma-1}{\gamma}}} , 
\end{equation} 
which implies 
\begin{equation}
\label{Laplsprim}
\int_0^\infty\!\!\! dx \,  |s^\prime_\gamma (x)| < \infty  \qquad \textrm{and} \qquad \cL_\lambda (s^\prime_\gamma)= \lambda e^{-\gamma \lambda^{\frac{\gamma -1}{\gamma}}}, 
\quad \lambda \! \in \! \bbR_+.   
\end{equation}
\item[(ii)] There exist $A, x_0 \! \in (0, \infty)$ such that 
\begin{equation}
\label{domlocspm}
\forall x \in [0, x_0], \qquad \sigma^{+} (x) \; \, \textrm{and} \; \, \sigma^-(x) \leq A x^{-\frac{3\gamma +1}{2}} e^{-b(x)} \; , 
\end{equation}
where we recall from (\ref{bbddeeff}) that $b(x)= \big( (\gamma \! -\! 1)/x\big)^{\gamma -1}$.  
\item[(iii)] We define the real valued sequence $(T^*_n)_{n\geq 0}$ by 
\begin{equation}
\label{TdedS}
 T^*_0 := (\gamma \! -\! 1)^\gamma S^*_0 \quad \textrm{and} \quad \forall n\geq 1, \quad T^*_n:=  (\gamma \! -\! 1)^\gamma 
S^*_n+ \big( n(\gamma \! -\! 1)- \frac{_{3\gamma -1}}{^2} \big) S^*_{n-1} \; .
\end{equation}
Then, for all positive integer $N$, as $x\rightarrow 0$, we have 
\begin{equation}
\label{erasprim}
s^\prime_\gamma (x)= \sum_{0\leq n<N} T^*_n \, x^{n(\gamma-1) -\frac{3\gamma+1}{2}} \, e^{-b(x)}  + \mathcal{O}_{N, \gamma} \big(x^{N(\gamma-1) -\frac{3\gamma+1}{2}} e^{-b(x)} \big) \; .
\end{equation}
\end{itemize}
 \end{lem}
\noi
\textbf{Proof.} We easily deduce from  (\ref{repsgam}), that $s_\gamma$ is differentiable on $\bbR_+$ and that $s^\prime_\gamma\! = \! \sigma^+\! -\! \sigma^-$. 
Using Fubini-Tonnelli and the change of variable $y\! = \! x^{-(\gamma -1)} m_\gamma (v)$, for fixed $v$, we get 
$$ \int_0^\infty \!\!\! dx \, \sigma^+(x) =   \int_0^\infty \!\!\! dx \, \sigma^-(x) = \frac{_\gamma}{^\pi} \Gamma_{\!e} \big(\frac{_\gamma}{^{\gamma -1}}\big) 
 \!\! \int_0^\pi \!\!\! dv\,  m_\gamma (v)^{-\frac{1}{\gamma -1}} < \infty  \; , $$
since $m_\gamma (v) \! \geq \! m_\gamma (0)\! >\! 0$ on $[0, \pi)$ and $\lim_{v\rightarrow \pi} m_\gamma (v)\! = \! \infty$; here, $\Gamma_{\!e}$ stands for Euler's gamma function. Thus, $\int_0^\infty
 dx \,  |s^\prime_\gamma (x)| \! < \! \infty$ and $\lambda \! \in \! \bbR_+ \mapsto \cL_\lambda (s^\prime_\gamma)$ is well-defined. Moreover, by Fubini,  
$$ \cL_\lambda (s^\prime_\gamma)= \int_0^\infty \!\!\! dx \, s^\prime_\gamma (x)\int_x^\infty \!\!\! dy \, \lambda e^{-\lambda y}= \lambda \!\! \int_0^\infty \!\! \! dy \, e^{-\lambda y} \int_0^y \!\! dx \, s^\prime_\gamma (x) = \lambda \cL_\lambda (s_\gamma)\; , $$   
which completes the proof of (\ref{Laplsprim}). Next, by Fubini-Tonnelli, we get 
\begin{equation}
\label{Laplhpl}
\int_0^\infty \!\!\! dx \, e^{-\lambda x}x^{-1} s_\gamma (x) = 
\int_0^\infty \!\!\! dx \, s_\gamma (x) \int_\lambda^\infty \!\!\! d\mu \, e^{-\mu x} =  \int_\lambda^\infty \!\!\! d\mu \, e^{-\gamma \mu^{\frac{\gamma -1}{\gamma}}} \; .
\end{equation}
which implies that $\cL_\lambda (\sigma^-) =\gamma  \int_\lambda^\infty 
d\mu \, e^{-\gamma \mu^{\frac{\gamma -1}{\gamma}}}$, since $\sigma^- (x)\! =\!  \gamma x^{-1} s_\gamma (x)$. 
This, combined with (\ref{Laplsprim}) entails (\ref{Laplsplms}), which completes the proof of $(i)$. 

\medskip

Laplace's method easily implies that there exists $c_+, c_- \in (0, \infty)$ such that 
$$ \sigma^+ (x) \underset{x\rightarrow 0}{\sim} c_+ x^{-\frac{3\gamma +1}{2}} e^{-b(x)} \quad \textrm{and} \quad
\sigma^- (x) \underset{x\rightarrow 0}{\sim} c_- x^{-\frac{\gamma +3}{2}} e^{-b(x)} \; , $$
which easily entails (\ref{domlocspm}) and which completes the proof of $(ii)$. 

More generally, the asymptotic expansion (\ref{eq: est}) of $s_\gamma$ is derived from (\ref{repsgam}) by an extension of Laplace's method proved in Zolotarev \cite{Zol}, Lemma 2.5.1, p.~97. When this method is applied to $\sigma^+$ and $\sigma^-$, one shows that $\sigma^+$ and $\sigma^-$ have an asymptotic expansion whose general term is $x^{n(\gamma -1) -\frac{3\gamma +1}{2}} e^{-b(x)}$. Thus, there exists a sequence $(T^*_n)_{n\geq 0}$ such that (\ref{erasprim}) holds true. It remains to prove (\ref{TdedS}). To that end, for any $n \! \in \! \bbN$, we set $a_n \! :=\! n(\gamma \! -\! 1) -\frac{3\gamma+1}{2}$. By Lemma \ref{asymJa}  we then get  
\begin{eqnarray*}
s_\gamma (x)& =& \sum_{0\leq n< N} T^*_n J_{a_n} (x) + \mathcal{O}_{N, \gamma} \big( J_{a_N}  (x)\big) \\
& =& \sum_{0\leq n <N} \!\!\!  \!\!\!\sum_{\quad 0\leq q <N-n}  \!\!\! \!\!\! \!\!\! T^*_n c_q (a_n , \gamma) x^{a_n +\gamma +q (\gamma -1)} e^{-b(x)} + 
\mathcal{O}_{N, \gamma} \big( x^{a_N +\gamma} e^{-b(x)} \big) \\
& =& \sum_{0\leq n\leq p <N} T^*_n c_{p-n} (a_n , \gamma) x^{p(\gamma-1)-\frac{\gamma +1}{2} } e^{-b(x)} +  
\mathcal{O}_{N, \gamma} \big( x^{N(\gamma-1)-\frac{\gamma +1}{2}} e^{-b(x)} \big) \; , 
\end{eqnarray*} 
which implies that $S^*_p \! =\!  \sum_{0\leq n\leq p} T^*_n c_{p-n} (a_n, \gamma)$, for all $p\! \in \! \bbN$. Then by (\ref{cnagamdef}), observe that 
\begin{eqnarray*}
 S_p^* &= & c_0 (a_p, \gamma) T^*_p+ \sum_{0\leq n\leq p-1} T^*_n c_{p-n} (a_n, \gamma)\\
 & =&   (\gamma\! -\! 1)^{-\gamma} T^*_p - (\gamma\! -\! 1)^{-\gamma} \big( p(\gamma\! -\! 1) -\frac{_{3\gamma -1}}{^2} \big)\!\!\!\! \sum_{0\leq n\leq p-1}\!\!\!  T^*_n c_{p-1-n} (a_n, \gamma) \; ,  
 \end{eqnarray*} 
which implies (\ref{TdedS}). This completes the proof of the lemma. \cqfd 

\bigskip

\noi
\textbf{Proof of Proposition \ref{derivss}.} Lemma \ref{erderivss} easily entails Proposition \ref{derivss}: indeed (\ref{Laplsprim}) entails (\ref{eLaplsprim}). We then set 
$$ \forall n \! \in \! \bbN, \quad T_n:=  (\gamma\! -\! 1)^{n (\gamma -1)} T^*_n /T^*_0 \; , $$
and we easily check that (\ref{TdedS}) entails (\ref{eTdedS}) and that (\ref{erasprim}) implies (\ref{asprim}). \cqfd  

\bigskip

We next introduce another function used in the asymptotic expansion of the height and the diameter of normalized stable tree.  
\begin{lem}
\label{thetpro} Let $\gamma \! \in \! (1, 2]$. Recall from (\ref{denstable}) (or from (\ref{rebsgam})) the definition of $s_\gamma$. We then introduce the following functions: for all $x\! \in \! \bbR_+$, 
 \begin{equation}
\label{thetadef}
h^{\! +} (x)\! = \! (\gamma \! -\! 1) \, x^{-1} \! s_\gamma (x), \quad h^{\! -}(x) 
\! = \! \frac{_{\gamma \! -\! 1}}{^\gamma} x^{-1-\frac{1}{\gamma}}\!\!  \int_0^x \!\!\! dy \, y^{\frac{1}{\gamma} -1} 
\! s_\gamma (y) \quad \textrm{and} \quad \theta(x) \! =\!  h^{\! +}  (x)\! - \! h^{\! -} (x) . 
\end{equation}
Then, the following holds true. 
\begin{itemize}
\item[(i)] $h^+, h^-$ and $\theta$ are well-defined and continuous, 
$h^+$ and $h^-$ are nonnegative and Lebesgue integrable, and 
for all $\lambda \! \in \! \bbR_+$, we have  
\begin{equation}
\label{Laplhplms}
\cL_\lambda(h^{\! +})=  (\gamma\! -\! 1) \!\! \int_\lambda^\infty \!\!\! d\mu \, e^{-\gamma \mu^{\frac{\gamma -1}{\gamma}}} 
\quad \textrm{and} \quad  \cL_\lambda(h^{\! -})= \cL_\lambda(h^{\! +})-\lambda^{\frac{1}{\gamma}} e^{-\gamma \lambda^{\frac{\gamma -1}{\gamma}}},
\end{equation}
which implies
\begin{equation}
\label{Lapltheta}
\int_0^\infty \!\!\! dx \,  |\theta (x)|  <  \infty \quad \textrm{and} \quad \cL_\lambda (\theta)= \lambda^{\frac{1}{\gamma}} e^{-\gamma \lambda^{\frac{\gamma -1}{\gamma}}}, 
\quad \lambda \! \in \! \bbR_+.   
\end{equation}
\item[(ii)] There exist $A, x_0 \! \in (0, \infty)$ such that 
\begin{equation}
\label{domlochpm}
\forall x \in [0, x_0], \qquad h^{\! +} (x)\quad \textrm{and} \quad   h^{\! -} (x) \; \,  \leq \;  Ax^{-\frac{\gamma +3}{2}} e^{-b(x)} \; , 
\end{equation}
where we recall from (\ref{bbddeeff}) that $b(x)= \big( (\gamma \! -\! 1)/x\big)^{\gamma -1}$.  

\item[(iii)]  Let $(V^*_n)_{n\geq 0}$ be a sequence of real numbers recursively defined by $V^*_0\! = \! (\gamma \! -\! 1) S^*_0 $ and for all $n \! \in \! \bbN$, 
\begin{equation}
\label{recVp}
 (\gamma \! -\! 1)^{\gamma -1} V^*_{n+1}=  (\gamma \! -\! 1)^{\gamma } S^*_{n+1} +  (\gamma \! -\! 1) \big(n \!-\! \frac{_1}{^2} \! -\! \frac{_1}{^{\gamma -1}}  \big) S^*_n - 
\big( n  \!-\! \frac{_1}{^2} \! -\! \frac{_1}{^{\gamma }}  \big) V^*_n \; . 
\end{equation}
Then for all positive integers $N$, as $x\rightarrow 0$, we get 
\begin{equation}
\label{theexpin}
\theta (x) = \sum_{0\leq n<N}  V^*_n \; x^{ n(\gamma-1)-\frac{\gamma+3}{2}} e^{-b(x)} + \mathcal{O}_{N, \gamma} \big( x^{ N(\gamma-1)-\frac{\gamma+3}{2} } e^{-b(x)}\big) . 
\end{equation}
\end{itemize} 
\end{lem}
\noi
\textbf{Proof.} The fact that $h^+$ and $h^-$ are well-defined is an easy consequence of the asymptotic expansion (\ref{expasyreb}) of $s_\gamma$ and observe that $h^+$, $h^-$ can be continuously extended by the value $0$ at 
$x\! = \!0$.  
Let $\lambda \! \in \! \bbR_+$; by (\ref{Laplhpl}) we get $\cL_\lambda(h^+)\! = \!  (\gamma\! -\! 1) 
 \int_\lambda^\infty d\mu \, \exp (-\gamma \mu^{\frac{\gamma -1}{\gamma}})$. 
Thus when $\lambda \! = \! 0$, we get 
$$ \int_0^\infty \!\!\! dx \, h^+(x)= \cL_0(h^+) =  (\gamma\! -\! 1) \!\! \int_0^\infty \!\!\! d\mu \, e^{-\gamma \mu^{\frac{\gamma -1}{\gamma}}} = \gamma^{-\frac{1}{\gamma -1}} \Gamma_{\! e} \big( \frac{_{\gamma}}{^{\gamma-1}}\big), $$
by an easy change of variable; here $\Gamma_{\! e}$ stands for Euler's Gamma function. By Fubini-Tonnelli and several linear changes of variable, we get 
\begin{eqnarray*}
\cL_\lambda (h^-) &=& \frac{_{\gamma \! -\! 1}}{^\gamma} \!\!\!  
\int_0^\infty \!\!\! dy \, y^{\frac{1}{\gamma} -1} s_\gamma (y) \!\! \int_y^\infty \!\!\! dx \,  x^{-1-\frac{1}{\gamma}} e^{-\lambda x} 
 = \frac{_{\gamma \! -\! 1}}{^\gamma} \lambda^{\frac{1}{\gamma}}
\!\!\! \int_0^\infty \!\!\! dy \, y^{\frac{1}{\gamma} -1} s_\gamma (y) \!\! \int_{\lambda y}^\infty \!\!\! d\mu \,  \mu^{-1-\frac{1}{\gamma}} e^{-\mu} \\
& =&\frac{_{\gamma \! -\! 1}}{^\gamma} \lambda^{\frac{1}{\gamma}}
\!\!\! \int_0^\infty \!\!\! dy \, y^{-1} s_\gamma (y) \!\! \int_{\lambda }^\infty \!\!\! d\nu \,  \nu^{-1-\frac{1}{\gamma}} e^{-\nu y} = \frac{_{\gamma \! -\! 1}}{^\gamma} \lambda^{\frac{1}{\gamma}}\!\!  \int_{\lambda }^\infty \!\!\! d\nu \,  \nu^{-1-\frac{1}{\gamma}} 
\!\!\! \int_0^\infty \!\!\! dy \, y^{-1} s_\gamma (y) e^{-\nu y} \\
& =& \frac{_{\gamma \! -\! 1}}{^\gamma} \lambda^{\frac{1}{\gamma}} \!\!  \int_{\lambda }^\infty \!\!\! d\nu \,  \nu^{-1-\frac{1}{\gamma}} 
\!\!\! \int_\nu^\infty \!\!\! d\mu  \, e^{-\gamma \mu^{\frac{\gamma -1}{\gamma}}} = (\gamma \! -\! 1)  \lambda^{\frac{1}{\gamma}} \!\! \int_\lambda^\infty \!\!\! d\mu  \, e^{-\gamma \mu^{\frac{\gamma -1}{\gamma}}} \big( \lambda^{-\frac{1}{\gamma}} \! -\! \mu^{-\frac{1}{\gamma}} \big) \\
& =& (\gamma\! -\! 1) \!\! \int_\lambda^\infty \!\!\! d\mu \, e^{-\gamma \mu^{\frac{\gamma -1}{\gamma}}} -(\gamma \! -\! 1)  \lambda^{\frac{1}{\gamma}} \!\! \int_\lambda^\infty \!\!\! d\mu  \, \mu^{-\frac{1}{\gamma}}e^{-\gamma \mu^{\frac{\gamma -1}{\gamma}}} \\
& =&  (\gamma\! -\! 1) \!\! \int_\lambda^\infty \!\!\! d\mu \, e^{-\gamma \mu^{\frac{\gamma -1}{\gamma}}}- \lambda^{\frac{1}{\gamma}} e^{-\gamma \lambda^{\frac{\gamma -1}{\gamma}}}. 
\end{eqnarray*}
Here we use (\ref{Laplhpl}) in the third line. When $\lambda \! = \! 0$, this proves that 
$$\int_0^\infty dx \, h^-(x)\! = \!  \gamma^{-\frac{1}{\gamma -1}} \Gamma_{\! e} \big( \frac{_{\gamma}}{^{\gamma-1}}\big) \; .$$
Thus, $\int_0^\infty dx \, |\theta (x)| \! < \! \infty$. It also 
implies (\ref{Lapltheta}) thanks to (\ref{Laplhpl}), which completes the proof of $(i)$. 

\medskip

We then prove ($ii$) and 
($iii$). To that end, we first observe that (\ref{expasyreb}) implies that $x^{-1} s_\gamma (x) \sim S^*_0 x^{-\frac{\gamma +3}{2}} e^{-b(x)}$ as $x\! \rightarrow \! 0$, which immediately entails (\ref{domlochpm}) for $h^+$.

We next find the asymptotic expansion of $h^-$ thanks to that of $s_\gamma$ and thanks to Lemma \ref{asymJa}. We first 
set $\alpha_n= \frac{1}{\gamma} -\frac{\gamma +3}{2} + n(\gamma \! -\! 1)$. From (\ref{expasyreb}) and Lemma  \ref{asymJa}, for all positive integer $N$, as $x\rightarrow 0$, we get 
\begin{eqnarray*}
h^-(x) &=& \sum_{0\leq n<N}  \frac{_{\gamma - 1}}{^\gamma} S^*_n \, x^{-1-\frac{1}{\gamma}}J_{\alpha_n} (x)+ \mathcal{O}_{N, \gamma} \big( x^{-1-\frac{1}{\gamma}}J_{\alpha_N} (x)\big)\\
& =&  \sum_{0\leq n<N}\!\!\! \!\!\!  \sum_{\quad 0\leq q <N-n}\!\!\! \!\!\! \!\!\! \frac{_{\gamma - 1}}{^\gamma} S^*_n c_q (\alpha_n, \gamma) x^{\alpha_n + \gamma -1-\frac{1}{\gamma}+ q(\gamma-1) } e^{-b(x)} + \mathcal{O}_{N, \gamma} \big( x^{\alpha_N +\gamma -1-\frac{1}{\gamma}} e^{-b(x)}\big) \\
& =&  \sum_{0\leq n<N}\!\!\! \!\!\!  \sum_{\quad 0\leq q <N-n}\!\!\! \!\!\! \!\!\! \frac{_{\gamma  - 1}}{^\gamma} S^*_n c_q (\alpha_n, \gamma) x^{ (n+q+1)(\gamma-1)-\frac{\gamma+3}{2} } e^{-b(x)} + \mathcal{O}_{N, \gamma} \big( x^{ (N+1)(\gamma-1)-\frac{\gamma+3}{2} } e^{-b(x)}\big) \\
& =& \sum_{0\leq p\leq N}
U_p \,  x^{ p(\gamma-1)-\frac{\gamma+3}{2}} e^{-b(x)} + \mathcal{O}_{N, \gamma} \big( x^{ (N+1)(\gamma-1)-\frac{\gamma+3}{2} } e^{-b(x)}\big) . 
\end{eqnarray*}
where the sequence $(U_p)_{p\geq 0}$ is given by 
$$ U_0= 0, \quad \textrm{and} \quad U_p = \sum_{0\leq n\leq p-1} \frac{_{\gamma \! -\! 1}}{^\gamma} S^*_n c_{p-1-n} (\alpha_n, \gamma)  , \quad p\geq 1 . $$
Observe that it implies (\ref{domlochpm}) for $h^-$, which completes the proof of $(ii)$. 
We next prove $(iii)$: to that end observe that  by (\ref{cnagamdef}), $c_{p-n} (\alpha_n, \gamma) \! =\! - (\gamma \! -\! 1)^{-\gamma} \big( \frac{1}{\gamma} \! -\! \frac{\gamma+1}{2} +p(\gamma \! -\! 1) \big) c_{p-1-n} (\alpha_n, \gamma)$. Thus we get  
\begin{eqnarray}
\label{recUp}
U_{p+1}& = & \sum_{0\leq n\leq p} \frac{_{\gamma \! -\! 1}}{^\gamma} S^*_n c_{p-n} (\alpha_n, \gamma) \; = \frac{_{\gamma \! -\! 1}}{^\gamma} S^*_p c_{0} (\alpha_p, \gamma) \; + \!\!\! \sum_{0\leq n\leq p-1}\!\!\!  \frac{_{\gamma \! -\! 1}}{^\gamma} S^*_n c_{p-n} (\alpha_n, \gamma) \nonumber  \\
& =&  \frac{_{1}}{^\gamma} (\gamma \! -\! 1)^{-(\gamma - 1)} S^*_p -  (\gamma \! -\! 1)^{-\gamma} \big( \frac{_1}{^\gamma} \! -\! \frac{_{\gamma+1}}{^2} +p(\gamma \! -\! 1) \big)\!\!\!  \sum_{0\leq n\leq p-1}\!\!\!  \frac{_{\gamma \! -\! 1}}{^\gamma} S^*_n c_{p-1-n} (\alpha_n, \gamma)  \nonumber  \\
& =& \frac{_{1}}{^\gamma} (\gamma \! -\! 1)^{-(\gamma - 1)} S^*_p -  (\gamma \! -\! 1)^{-\gamma} \big( \frac{_1}{^\gamma} \! -\! \frac{_{\gamma+1}}{^2} +p(\gamma \! -\! 1) \big) U_p  \nonumber  \\
& =&  (\gamma \! -\! 1)^{-(\gamma - 1)} \Big( \frac{_{1}}{^\gamma}  S^*_p - \big( p\! -\! \frac{_1}{^2} \! -\! \frac{_1}{^\gamma}\big) U_p  \Big) . 
\end{eqnarray} 
We then set $V^*_p= (\gamma \! -\! 1) S^*_p -U_p$ for all $p\! \in \! \bbN$, so that for all positive integer $N$, as $x\rightarrow 0$, (\ref{theexpin}) holds true. Moreover, (\ref{recUp}) easily entails that 
$(V^*_p)_{p\geq 0}$ satisfies (\ref{recVp}), which completes the proof of the lemma. \cqfd

\bigskip

\noi
\textbf{Proof of Proposition \ref{thetproo}.} Lemma \ref{thetpro} easily entails Proposition \ref{thetproo}. Indeed, (\ref{Lapltheta})
implies (\ref{eLapltheta}). We set 
$$ \forall n\in \bbN, \quad V_n= (\gamma \! -\! 1)^{n(\gamma -1)} V^*_n/V^*_0 \; .$$
Then, (\ref{recVp}) entails (\ref{erecVp}) and (\ref{theexpin}) implies (\ref{eq: esth}), which completes the proof of Proposition \ref{thetproo}.
\cqfd

\begin{lem}
\label{quintlaLa} There exist $\lambda_0, A\! \in \! (0, \infty)$ such that 
$$ \forall \lambda \in [\lambda_0, \infty), \quad   \!\! \int_\lambda^\infty \!\!\! d\mu \, e^{-\gamma \mu^{\frac{\gamma -1}{\gamma}}} \leq A \lambda^{\frac{1}{\gamma}} e^{-\gamma \lambda^{\frac{\gamma -1}{\gamma}}}\; .$$
\end{lem}
\noi
\textbf{Proof.} Integration by part implies 
$$  (\gamma \! - \! 1) \!\! \int_\lambda^\infty  \!\!\! \!\!\! d\mu \, e^{-\gamma \mu^{\frac{\gamma -1}{\gamma}}} \! = \lambda^{\frac{1}{\gamma}} e^{-\gamma \lambda^{\frac{\gamma -1}{\gamma}}} \!\! +     \frac{_1}{^\gamma}  \! \int_\lambda^\infty \!\!\!\!\!\! d\mu \, \mu^{-\frac{\gamma-1}{\gamma}}e^{-\gamma \mu^{\frac{\gamma -1}{\gamma}}} \leq \lambda^{\frac{1}{\gamma}} e^{-\gamma \lambda^{\frac{\gamma -1}{\gamma}}} \!\! +  \frac{_1}{^\gamma} 
\lambda^{-\frac{\gamma-1}{\gamma}} \! \! \int_\lambda^\infty \!\!\!\!\!\! d\mu \, e^{-\gamma \mu^{\frac{\gamma -1}{\gamma}}},  
$$
which immediately entails the lemma. \cqfd

\paragraph{Asymptotic expansion of $w \!-\!1$.}
Recall from (\ref{wdef1}) the definition of $w$. We next introduce 
\begin{equation}
\label{phideff}
\forall y \! \in \! (0, \infty), \quad \phi(y):= w(y)\! -\! 1, \quad \textrm{that satisfies} \quad \int_{\phi(y)}^\infty \!\! \frac{du}{(u+1)^\gamma \! -\! 1 }= y \; , 
\end{equation}
by (\ref{wdef1}). We easily see that $\lim_{y\rightarrow \infty} \phi(y)\! = \! 0$ and $\lim_{y\rightarrow 0} \phi(y)\! = \! \infty$ and that $\phi$ is a $C^\infty$ decreasing function. The following lemma asserts that $\phi$ decreases exponentially fast as $y\! \rightarrow \! \infty$. 
\begin{lem}
\label{vraiexphi}Let $\gamma \! \in \! (1, 2]$. Let $\Psi (\lambda)\! = \! \lambda^\gamma$, $\lambda \! \in \! \bbR_+$. 
Recall from (\ref{phideff}) the definition of $\phi$. We set 
\begin{equation}
\label{defGG}
y_0:= \int_{1}^\infty \!\! \frac{du}{(u+1)^\gamma \! -\! 1 } \quad \textrm{and} \quad \forall y \! \in \! [-1, \infty) , \quad G(y) := \int_0^y \! \frac{du}{u}\, \frac{(u+1)^\gam\! -\! 1\! -\! \gam u}{(u+1)^\gam \! -\! 1}\; .
\end{equation}
Then, 
\begin{equation}
\label{zahcgdfd}
\forall y \! \in [-1, 1], \quad \exp (G(y))= 1+\sum_{n\geq 1} A_n y^n \quad  \textrm{and} \quad  1+\sum_{n\geq 1} |A_n| < e^{\gamma -1} . 
\end{equation}
Moreover, for $y \! \in \! [y_0, \infty)$, 
\begin{equation} 
\label{cphi1} 
e^{\gamma y-C_0 }\phi (y)= \exp \big( G(\phi(y)) \big)= 1+ \sum_{n\geq 1} A_n \phi(y)^n \; , 
\end{equation}
where $C_0$ is given by (\ref{eq: Cgam}). Then, there exists a real valued sequence $(\beta_n)_{n\geq 1}$ and  $y_1 \! \in \! [y_0, \infty)$ such that 
\begin{equation}
\label{echtexphi}
\sum_{n\geq 1} |\beta_n| e^{-\gamma ny_1} < \infty \quad \textrm{and} \quad 
\forall y \in [y_1, \infty), \qquad \phi (y) =\sum_{n\geq 1} \beta_n e^{-\gamma n y} \; . 
\end{equation}
Here $\beta_1= e^{ C_0}$ and $\beta_2= \frac{\gamma -1}{2} e^{2C_0}$. 
\end{lem}
\noi
\textbf{Proof.} For all $y \! \in \! (0, \infty)$, we first set $F(y)\! :=\! \int_y^\infty \frac{du}{(u+1)^\gam\! -\! 1}$ that is such that $F(\phi (y))\! = \! y$. 
Observe that 
$$F(y)\! =\! \int_1^\infty \!\!\!\!\! \frac{du}{(u+1)^\gam \! -\! 1}+ \frac{_1}{^\gamma} \! \int_y^1\! \frac{du}{u}-\! 
\frac{_1}{^\gamma} \! \int_0^1\! \frac{du}{u}\, \frac{(u+1)^\gam \! -\! 1\! -\! \gam u}{(u+1)^\gam \! -\! 1}
+  \frac{_1}{^\gamma} \! \int_0^y \! \frac{du}{u}\, \frac{(u+1)^\gam\! -\! 1\! -\! \gam u}{(u+1)^\gam \! -\! 1} ,$$
which makes sense since 
$\frac{1}{u}\, \frac{(u+1)^\gam-1-\gam u}{(u+1)^\gam-1}\to \frac{\gam-1}{2}$ as $u\to 0+$. We then set 
$$ C_0:= \gamma \! \int_1^\infty \!\!\!\!\! \frac{du}{(u+1)^\gam \! -\! 1}- \int_0^1\! \frac{du}{u}\, \frac{(u+1)^\gam \! -\! 1\! -\! \gam u}{(u+1)^\gam \! -\! 1} $$
and we get 
$$ \forall y \in (0, \infty) , \quad \gamma F(y)= C_0- \log y+ G(y), \quad \textrm{where} \quad G(y)\! :=\!\!   \int_0^y \! \frac{du}{u}\, \frac{(u+1)^\gam\! -\! 1\! -\! \gam u}{(u+1)^\gam \! -\! 1}. $$
Since $F(\phi(y))\! =\!  y$, this implies 
\begin{equation}
\label{logcL}
  \forall y \in (0, \infty) , \quad \log \phi(y)= C_0-\gamma y +  G(\phi(y))\; .
\end{equation}  
Let us show that $G(y)$ (and therefore $\exp(G(y))$) is analytic in a neighborhood of $0$. We set
$$
a_n=\igam\binom{\gam}{n\! +\! 1}=\frac{(-1)^{n-1}}{(n+1)!}\prod_{k=1}^n|k-\gam|= \frac{(\gamma \! -\! 1) (-1)^{n-1}}{n(n+1)} \prod_{k=1}^{n-1} \Big( 1-\frac{\gam \! -\! 1}{k} \Big)  , \quad n\ge 1.
$$
We observe that  $|a_n|<\frac{\gam-1}{n(n+1)}$. Then for all $u \!\in \! [-1, 1]$, we set  
$$T(u):=\sum_{n\geq 1} |a_n| u^n\quad \textrm{and} \quad S(u) \! :=\! \frac{(1+u)^\gam \! -\! 1\! -\! \gam u}{\gamma u}= \sum_{n\ge 1}a_n u^n=-T(-u)\; , $$
since $(-1)^{n-1} a_n= |a_n|$. The power 
series $T$ and $S$ are absolutely convergent for $|u| \! \leq\! 1$. Moreover, $|S(u)|\leq T(|u|)\leq T(1)\! =\! -S(-1)= \frac{\gamma -1}{\gamma}\! \le  \! 1$.  
Thus, for all $u \! \in \! [-1, 1]$,
$$
\frac{(1+u)^\gam-1-\gam u}{(1+u)^\gam-1}=\frac{S(u)}{1+S(u)}=\sum_{p\ge 1}(-1)^{p-1} S(u)^p= \sum_{n\geq 1} (-1)^{n-1}nB_n u^n$$
is analytic for $|u|\! \leq \! 1$, where $nB_n \geq 0$ and can be derived explicitly from the $a_n$.  
Note that $\sum_{n\geq 1} nB_n \! = \! T(1)/(1-T(1))\! =\! \gamma - 1\! \le \! 1$. Therefore, for all $y \ino [-1, 1]$,  $G(y) \! =\! \sum_{n\geq 1} (-1)^{n-1} B_n y^n$,  
which is absolutely convergent; moreover $|G(y)|\! \leq \! -G(-1) \! < \! \sum_{n\geq 1} nB_n \! =\! \gamma \! -\! 1 \! < \! 1$.  Thus,  
$$ \forall y \! \in \! [-1, 1], \quad \exp(G(y))\! =\! 1+\sum_{n\ge 1}A_ny^n \quad \textrm{where} \quad A_n\! =\! (-1)^n \!\!\!\! \!\! \!\! \!\! \!\! \!\! \!\! \! \!\!  \!\! \sum_{ \qquad \begin{subarray}{c}
 p_1, \ldots , p_n\geq 0\\
 p_1+2p_2+\cdots+np_n=n
\end{subarray}} \!\!\!\! \!\! \!\! \!\! \!\! \!\! \!\! \! \!\!    \frac{(-B_1)^{p_1} \ldots (-B_n)^{p_n} }{p_1! \ldots p_n!} .$$
%
We easily see that $1+ \sum_{n\geq 1} |A_n| \! \le \! \exp (-G(-1))\! < \! \exp (\gamma \! -\! 1)$. 
Observe that $\phi(y_0)\! = \! 1$. Then (\ref{cphi1}) follows from (\ref{logcL}) for all $y\! \in \! [y_0, \infty)$. 

We next set $H(y)\! := \! \exp (C_0+ G(y))$. By (\ref{zahcgdfd}), $H$ has a power expansion whose radius of convergence is larger that $1$. 
By Lagrange inversion (recalled in Proposition \ref{Lagraninv}, in Appendix) 
there exists $x_0 \! \in \! (0, \infty)$ such that for all $x\ino [-x_0, x_0]$, the equation $z\! = \!x H(z)$ has a unique unique solution $z\! =:\! f(x)$ 
in $[-1/2, 1/2]$; moreover, for all $ x \ino [-x_0, x_0]$
\begin{equation}
\label{dycgbdgx}
 f(x):= \sum_{n\geq 1} \beta_n x^n \quad \textrm{where} \quad \forall n\! \ge \! 1\, , \quad  \beta_n :=\frac{1}{n!} \, \frac{d^{n-1}}{dy^{n-1}} \big( H^n \big)_{\big| y= 0} \quad \textrm{and} \quad 
\sum_{n\geq 1} |\beta_n| x_0^n\! < \! \infty \; .
\end{equation}
Next observe that (\ref{cphi1}) implies that 
$\phi (y)\! = \! e^{-\gamma y} H(\phi(y)$, for all $y\ino [y_0, \infty)$. Since $\lim_{y\rightarrow \infty} \phi(y) \! = \! 0$, there is 
$y_1\ino [y_0, \infty)$ such that $\phi(y)\ino [0, 1/2]$ for all $y\ino [y_1, \infty)$ and we clearly get $\phi (y)\! = \! f(e^{-\gam y})$, which proves (\ref{echtexphi}). 
An easy computation entails $\beta_1= e^{ C_0}$ and $\beta_2= \frac{\gamma -1}{2} e^{2C_0}$. \cqfd 

\bigskip

We next derive from the previous lemma a similar asymptotic expansion for the function ${\rm L}_1(y, 0)$ that is connected to the diameter of 
$\gamma$-stable normalized trees.  
\begin{lem}
\label{asfundia} Let $\gamma \! \in \! (1, 2]$. Let $\Psi (\lambda)\! = \! \lambda^\gamma$, $\lambda \! \in \! \bbR_+$. 
Recall from (\ref{htdianor}) the definition of ${\rm L}_1(y, 0)$ and recall from 
(\ref{eq: Cgam}) the definition of $C_0$. Then, there exist $y_2 \! \in \! (0, \infty)$, and two real valued sequences 
$(\gamma_n)_{n\geq 2}$, $(\delta_n)_{n\geq 2}$ such that 
\begin{equation}
\label{intvalgd}
\gamma_2\! = \! \frac{_1}{^2} \gamma (\gamma \! -\! 1) e^{2 C_0}, \quad  \delta_2\! =\! - \frac{_1}{^2} (\gamma +1) e^{2C_0} \quad \textrm{and} \quad  \sum_{n\geq 2} \big(n |\gamma_n|+|\delta_n| \big) e^{-\gamma ny_2}  \! < \! \infty 
\end{equation}
and 
\begin{equation}
\label{echtexpLdi}
\forall y\in [y_2, \infty) , \quad {\rm L}_1(y, 0)= \sum_{n\geq 2} (n\gamma_n y + \delta_n) e^{-\gamma n y} \; .
\end{equation}
\end{lem}
\noi
\textbf{Proof.} Recall that $\phi(y)\! = \! w(y)\! -\! 1$. Then (\ref{htdianor}) and an elementary computation entails 
\begin{eqnarray}
\label{formumu}
\mathrm{L}_1(y,0) &= & \phi(y) - \frac{_1}{^\gamma} \big[ (1\! +\! \phi(y))^\gamma \! -\! 1\big] (1\! +\! \phi(y)) + 
\frac{_{\gamma -1}}{^\gamma} \, y \,  \big[ (1\! +\! \phi(y))^\gamma \! -\! 1\big]^2  \nonumber \\
& =&\gamma (\gamma \! -\! 1)y \phi(y)^2 K(\phi(y)) - \frac{_1}{^2} (\gamma +1) \phi(y)^2 M(\phi(y)) \\
\textrm{where for all $u \ino [-1, \infty)$,} & K(u) \!\!\! \!\! & = \! \frac{\big( (u+1)^\gamma \! -\! 1\big)^2}{(\gamma u)^2} \quad \textrm{and} \quad M(u)\! = \! \frac{(u+1)^{\gamma +1} 
\! -\! 1\! -\! (\gamma +1) u}{ \frac{_1}{^2} \gamma(\gamma +1) u^2} . \nonumber
\end{eqnarray}
Recall that $H(y) \! = \! \exp (C_0+G(y))$ and recall from (\ref{cphi1}) that
for all $y \! \in  \! [y_0, \infty)$, $\phi(y)= e^{-\gamma y} H(\phi(y))$. 
This, combined with (\ref{formumu}), entails that 
\begin{equation}
\label{reformu} 
{\rm L}_1(y,0) =  \gamma (\gamma \! -\! 1) e^{-2\gamma y} y H(\phi(y))^2 K(\phi(y)) -\frac{_1}{^2} (\gamma +1) e^{-2\gamma y} H(\phi(y))^2  M(\phi(y)) \; .
\end{equation}
Recall from (\ref{dycgbdgx}), the definition of $f$ and that of $(\beta_n)_{n\geq 1}$. Note that there exists $x_1 \ino (0, x_0)$ such that for all $x\ino [0, x_1]$, 
\begin{align}
\label{expauxii}
  \gamma (\gamma \! -\! 1)& H(f(x))^2 K(f(x))= \sum_{n\geq 0} \gamma^\prime_n x^n \quad \textrm{and} \quad -\frac{_1}{^2} (\gamma +1) H(f(x))^2  M(f(x))=  \sum_{n\geq 0} \delta^\prime_n x^n , \\ \label{intvalgdd}
\textrm{with} & \gamma^\prime_0\! = \! \gamma (\gamma \! -\! 1) e^{2 C_0}, \quad  \delta^\prime_0\! =\! - \frac{_1}{^2} (\gamma +1) e^{2 C_0} \quad \textrm{and} \quad  \sum_{n\geq 0} \big( |\gamma^\prime_n|+|\delta^\prime_n| \big) 
x_2^n   \! < \! \infty , 
\end{align} 
%
%
since $K(0)\! = \! M(0)\! = \! 1$ and since $H(0)^2\! = \! e^{2C_0}$. 
Next by (\ref{echtexphi}) in Lemma \ref{vraiexphi}, we have $\phi(y)\! =\! f(e^{-\gamma y})$, for all $y \! \in \! [y_1, \infty)$. 
Then, we set $ y_2\! : =\!  y_1 \! \vee \! (\! - \tfrac{1}{\gam} \log x_1)$,  and for all $n\! \geq \! 2$, $\gamma_n\! := \! n^{-1}\gamma^\prime_{n-2}$ and $\delta_n\! := \! \delta^\prime_{n-2}$. 
We then see that (\ref{intvalgdd}) implies (\ref{intvalgd}) and that (\ref{expauxii}) and (\ref{reformu}) imply 
(\ref{echtexpLdi}), which completes the proof of the lemma. \cqfd

\subsection{Proof of Theorem \ref{thm: ht_asy}. }
\label{pfasyht}
We first set 
\begin{equation}
\label{fgamdef}
\forall x \in (0,\infty) , \quad f_\Gamma (x) := c_\gamma x^{-1-\frac{1}{\gamma}} \bN_{\! {\rm nr}}
\big( \Gamma \! >\! x^{-\frac{\gamma-1}{\gamma}}\big)\; .
\end{equation}
Then, Proposition \ref{prop: lp}, (\ref{eq: def_bL}), (\ref{scalebL}) and (\ref{htdianor}) imply for all $\lambda \! \in \! (0, \infty)$, 
\begin{equation} 
\label{LapGamnr}
\cL_\lambda (f_\Gamma) \! = \!  
\! \int_0^\infty \!\! \!\! dx \, e^{-\lambda x}\!  f_\Gamma (x)  \! = \! {\rm L}_\lambda (0, 1)\! =\!  \lambda^{\frac{1}{\gamma}} {\rm L}_1 \big(0, \lambda^{\frac{\gamma-1}{\gamma}} \big) 
 \! =\!   \lambda^{\frac{1}{\gamma}}  \big( w\big(\lambda^{\frac{\gamma-1}{\gamma}} \big) \! -\! 1\big)=  
 \lambda^{\frac{1}{\gamma}} \phi \big(\lambda^{\frac{\gamma-1}{\gamma}} \big), 
\end{equation}
where we recall from (\ref{phideff}) that $\phi (y)\! = \! w(y)\! -\! 1$.  We next use Lemma \ref{vraiexphi}: let $\lambda_1$ be such that $\lambda_1^{\frac{\gamma -1}{\gamma}}\! = \! y_1$;  then the sequence $(\beta_n)_{n\geq 1}$ satisfies 
\begin{equation}  
\label{expanfga} 
\forall \lambda \! \in \! [\lambda_1, \infty ), \quad \sum_{n\geq 1} |\beta_n| \lambda^{\frac{1}{\gamma}} e^{-\gamma n\lambda^{\frac{\gamma -1}{\gamma}} } < \infty \quad \textrm{and} \quad  \cL_\lambda (f_\Gamma) = \sum_{n\geq 1} \beta_n \lambda^{\frac{1}{\gamma}} e^{-\gamma n\lambda^{\frac{\gamma -1}{\gamma}} } \; .
\end{equation}
Recall from Lemma \ref{thetpro} the definition of the functions $\theta, h^+$ and $h^-$. Then for all integer $n\! \geq \! 1$, and all $x\! \in \! \bbR_+$, we set 
$$ \theta_n (x) \! = \! n^{-\frac{\gamma+1}{\gamma-1}} \theta \big( n^{-\frac{\gamma}{\gamma-1}} x\big) , \quad 
h^+_n (x) \! =\!  n^{-\frac{\gamma+1}{\gamma-1}} h^+ \big( n^{-\frac{\gamma}{\gamma-1}} x\big) 
\quad \textrm{and} \quad h^-_n (x) \! =\!  n^{-\frac{\gamma+1}{\gamma-1}} h^- \big( n^{-\frac{\gamma}{\gamma-1}} x\big). $$
Lemma \ref{thetpro} implies that $h^+_n$, $h^-_n$ are Lebesgue integrable, nonnegative and continuous. Moreover, $\theta_n\! = \! h^+_n\! -\! h^-_n$. Consequently, $\theta_n$ is also nonnegative continuous and Lebesgue integrable, and 
(\ref{eLapltheta}) entails that $\cL_\lambda (\theta_n)= \lambda^{\frac{1}{\gamma}} e^{-\gamma n\lambda^{\frac{\gamma -1}{\gamma}} } $. Thus, by (\ref{expanfga})
\begin{equation}
\label{identfga}
\forall \lambda \! \in \! [\lambda_1, \infty ), \quad  \cL_\lambda (f_\Gamma) = \sum_{n\geq 1} \beta_n \cL_\lambda (\theta_n) . 
 \end{equation}

We next prove that the assumptions ($a$), ($b$), ($c$) of Lemma \ref{identiLa} hold true with 
$$f\! :=\!  f_\Gamma, \quad f^+_n\! := \! h^+_n, \quad  f^-_n\! := \! h^-_n, \quad \textrm{and} \quad q_n\! :=\! \beta_n\; .$$ 
To that end, we first observe that by an easy change of variable and by 
(\ref{Laplhplms}) in Lemma \ref{thetpro}, we get 
$$\forall \lambda \in (0, \infty), \; \forall n \geq 1, \qquad  \cL_{\lambda} (h^+_n)  \quad  \textrm{and} \quad   \cL_{\lambda} (h^-_n) \; \,  \leq \; (\gamma\! -\! 1)\,  n^{-\frac{1}{\gamma-1}} \!  \! \!  
\int_{n^{\frac{\gamma }{\gamma-1}}\lambda}^\infty \!\!\! \! \!\! d\mu \, e^{-\gamma   \mu^{-\frac{\gamma -1}{\gamma}}} . $$
Thus, by Lemma \ref{quintlaLa}, for all $\lambda \! \in \! (0, \infty)$ and for all sufficiently  large $n$, $\cL_\lambda (h^+_n)$ and $\cL_{\lambda} (h^-_n)$ are bounded by $A\lambda^{\frac{1}{\gamma}} \exp( - \gamma n\lambda^{\frac{\gamma -1}{\gamma}} )$, where $A$ is a positive constant.  Thus, 
  \begin{equation}  
\label{hyposupph} 
\forall \lambda \! \in \! [\lambda_1, \infty ), \quad \sum_{n\geq 1} |\beta_n| \big(  \cL_{\lambda} (h^+_n) + \cL_{\lambda} (h^-_n) \big) \leq \, 2A\!    \sum_{n\geq 1} |\beta_n| 
\lambda^{\frac{1}{\gamma}} e^{-\gamma n\lambda^{\frac{\gamma -1}{\gamma}} } < \infty, 
\end{equation}
the last inequality being a consequence of (\ref{expanfga}).

  Next, deduce from (\ref{domlochpm}) in Lemma \ref{thetpro}  
that for all fixed $x \! \in \! (0, \infty)$ and for all sufficiently large $n$, 
$$  \sup_{y\in [0, x]} h^{\! +}_n \qquad \textrm{and} \quad    \sup_{y\in [0, x]} h^{\! -}_n \quad 
   \leq \quad   B n^{q} x^{-\frac{\gamma +3}{2}} \exp \big( \!\! -\! (\gamma \! -\! 1)^{\gamma -1}n^{\gamma} x^{-(\gamma -1)}  \big) ,  $$
where $q= \frac{\gamma (\gamma+3)}{2(\gamma-1)} -\frac{\gamma+1}{\gamma-1}$ and where $B$ is a positive 
constant only depending on $\gamma$. Since $\gamma \! >\! 1$, $n^\gamma \! \geq \! n$; this combined with (\ref{expanfga}) entails that for all $x\! \in \! \bbR_+$, 
\begin{equation}
\label{autrhyp}
 \sum_{n\geq 1} |\beta_n| \big( \! \sup_{\; y\in [0, x]} \!\! h^{\! +}_n+  \sup_{y\in [0, x]} h^{\! -}_n\, \big) < \infty \; .
\end{equation} 
By (\ref{identfga}), (\ref{hyposupph}) and (\ref{autrhyp}), Lemma \ref{identiLa} applies and we get 
$$\forall x\in \bbR_+, \quad f_\Gamma (x)=c_\gamma x^{-1-\frac{1}{\gamma}} \bN_{\! {\rm nr}}
\big( \Gamma \! >\! x^{-\frac{\gamma-1}{\gamma}}\big)= \sum_{n\geq 1} \beta_n \theta_n (x) \; . $$
This proves  
\begin{equation}
\label{expconcl}
\forall r \in (0, \infty) , \quad  c_\gamma \bN_{\! {\rm nr}}
\big( \Gamma \! >\! r)= \sum_{n\geq 1}  \beta_n \, (nr)^{-\frac{\gamma+1}{\gamma-1}} \theta \big( (nr)^{-\frac{\gamma}{\gamma-1}} \big) , 
\end{equation}  
which implies (\ref{exactoht}). Note that (\ref{autrhyp}) and (\ref{echtexphi}) with $x_1\! =\!  e^{-\gamma y_1}$ in Lemma \ref{vraiexphi} imply (\ref{propbetan}) in Theorem \ref{thm: ht_asy}. 

It remains to prove the asymptotic expansion (\ref{eq: ht_asy}). 
To that end, recall that $ \xi(r)\! =\!  r^{-\frac{\gamma+1}{\gamma-1}} \theta \big( r^{-\frac{\gamma }{\gamma-1}} \big)$, for all $r \! \in \! \bbR_+$. Then (\ref{eq: esth}) in Proposition \ref{thetproo} easily entails that 
for any integer $N \! \geq \! 1$, as $r\to \infty$, 
\begin{equation}\label{xi_asy}
\frac{1}{C^*_1}\, r^{-1-\frac{\gamma}{2}} \, e^{r^\gamma} \, \xi\big( r (\gamma \! -\! 1)^{-\frac{\gamma -1}{\gamma}} \big)\,  = \, 1+ \!\! \sum_{1\leq n<N} \!\!\! V_n \, r^{-n \gamma} \; + \cO_{\! N, \gamma } \big( r^{-N\gamma }\big)\; , 
\end{equation}
where $C^*_1:= (2\pi)^{-\frac{1}{2}} (\gamma \! -\! 1)^{\frac{1}{2}+ \frac{1}{\gamma}} 
\gamma^{\frac{1}{2}}$ and where the sequence $(V_n)_{n\geq 1}$ is recursively defined by (\ref{erecVp}) in Proposition \ref{thetproo}.  This first implies that there exist $A, r_1 \in (0, \infty)$ that only depend on $\gamma$ such that 
\begin{equation}
\label{sblnrxi}
 \forall r \in (r_1, \infty ), \; \forall n\geq 2, \quad  \big| \xi \big( nr (\gamma \! -\! 1)^{-\frac{\gamma -1}{\gamma}} \big)  \big|  \leq A r^{1+ \frac{\gamma}{2}} e^{-n2^{\gamma-1}r^\gamma} \; .
\end{equation} 
Recall from Proposition \ref{thetproo} that there exists $x_1 \! \in \! (0, \infty)$ such that $\sum_{n\geq 1}|\beta_n| x_1^n\! < \! \infty$. Without loss of generality, we can choose $r_1$ such that $\exp (-2^{\gamma-1}r_1^\gamma)\! \leq \! x_1$. Then (\ref{expconcl})    
and (\ref{sblnrxi}) imply that 
$$  \bN_{\! {\rm nr}} \Big( \Gamma \! >\! r (\gamma \! -\! 1)^{-\frac{\gamma -1}{\gamma}} \Big)= c_\gamma^{-1}\!  \beta_1 \, \xi \big(  r (\gamma \! -\! 1)^{-\frac{\gamma -1}{\gamma}} \big)+ \cO_{\! \gamma} \big( r^{1+\frac{\gamma}{2}} e^{-2^{\gamma}r^\gamma}\big), \quad \textrm{as $r\rightarrow \infty$,}$$
and (\ref{xi_asy}) implies  (\ref{eq: ht_asy}) since $C_1\! = \! c_\gamma^{-1} \beta_1 C^*_1$, where we recall from (\ref{eq: defp}) that $c^{-1}_\gamma\! = \! \gamma \Gamma_{\! e} \big( \frac{\gamma -1}{\gamma}\big)$ and where we recall from Lemma \ref{vraiexphi} that $\beta_1\! =\! \exp (C_0)$. This completes the proof of Theorem \ref{thm: ht_asy}.

\subsection{Proof of Theorem \ref{thm: dm_asy}. }
\label{pfasydm}
We first set 
\begin{equation}
\label{fdiadef}
\forall x \in (0,\infty) , \quad f_D (x) := c_\gamma x^{-1-\frac{1}{\gamma}} \bN_{\! {\rm nr}}
\big( D \! >\! 2x^{-\frac{\gamma-1}{\gamma}}\big)\; .
\end{equation}
Then, Proposition \ref{prop: lp}, (\ref{eq: def_bL}) and (\ref{scalebL}) imply for all $\lambda \! \in \! (0, \infty)$, 
\begin{equation} 
\label{LapGamnr}
\cL_\lambda (f_D) =   
\! \int_0^\infty \!\! \!\! dx \, e^{-\lambda x}\!  f_D (x)   =  {\rm L}_\lambda (1, 0) =  \lambda^{\frac{1}{\gamma}} {\rm L}_1 \big(\lambda^{\frac{\gamma-1}{\gamma}} , 0 \big) \; .
\end{equation}
We next use Lemma \ref{asfundia}: let $\lambda_2$ be such that $\lambda_2^{\frac{\gamma -1}{\gamma}}\! \!\! = \! y_2$; then the sequences  $(\gamma_n)_{n\geq 2}$ and $(\delta_n)_{n\geq 2}$ satisfy 
\begin{align}
\label{expanfdia} 
\forall \lambda \! \in  & [\lambda_2, \infty ), \quad \sum_{n\geq 2} \big( n |\gamma_n|  \lambda^{\frac{\gamma -1}{\gamma}} + 
|\delta_n| \big) \lambda^{\frac{1}{\gamma}} e^{-\gamma n\lambda^{\frac{\gamma -1}{\gamma}} } < \infty \nonumber \\
\textrm{and} & \quad  
\cL_\lambda (f_D) = \sum_{n\geq 2} n\gamma_n  \lambda e^{-\gamma n\lambda^{\frac{\gamma -1}{\gamma}} }  
+ \sum_{n\geq 2}\delta_n \lambda^{\frac{1}{\gamma}} e^{-\gamma n\lambda^{\frac{\gamma -1}{\gamma}} } \; .
\end{align}
Recall from (\ref{ethetadef}) in Proposition \ref{thetproo} the definition of $\theta$ and recall Proposition \ref{derivss} that provides properties of the derivative $s^\prime_\gamma $ of the density $s_\gamma$ given by (\ref{denstable}). For all $n \! \geq \! 2$, and all $x\! \in \! (0, \infty)$, we set 
$$ \overline{\theta}_n (x)= n^{-\frac{2\gamma}{\gamma-1}} s_\gamma^\prime \big(  n^{-\frac{\gamma}{\gamma-1}} x\big)  \quad \textrm{and} \quad \theta_n (x)  =  n^{-\frac{\gamma+1}{\gamma-1}} \theta \big( n^{-\frac{\gamma}{\gamma-1}} x\big)  . $$
Then, Proposition \ref{derivss} and Proposition \ref{thetproo} imply that $\overline{\theta}_n$ and $\theta_n$ are continuous and Lebesgue integrable, and that 
$$ \forall \lambda \in \bbR_+, \quad \cL_\lambda(\overline{\theta}_n) = \lambda 
e^{-\gamma n \lambda^{\frac{\gamma-1}{\gamma}}} \quad \textrm{and}  \quad  
\cL_\lambda(\theta_n) = \lambda^{\frac{1}{\gamma}}  
e^{-\gamma n \lambda^{\frac{\gamma-1}{\gamma}}} \; .$$
Thus, 
$$  \forall \lambda \in \bbR_+, \quad \cL_\lambda(f_D)= \sum_{n\geq 2} 
\cL_\lambda \big( n\gamma_n \overline{\theta}_n  + \delta_n \theta_n \big) . $$
We argue as in the proof of Theorem \ref{thm: ht_asy} using Lemma \ref{identiLa} to deduce that 
$$\forall x\in \bbR_+, \quad f_D (x)=c_\gamma x^{-1-\frac{1}{\gamma}} \bN_{\! {\rm nr}}
\big( D \! >\! 2x^{-\frac{\gamma-1}{\gamma}}\big)= \sum_{n\geq 2} n \gamma_n \overline{\theta}_n (x) + \delta_n \theta_n (x) \; , $$
the sum of functions being normally convergent on every compact subset of $\bbR_+$. This easily entails that 
\begin{equation}
\label{expconclD}
\forall r \in (0, \infty) , \quad  c_\gamma \bN_{\! {\rm nr}}
\big(D \! >\! 2r)=
 \sum_{n\geq 2}  \gamma_n \, (nr)^{\! -\frac{\gamma+1}{\gamma-1}} s^\prime_\gamma 
 \big( (nr)^{\! -\frac{\gamma}{\gamma-1}} \big) + \delta_n \, (nr)^{\! -\frac{\gamma+1}{\gamma-1}} \theta \big( (nr)^{\! 
 -\frac{\gamma}{\gamma-1}} \big) , 
\end{equation}  
which is (\ref{exactodiam}). Note that (\ref{propgadeln}) is an easy consequence of the estimate (\ref{asprim}) in Proposition \ref{derivss}, of (\ref{eq: esth}) in Proposition  \ref{thetproo}  and of Lemma \ref{asfundia} with $x_2\! = \! e^{-\gamma y_2}$. 
  Recall from (\ref{edefxibar}) and (\ref{edefxi}) the following notation, 
$$ \forall r \in \bbR_+, \qquad  \overline{\xi} (r)=r^{-\frac{\gamma+1}{\gamma-1}} s^\prime_\gamma  \big( r^{-\frac{\gamma }{\gamma-1}} \big) \quad \textrm{and} \quad \xi (r)=r^{-\frac{\gamma+1}{\gamma-1}} \theta \big( r^{-\frac{\gamma }{\gamma-1}} \big) \; .$$
Note that (\ref{exactodiam}) implies 
\begin{equation}
\label{deucanc}
 c_\gamma \bN_{\! {\rm nr}}
\big(D \! >\! r)= \gamma_2 \overline{\xi} (r) + \delta_2 \xi(r) + \sum_{n\geq 3} \gamma_n 
 \overline{\xi} (nr/2) + \delta_n \xi(nr/2) \; .
\end{equation}
Then, recall from (\ref{xi_asy}) the asymptotic expansion of $\xi$ and deduce from (\ref{asprim}) in Proposition \ref{derivss} that 
\begin{equation}\label{xibar_asy}
\frac{1}{C^*_1}\, r^{-1-\frac{3\gamma}{2}} \, e^{r^\gamma} \, \overline{\xi}\big( r (\gamma \! -\! 1)^{-\frac{\gamma -1}{\gamma}} \big)\,  = \, 1+ \!\! \sum_{1\leq n<N} \!\!\! T_n \, r^{-n \gamma} \; + \cO_{\! N, \gamma } \big( r^{-N\gamma }\big)\; , 
\end{equation}
where $C^*_1:= (2\pi)^{-\frac{1}{2}} (\gamma \! -\! 1)^{\frac{1}{2}+ \frac{1}{\gamma}} 
\gamma^{\frac{1}{2}}$ and where the sequence $(T_n)_{n\geq 1}$ is recursively defined by (\ref{eTdedS}) in Proposition \ref{derivss}. We easily deduce from the asymptotic expansions 
(\ref{xi_asy}) and (\ref{xibar_asy}) that there exists $B, r_2 \! \in \! (0, \infty)$ such that for all $r\! \in \! (r_2, \infty)$ and for all $n\! \geq \!  3$, 
\begin{equation}
\label{sblnrxibar}
 \big| \overline{\xi} \, \big(\frac{_{_1}}{^{^2}}nr   (\gamma \! -\! 1)^{-\frac{\gamma -1}{\gamma}} \big) \big|  \; \textrm{and} \; 
 \big| \xi \, \big(\frac{_{_1}}{^{^2}}nr   (\gamma \! -\! 1)^{-\frac{\gamma -1}{\gamma}} \big) \big|  \;   \leq B r^{1+ \frac{3\gamma}{2}} e^{-n3^{\gamma-1}2^{-\gamma} r^\gamma} \; .
\end{equation} 
This combined with (\ref{deucanc}) implies that 
$$   \bN_{\! {\rm nr}}
\big(D \! >\! r(\gamma \! -\! 1)^{-\frac{\gamma -1}{\gamma}}  \big)=  c_\gamma^{-1} \gamma_2 \overline{\xi} \big( 
r (\gamma \! -\! 1)^{-\frac{\gamma -1}{\gamma}} \big) + c_\gamma^{-1}  \delta_2 \xi \big( r (\gamma \! -\! 1)^{-\frac{\gamma -1}{\gamma}} \big)+ \cO_{\! \gamma} \big( r^{1+ \frac{3\gamma}{2}} e^{-n (3/2)^{\gamma} r^\gamma} \big),  $$
as $r\rightarrow \infty$. Then (\ref{xi_asy}) and (\ref{xibar_asy}) imply 
\begin{align}
\label{beurk}
  \bN_{\! {\rm nr}}
\big(D \! >\! r(\gamma \! -\! 1)^{-\frac{\gamma -1}{\gamma}}  \big)= & c_\gamma^{-1} \gamma_2    C^*_1
 r^{1+ \frac{3\gamma}{2}} e^{-r^\gamma} \nonumber  \\
 + & \sum_{1\leq n< N} c_\gamma^{-1} C^*_1 (\gamma_2 T_n+ \delta_2 V_{n-1} ) 
 r^{-n\gamma +1+\frac{3\gamma}{2}}  e^{-r^\gamma}  + \cO_{\! N, \gamma} \big(r^{-N\gamma +1+\frac{3\gamma}{2}}  e^{-r^\gamma}  \big) 
 \end{align}
Recall from (\ref{intvalgd}) in Lemma \ref{asfundia} that 
$\gamma_2\! = \! \frac{_1}{^2} \gamma (\gamma \! -\! 1) e^{2 C_0}$ and $ \delta_2\! =\! - \frac{_1}{^2} (\gamma +1) e^{2C_0}$. 
This implies (\ref{eq: dm_asy}) with 
$$ C_2\! =\!  c_\gamma^{-1} C_1^*  \gamma_2 \quad \textrm{and} \quad \forall n\geq 1, \quad U_n= T_n + \frac{_{\delta_2}}{^{\gamma_2}} V_{n-1}= T_n-\frac{_{\gamma+1}}{^{\gamma (\gamma-1)}} V_{n-1} \; .$$ 
This completes the proof of Theorem \ref{thm: dm_asy}.

\section{Proof of Theorem \ref{htdm0}.}
\label{Pfhtdm0}
In this section, we fix $\gamma \ino (1, 2)$. Recall that $1/c_\gamma\! =\! \gamma\Gamma_{\! {\rm e}}(1\! -\! \frac{1}{\gamma})$. We set 
\begin{equation}
\label{def: g_ht}
\forall r\in (0, \infty), \quad g_\Gamma(r)\! :=\! c_\gamma r^{-\frac{1}{\gamma}}\nr \big(\Gamma \! \le \! r^{-\frac{\gam-1}{\gam}} \big) \quad \textrm{and} \quad \forall \lambda \ino \bbR, \quad p(\lambda) \! := \! \int_0^\infty \!\!\! \! e^{-\lambda r} \! g_\Gamma (r) \, dr \; . 
\end{equation}
Note that the Laplace transform $p$ is decreasing and that $p (\lambda) \! < \! \infty$ for all $\lambda \ino (0, \infty)$. We next set:  
\begin{equation}
\label{lamsta2}
\llcr\! := \! \sup \Big\{ \lambda \ino \bbR : p( - \lambda) \! = \! \int_0^\infty \!\!\!  \! e^{\lambda r} \! g_\Gamma (r) \, dr \! < \! \infty  \Big\} \quad \textrm{and} \quad \bH\! := \! \big\{ z\ino \bbC :{\rm Re} (z) \! >\! -\! \llcr\big\} \; .
\end{equation}
Clearly $\llcr \! \geq \! 0$. We shall actually prove that $\llcr\ino (0, \infty)$ and that $\int_0^\infty e^{(\llcr-\lam)r}r^2 g_\Gamma (r) \,dr \sim A\lam^{\gam -2}$, for a certain $A\ino (0, \infty)$, as $\lam \! \to\! 0 $. 
However, Karamata's theorem seems to be ineffective to derive asymptotics on $e^{\llcr r}r^2 g_\Gamma (r)$ because this function has no clear monotony properties. 
Thus, we proceed more carefully and we shall use a variant of Ikehara-Ingham Tauberian Theorem to prove Theorem \ref{htdm0}. This requires analytic continuation of $p$. More precisely,  
 standard results on Laplace transform (see for instance Widder \cite{Widder}, Chapter 1) 
imply that $p$ can be analytically extended to $\bH$ by $p(z)\! = \! \int_0^\infty e^{-z r} g_\Gamma (r)\, dr$, for all $z\ino \bH$. We first prove the following lemma. 
\begin{lem}
\label{szcxuca} There exists a real number $\epp_0 \ino (0, \infty)$ and a decreasing analytic function $q\! :(-\epp_0, \infty) \! \rightarrow \! (0, \infty)$ such that 
\begin{equation}
\label{fnksvcq}
\forall \lambda \ino (-\epp_0, \infty) \, , \quad \int_{q(\lambda)}^\infty \! \frac{du}{u^\gamma \! -\! \lam} \! = \! 1  \quad \textrm{and} \quad q^\prime (\lambda) \! = \!  p(\lambda) \! = \! \int_0^\infty \!\!\!  \! e^{-\lambda r} \! g_\Gamma (r) \, dr \; , 
\end{equation}
which implies that $\llcr \! \geq \! \epp_0 $. 
\end{lem}
\noi
\textbf{Proof.} Recall from \eqref{wlam} the definition of $w_\lam(y)$. For all $\lam \ino [0, \infty)$, we set $q(\lam)\! := \! w_\lam(1)\! =\! \bN [1\! -\! e^{-\lam \zeta}\indi_{\{\Gamma\le 1\}}]$. Then, $q$ is clearly decreasing and $C^1$ on $[0, \infty)$. By \eqref{inteqw}, $q$ satisfies 
\begin{equation}
\label{inteqphi}
\forall \lam \ino [0, \infty) \, ,\quad \int_{q(\lam)}^\infty\!  \frac{du}{u^\gam-\lam}\! =\! 1\; . 
\end{equation}
Recall that $\bN(\zeta\ino dr)\! =\! c_\gamma r^{-1-\frac{1}{\gamma}}\, dr$. Thus, by \eqref{zetades}, we get 
$$
\int_0^\infty \!\!\!\! e^{-\lam r}\! g_\Gamma(r)\, dr\! =\! \int_0^\infty \bN(\zeta \ino dr) \, r e^{-\lam r} \nr \big( r^{\frac{\gam-1}{\gam}} \Gam\le 1\big)
\! = \! \bN\big[ \zeta e^{-\lam \zeta} \indi_{\{\Gam\le 1\}} ] \! = \! q^\prime (\lambda) \, .
$$
By (\ref{inteqphi}) we get (\ref{fnksvcq}) for all $\lambda \ino [0, \infty)$ and it is also easy to see that $q(0)\! = \! \bN (\Gamma \! >\! 1)\! = \! (\gam \! -\! 1)^{-\frac{1}{\gam-1}}$.

 Next observe that $q(\lambda) \! >\! \bN [1\! -\! e^{-\lam \zeta}]\! = \! \lambda^{\frac{1}{\gam}}$, which implies $\lambda q(\lam)^{-\gam} \! < \! 1 $, for all $\lambda \ino [0, \infty)$. The change of variable $v\! :=\! u^{-\gam}$ in (\ref{inteqphi}) and the expansion $(1\! -\! \lambda v)^{-1}\! = \! \sum_{n\geq 0} (\lam v)^n$ imply the following. 
 \begin{equation}
 \label{zgfyrbx}
 1\! = \! \frac{1}{\gam} \! \int_0^{q(\lambda)^{-\gam}} \! \!\frac{v^{-\frac{1}{\gam}}  dv}{1\! -\! \lam v}\! = \! \! \frac{1}{\gam} \! \sum_{n\in \bbN} \lambda^n \! \int_0^{q(\lambda)^{-\gam}} \! \!\!\!\!\! v^{n-\frac{1}{\gam}}  dv \! = \! \frac{1}{\gam} q(\lam )^{-(\gam -1)} \sum_{n\in \bbN} \frac{(\lambda q(\lambda)^{-\gam})^n }{n+1\! -\! \tfrac{1}{\gam}} .
\end{equation}  
 This easily implies that for all $\lambda \ino [0, \infty)$, 
 \begin{equation}
 \label{sztxuvc}
\lam q(\lam)^{-\gam} \! = \! \lambda \,  (\gam \! -\! 1)^{\frac{\gam}{\gam -1}}  \Big( 1+ \sum_{n\geq1} \frac{\gam \! -\! 1}{\gam( n+ 1)\! -\! 1} (\lam q(\lam)^{-\gam} )^n \Big)^{-\frac{\gam}{\gam -1}} \; . 
\end{equation}
First note that there is $\epp_1\ino (0, \infty)$ such that the function $H(x) \! := \!   (\gam \! -\! 1)^{\frac{\gam}{\gam -1}} 
( 1+ \sum_{n\geq1} \frac{\gam  - 1}{ \gam (n+ 1) - 1} x^n )^{-\frac{\gam}{\gam -1}}$ has an absolutely convergent power expansion for all $x\ino (-\epp_1, \epp_1)$  
and next observe that (\ref{sztxuvc}) implies 
that $\lambda q(\lam)^{-\gam} \! = \! \lambda H(\lambda q(\lam)^{-\gam} )$ in a right neighbourhood of $0$. Lagrange inversion (as recalled in Theorem \ref{Lagraninv}) 
implies that there is $\epp_2 \ino (0, \infty)$ such that $\lambda \!  \mapsto\!  \lambda q (\lambda)^{-\gam}$ extends  analytically on $(-\epp_2, \epp_2)$. This implies that there exists a sequence of real numbers $(a_n)_{n\in \bbN}$ and a real number $\epp_0 \ino (0, \infty)$ such that the following power expansion 
\begin{equation}
\label{dhchcdg}
q(\lambda) = \sum_{n\in \bbN} a_{n} \lambda^n  \; , \quad \lambda \ino (-\epp_0, \epp_0) \; , 
\end{equation} 
is absolutely convergent and such that $-\lambda q(-\lam)^{-\gam} \! = \! -\lambda H(-\lambda q(-\lam)^{-\gam} )$, for all $\lambda \ino [0, \epp_0)$. This equality easily implies that (\ref{zgfyrbx}) holds true with 
$-\lambda$ instead of $\lam$, namely: 
\begin{equation}
\label{ztxdov}
\forall \lambda \ino [0, \epp_0)\, , \quad \int_{q(-\lambda)}^\infty \! \frac{du}{u^\gamma +\lam} \! = \! 1 \; .
\end{equation}
Since $q^\prime(\lambda)\! = \! \int_{0}^\infty e^{-\lambda r} \! g_\Gamma (r) \, dr $, for all $\lam \ino [0, \infty)$, (\ref{dhchcdg}) and 
standard results on the Laplace transform (see for instance Widder \cite{Widder}, Chapter 1) imply that 
$$ \forall n\ino \bbN, \quad \frac{1}{n!} \! \int_0^\infty \!\!\!\!  r^n g_\Gamma (r) \, dr \! =\!  (-1)^n (n+1)a_{n+1} \; .$$
Since $\sum_{n \in \bbN} |(n+1)a_{n+1} \lam^n| \! < \! \infty $, for all $\lambda \ino (-\epp_0, \epp_0)$, this implies that 
$$ \forall \lam \ino [0, \epp_0)\, , \quad q^\prime (-\lam)\! = \! \sum_{n\in \bbN} (-1)^n (n+1)a_{n+1} \lambda^n \! =\! \sum_{n\in \bbN} \frac{\lambda^n}{n!} \! \int_0^\infty \!\!\!\!  r^n g_\Gamma (r) \, dr \! = \!  \int_0^\infty \!\!\!  \! e^{\lambda r} \! g_\Gamma (r) \, dr \; .$$
This, combined with (\ref{ztxdov}), completes the proof of (\ref{fnksvcq}).  \cqfd 

\bigskip

We next set $D_-\! := \! \{ z\ino \bbC: \textrm{${\rm Re} (z) \! \le \! 0$ and ${\rm Im} (z)\! = \! 0$}\}$, the negative axis of the complex plane. For any $b\ino \bbC$, we use the following notation 
\begin{equation}
\label{fbvuevcx}
 \forall z\ino \bbC \backslash D_-\, , \quad z^b\! := \! \exp (b \log z ) \; , 
\end{equation}  
 where $\log$ is the usual determination of the logarithm in $\bbC \backslash D_-$. Standard results in complex analysis assert that $z\! \mapsto \! z^b$ is analytic in the domain $\bbC\backslash D_-$. 
The following lemma concerns the analytic continuation of $q$ introduced in Lemma \ref{szcxuca}. Recall from (\ref{lamsta2}) the definition of $\llcr$ and that of 
the right half-plane $\bH$. 
\begin{lem}
\label{zeuccp} There exists a connected open subset $U$ containing $\overline{\bH}\backslash \{ \! -\! \llcr\}$ such that the function $q$ (introduced in Lemma \ref{szcxuca}) has an analytic continuation to 
$U$ that is $C^2$ on $\overline{\bH}$ and such that 
$q^\prime (z)\! = \! \int_0^\infty e^{-zr} g_\Gamma (r) \, dr$, for all $z\ino \overline{\bH}$. Moreover, $q$ satisfies the following properties. 
\begin{itemize}
\item[(i)] Let $U_0$ denote the open strip $\{ -\llcr \! < \! {\rm Re} (z) \! < \! 0 \}$. Then $q$ satisfies 
\begin{equation}
\label{diff_eq}
\forall z\ino U_0\, , \quad q(z) \ino \bbC\backslash D_-\quad \textrm{and} 
\quad zq'(z)=-\frac{_{\gam-1}}{^\gam}q(z)^\gam+\frac{_1}{^\gam}q(z)+\frac{_{\gam-1}}{^\gam}z \; .
\end{equation} 
\item[(ii)] $q(-\llcr)\!= \! 0$ and as $z\! \to \! 0$ with ${\rm Re} (z) \! >\! 0$, 
\begin{align}
\label{phi_ddd}
&q^{(3)}(-\llcr+z)=\frac{{(\gam\! -\! 1)^{\gam+2}}}{{\gam^\gam \llcr}}\, z^{\gam-2}-\frac{{(2\gam\! -\! 1)(\gam\! -\! 1)}}{{\gam^3\llcr^2}}+ {\rm o}(1), \\
\label{phi_d4}
&q^{(4)}(-\llcr+z)=\frac{{(\gam \! -\! 1)^{\gam+2}(\gam-2)}}{{\gam^\gam \llcr}}\, z^{\gam-3}+\frac{{(\gam\! -\! 1)^{\gam+3}(\gam+2)}}{{2\gam^{\gam+1}\llcr^2}}\, z^{\gam-2}+ {\rm o}(z^{\gam-2}).
\end{align}
\item[(iii)] $-\llcr$ is the only singular point of $q$ in $U$ and $\llcr\! = \! (\frac{\pi/\gam}{\sin (\pi / \gam)})^{\frac{\gamma}{\gam -1}} $. 
\end{itemize}
\end{lem}
\begin{rem}
\label{fufluns}
The statement in Lemma \ref{zeuccp} is not valid for $\gam\! =\! 2$. Indeed, if $\gam\! =\! 2$, for all $\lam\ino (0, \infty)$, 
$q(\lam)\! =\! \sqrt\lam\coth\sqrt\lam$ and $q(-\lam)\! =\! \sqrt\lam\cot\sqrt\lam$. Therefore, $q$ is analytic on $(-\pi^2, \infty)$. But note that $(\frac{\pi/\gam}{\sin (\pi / \gam)})^{\frac{\gamma}{\gam -1}}\! = \! \pi^2/4$ when $\gam \! = \! 2$. The reason for the distinct behaviour of $q$ when $\gam\! = \! 2$ boils down to the elementary fact that $0$ is a singular point for $z\! \mapsto \! z^\gam$ when $\gam\ino (1, 2)$. It is not the case when $\gam \! = \! 2$. \cq 
\end{rem}

\noi
\textbf{Proof of Lemma \ref{zeuccp}.} Let $\lam \ino (0, \infty)$. By the change of variable $v\! := \! \lambda u^{-\gam}$, we get 
$$ \int_0^\infty \!\! \! \frac{du}{u^\gam + \lam}  =  \tfrac{1}{\gam} \lambda^{-\frac{\gam -1}{\gam}}\!\!  \int_0^\infty \!\frac{v^{-\frac{1}{\gam}}dv}{1+v} =  \frac{\pi/ \gam}{\sin (\pi /\gam)} \lambda^{-\frac{\gam -1}{\gam}}\; .$$
Here, we use E.~Schl$\overset{..}{{\rm a}}$fli's identity $\int_0^\infty \! v^{-s}\!  / (1+v) \, dv \! = \! \pi / \sin (\pi s)$, that is valid for all $s\ino \bbC$ such that $0\! < \! {\rm Re} (s) \! < \! 1$ (see for instance I.~Gradshteyn \& I.~Ryzhik \cite{GraRyz07}, Chapter 17, Section 43, p.~1131, Table of Mellin transform, formula 6). We then set 
$$ \lambda_1 \! := \!  \Big(\frac{\pi/\gam}{\sin (\pi / \gam)} \Big)^{\frac{\gamma}{\gam -1}} \quad \textrm{that satisfies} \quad \int_0^\infty \!\! \! \frac{du}{u^\gam + \lam_1} =1 \; .$$
Therefore, there exists a strictly decreasing continuous function $r\! : [0, \lambda_1] \! \mapsto \! [0, q(0)]$ that satisfies 
\begin{equation}
\label{ducvcjd}
\forall \lambda \ino [0, \lambda_1]\, , \quad  \int_{r(\lambda)}^\infty   \!\frac{du}{u^\gam + \lam} =1 \; .
\end{equation}
Note that $r(0)\! = \! q(0)\! = \! (\gam \! -\! 1)^{-\frac{1}{\gam -1}}$, that $r(\lambda_1)\! = \! 0$. By Lemma \ref{szcxuca}, $\epp_0 \! \le \! \lambda_1$ and 
$r(\lambda)\! = \! q(-\lambda)$, for all $\lambda \ino [0, \epp_0)$. An easy linear change of variable in (\ref{ducvcjd}) entails 
$$\lambda^{\frac{\gamma -1}{\gamma}}\! = \! \int^\infty_{\lambda^{-\frac{1}{\gam}} r(\lam)} \! \frac{dv }{v^\gam +1} \quad \textrm{and thus} \quad - r^\prime(\lam)\! = \! \frac{_{\gam -1}}{^\gam} \lambda^{-1} r(\lam)^\gam -\frac{_1}{^\gam} \lambda^{-1} r(\lam) + \frac{_{\gam -1}}{^\gam} \, , \quad \lam \ino (0, \lam_1) \; .$$ 
Thus, we have proved that $q$ can be extended uniquely on $[-\lam_1, \infty)$ in such a way that $\int_{q (\lambda)}^\infty du / (u^\gam \! -\! \lam)\! = \! 1$ for all 
$\lam \ino [-\lam_1, \infty)$ and we have 
$$ \forall \lambda \ino (-\lambda_1, 0)\, , \quad  q^\prime (\lambda)= F\big( \lambda , q(\lambda)\big) \; , $$
where we have set 
\begin{equation}
\label{ucsjdbc}
 \forall (z,v) \ino V \, , \quad F(z, v)\! :=\!  -\tfrac{\gam -1}{\gam} z^{-1}v^\gam +\tfrac{1}{\gam} z^{-1} v  + \tfrac{\gam -1}{\gam} \quad \textrm{where} \quad V\! :=\! (\bbC \backslash \{ 0\}) \! \times \! (\bbC \backslash D_-)
\end{equation}
Note that $V$ is an open subset of $\bbC^2$ and we recall the convention specified by (\ref{fbvuevcx}) for the power of complex numbers. Recall that $D(z_0, r)$ stands for 
the open disk in $\bbC$ with centre $z_0\ino \bbC$ and radius $r\ino (0, \infty)$; to simplify notation we identify $\bbR$ with the set of complex numbers whose imaginary part is null. 
We next use Proposition \ref{ComplODE} (see Appendix, Section \ref{comanarec}). First, we easily check that $F$ is analytic in the two variables $z$ and $v$ on $V$.  
Then, for all $\lambda \ino (-\lambda_1, 0)$, since $(\lambda, q(\lambda)) \ino V$, Proposition \ref{ComplODE} implies that there exists $r_\lam\ino (0, \infty)$ and an analytic function 
$f_\lam: D(\lambda , r_\lam) \rightarrow \bbC\backslash D_-$ such that $f_\lambda$ is the unique solution of 
$$ \forall z\ino D(\lam, r_\lam)\, , \quad f_\lam (z) \ino V\, , \quad  f^\prime_\lam (z) \! =\!  F\big( z, f_\lam (z)\big) \quad \textrm{and} \quad f_\lambda (\lam)\! =\!  q(\lam)\; .$$
The restriction of $f_\lam$ on the real interval $(\lam - r_\lam, \lam + r_\lam)$ clearly satisfies the same (real time parameter) ordinary differential equation as $q$; since this ODE is locally Lipschitz, uniqueness in 
the Picard-Lindel$\overset{..}{\rm o}$f Theorem (also known as Cauchy-Lipschitz Theorem) implies that $f_\lam$ and $q$ coincide on the real interval $(\lam\! -\! r_\lam, \lam + r_\lam)$. Let $\lambda, \lambda^\prime \ino (-\lam_1, 0)$ be such that 
$W\! :=\! D(\lam, r_\lam)\cap D (\lambda^\prime, r_{\lam^\prime})\! \neq \emptyset$; since $W$ is connected and since $f_\lam$ and $f_{\lam^\prime}$ are equal to $q$ on the real interval $W\cap \bbR$, the principle of isolated zeroes for analytic functions implies that $f_\lam$ and $f_{\lam^\prime}$ coincide on $W$. 
This implies that $q$ can be extended uniquely on the open subset $U_1\! := \! \bigcup_{\lambda \in (-\lam_1, 0)} D(\lam, r_\lam)$, that 
$q: U_1\rightarrow \bbC\backslash D_-$ is analytic and that $q$ satisfies the complex differential equation: 
\begin{equation}
\label{iopkpikk}
 \forall z \ino U_1, \quad q(z) \ino V \, , \quad q^\prime (z)= F\big( z, q(z) \big) \; .
\end{equation} 
Since $(-\lam_1, 0) \! \subset \! U_1$, this implies that the restriction of $q$ on $(-\lam_1, \infty)$ is analytic. We next prove that it entails that 
\begin{equation}
\label{decxtqs}
\llcr \! \ge \! \lam_1 \quad \textrm{and} \quad \forall \lam \ino [-\lam_1, \infty) , \quad q^\prime(\lam) \! =\!  \int_0^\infty \!\!\!\!  e^{-\lambda r}\! g_\Gamma (r) \, dr \; .
\end{equation}
\textit{Indeed}, suppose that $\llcr \!  < \! \lam_1$. By standard results on Laplace transform 
$\lambda \! \mapsto \! \int_0^\infty \! e^{-\lam r} g_\Gamma (r) \, dr $ is analytic on $(-\llcr, \infty)$. Lemma \ref{szcxuca} implies that it coincides with $q^\prime$ on 
$(-\epp_0, \infty)$. Since $q^\prime$ is also analytic on the interval $(-\llcr, \infty)$ (supposedly included in $(-\lam_1, \infty)$), 
the principle of isolated zeroes for analytic functions entails that $q^\prime (\lambda)\! = \! \int_0^\infty e^{-\lam r} g_\Gamma (r) \, dr $, for all $\lambda\ino (-\llcr, \infty)$. Standard results on Laplace transform also imply that for all $n\ino \bbN$, $\int_0^\infty e^{-\lam r} r^ng_\Gamma (r) \, dr\! = \! (-1)^nq^{(n+1)} (\lam)$, for all $\lam\ino (-\llcr, 0)$. By continuity of $q^{(n+1)}$ and the monotone convergence theorem, 
we get $\int_0^\infty e^{\llcr r} r^n g_\Gamma (r) \, dr\! = \! (-1)^n q^{(n+1)} (-\llcr)$. Since $\lam_1 \! >\! \llcr$, $q^\prime$ is analytic at $\llcr$ and there exists $\epp \ino (0, \lambda_1\! -\! \llcr)$ such that 
$$ \int_0^\infty \!\!\!\!  e^{(\llcr+ \epp) r}\! g_\Gamma (r) \, dr= \sum_{n\in \bbN } \frac{\epp^n}{n!} \int_0^\infty \!\!\!\! e^{\llcr r}\! r^n g_\Gamma (r) \, dr= \sum_{n\in \bbN} \frac{(-\epp)^n}{n!}q^{(n+1)} (-\llcr)= q^\prime (-\llcr\! -\! \epp) <\infty \; , $$
which contradicts the definition (\ref{lamsta2}) of $\llcr$. Thus $\lam_1 \! \le \! \llcr$ and (\ref{decxtqs}) holds true. 

We set $\bH_1\! := \{ z\ino \bbC: {\rm Re} (z) \! >\! -\! \lam_1 \}$ and we next prove that $q$ can be extended analytically on $\bH_1$, that $q$ is continuous on $\overline{\bH}_1$ and that ${\rm Re} (q(z)) \! >\! 0$, for all $z\ino \overline{\bH}_1\backslash \{ -\lam_1\}$. 
\textit{Indeed}, (\ref{decxtqs}) implies that $q^\prime$ can be extended analytically on $\bH_1$ and that $q^\prime (z)\! = \! \int_0^\infty e^{-zr} g_\Gamma (r) \, dr$, for all $z\ino \bH_1$. Thus, $q$ can be extended analytically on $\bH_1$ and we easily get $q(z)\! = \! q(0)\! -\! \int_{0}^\infty g_\Gamma (r) r^{-1} (e^{-zr} \! -\! 1) \, dr$, for all $z\ino \bH_1$. Since $\lam\! \mapsto \! q(\lam)$ decreases to $q(\! -\! \lam_1)\! = \! 
0$ as $\lam\!  \downarrow \! -\! \lam_1$, monotone convergence theorem implies that $\int_{0}^\infty g_\Gamma (r) r^{-1} (e^{\lam_1r} \! -\! 1) \, dr\! = \! q(0) \! < \! \infty$. It thus implies that 
\begin{equation}
\label{dhcxiahd}
\forall z\ino \overline{\bH}_1\, , \quad q (z) \! = \!\! \int_0^\infty \!\!\! \! dr \, g_\Gamma (r)  r^{-1} \big( e^{\lam_1 r} \! -\! e^{-zr} \big) ,
\end{equation}
and $q$ is continuous on $\overline{\bH}_1$. For all $\lam \ino [-\lam_1, \infty)$ and all $t\ino \bbR$, we also get 
$$ {\rm Re} (q(\lam + it))\! = \! \int_0^\infty \!\!\! \! dr \, g_\Gamma (r)  r^{-1} e^{-\lam r}\big( e^{(\lam_1 +\lam) r} \! -\! \cos (tr) \big) \; .$$
If $t\! \neq \! 0$ or $\lam\! \neq \! -\! \lam_1$, then $r\mapsto g_\Gamma (r)  r^{-1} e^{-\lam r}( e^{(\lam_1 +\lam) r} \! -\! \cos (tr)) $ is nonnegative and strictly positive on a non-empty interval. Thus, 
${\rm Re} (q(z)) \! >\! 0$, for all $z\ino \overline{\bH}_1\backslash \{ -\lam_1\}$.

We denote by $U_2$ denote the open strip $\{ -\lambda_1 \! < \! {\rm Re} (z) \! < \! 0 \}$. We next prove that 
\begin{equation}
\label{dcijfplm}
\forall z\ino U_2\, , \quad q^\prime (z) = F(z, q(z)) \; , 
\end{equation}
where we recall from (\ref{ucsjdbc}) the definition of the open set $V$ and the function 
$F\! : V \! \rightarrow \! \bbC \backslash D_-$. 
We then fix $\lambda \ino (-\lam_1, 0)$ and we consider $y: I\rightarrow \bbC \backslash D_-$, the maximal solution of the (real time parameter) ordinary differential equation 
\begin{equation}
\label{djcvszkf}
\forall t \ino I \, , \quad  y^\prime (t)= iF(\lam + it, y(t)) \quad \textrm{and} \quad y(0)=q(\lam) \; .
\end{equation}
Here, $I$ is the maximal (open) interval of definition for (\ref{djcvszkf}). Existence and uniqueness of such a maximal solution is a consequence of Picard-Lindel$\overset{..}{{\rm o}}$f Theorem. 
Recall (\ref{iopkpikk}) and recall that by definition $(-\lam_1, 0) \! \subset \!  U_1$. Thus, there exists $\epp\! >\! 0$ such that $(-\epp, \epp) \! \subset \! I$ and $y(t)\! = \! q(\lam +it)$, for all $t\ino (-\epp, \epp)$. Next, observe that $(\lam+ is, y(s))\ino V$ for all $s\ino I$; then by Proposition \ref{ComplODE}, there exist $\eta_s\ino (0, \infty)$ and an analytic function 
$h_s\! : D(\lam+is, \eta_s )\!  \rightarrow \!  \bbC \backslash D_-$ such that $h^\prime_s (z)\! =\!  F(z, h_s (z))$, for all $z\ino  D(\lam+is, \eta_s )$ and $h_s (\lam +is)\! = \! y(s)$. Thus $t\ino (s\! -\! \eta_s, s+ \eta_s) 
\! \mapsto \!  h_s (\lam+it)$ satisfies the same (real time parameter) ODE as $y$ and thus $h_s (\lam+it)\! = \! y(t)$, for all $t\ino (s\! -\! \eta_s, s+ \eta_s)$. Let $s, s^\prime \ino I$ be such that 
$W\! :=\! D(\lam+is, \eta_s ) \cap D(\lam+is^\prime, \eta_{s^\prime} )\! \neq \! \emptyset$; since $W$ is connected and since $h_s$ and $h_{s^\prime}$ are equal to $y$ on $W\cap (\lam +i\bbR)$ (with an obvious notation), the principle of isolated zeroes for analytic functions implies that $h_s$ and $h_{s^\prime}$ coincide on $W$. Thus, there is an analytic function $w$ from the open set $O\! := \! \bigcup_{s\in I}  D(\lam+is, \eta_s )$ to $\bbC \backslash D_-$ such that $w^\prime (z)\! = \! F(z, w(z))$ and such that $w(\lam+it)\! = \! y(t)$, for all $t\ino I$. Note that 
$O$ is connected and that $O \! \subset \! \bH_1$; since $w(\lam+it)\! = \! y(t)\! = \! q(\lam + it)$, for all $t\ino (-\epp, \epp)$, the principle of isolated zeroes for analytic functions implies that $q$ and 
$w$ coincide on $O$. This proves that $q^\prime (z)\! = \! F(z, q(z))$ for all $z\ino O$ and that $q(\lam +it)\! = \! y(t)$, for all $t\ino I$.   
If we prove that $I\! = \! \bbR$, then the previous arguments entail $q(\lam +it)\! = \! y(t)$ for $t\ino \bbR$, and $q^\prime (\lam+ it)\! = \! F(\lam+it, q(\lam +it))$, $t\ino \bbR$, which implies (\ref{dcijfplm})
since $\lam$ is arbitrarily chosen in $(-\lam_1, 0)$. 

Let us prove that $I\! = \! \bbR$. We argue by contradiction: assume first that $I$ has a bounded right end denoted by $a$, namely $I\cap [0, \infty)\! = [0, a)$. By continuity of $q$, $\lim_{t\rightarrow a-} y(t)\! = \! q(\lam + ia)$; since ${\rm Re} (q(z)) \! >\! 0$, for all $z\ino \overline{\bH}_1\backslash \{ -\lam_1\}$, we get $(\lam+ia, q(\lam +ia))\ino V$ and by Proposition \ref{ComplODE}, there exist $\eta\ino (0, \infty)$ and an analytic function 
$h\! : D(\lam+ia, \eta )\!  \rightarrow \!  \bbC \backslash D_-$ such that $h^\prime (z)\! =\!  F(z, h (z))$, for all $z\ino  D(\lam+ia, \eta )$ and $h (\lam +ia)\! = \! q(\lam+ ia)\! = \! y(a-)$. Then set $x(t)\!= \! y(t)$, $t\ino I$ and $x(t)\! = \! h(\lam +i t)$ for all $t\ino [a, a+ \eta)$; we observe that $x$ satisfies the same (real time parameter) ODE as $y$ and that it strictly extends $y$,
which contradicts the definition of $I$. Thus $I$ is unbounded from the right. We argue in the same way to prove that $I$ is unbounded from the left, which proves that $I\! = \! \bbR$ and (\ref{dcijfplm}) as already mentioned.  

We thus have proved that $q$ can be extended analytically on $\bH_1$, that $q$ is continuous on $\overline{\bH}_1$ and that $q$ satisfies (\ref{dcijfplm}). Recall that $q(\!- \! \lam_1)\! = \! 0$, which implies by  (\ref{dcijfplm}) that $q^\prime (-\lam_1 +z)$ tends to $F(\! -\! \lam_1, 0)\! = \! (\gam\! -\! 1)/ \gam$ as $z\! \rightarrow \! 0$ with ${\rm Re} (z) \! >\! 0$. We then set $q^\prime (\! -\! \lam_1)\! := \!  (\gam\! -\! 1)/ \gam$; 
(\ref{dhcxiahd}) and monotone convergence entail $\int_0^\infty \! e^{\lam_1 r} g_\Gamma (r) \, dr \! = \! \lim_{\lam \downarrow -\lam_1} q^\prime (\lam)\! =\!  (\gam\! -\! 1)/ \gam $. This also proves that $q^\prime$ is continuous on $\overline{\bH}_1$. Therefore $q$ is $C^1$ on $\overline{\bH}_1$. We also derive from (\ref{dcijfplm}) that 
\begin{equation}
\label{fkvudc}  
\forall z\ino U_2\, , \quad -zq^{\prime \prime} \! (z)\! = \! \big(1\!-\! \tfrac{1}{\gam} + (\gam \! -\! 1) q(z)^{\gam -1} \big) q^\prime (z) +  \tfrac{1}{\gam} \! -\! 1 \; .
\end{equation}
Thus, $q^{\prime \prime} (-\lam_1 +z)$ tends to $q^{\prime \prime}(\! -\! \lam_1) \! := \! -(\gam  - 1)/(\lam_1 \gam^2)$ as $z\! \rightarrow \! 0$ with ${\rm Re} (z) \! >\! 0$ and monotone convergence entails that $\int_0^\infty \! re^{\lam_1 r} g_\Gamma (r) \, dr \! = \! -q^{\prime \prime}(\! -\! \lam_1)$, which implies that $q^\prime$ is $C^1$, and therefore that $q$ is $C^2$ on $\overline{\bH}_1$. 
We next observe that for all $z \ino \bbC $ such that ${\rm Re} (z) \! >\! 0$, we get 
\begin{eqnarray}
\label{djxfbpl}
q(\! -\! \lam_1+z)& = & zq^\prime (\! -\! \lam_1) + \tfrac{1}{2}z^2 q^{\prime \prime} (\! -\! \lam_1) + z^2 \!\! \int_0^1 \!\!\! dt \!  \int_0^t \! \!\! ds \, \big( q^{\prime \prime} (\! -\! \lam_1+ sz)\! -\! q^{\prime \prime} (\! -\! \lam_1) \big) \nonumber \\
& =& \frac{\gam \! -\! 1}{\gam} \,  z- \frac{\gam \! -\! 1}{2\lam_1\gam^2} \, z^2 + \mathrm{o} (z^2) ,  
\end{eqnarray}
as $z\! \to \! 0$. A similar argument entails that 
\begin{equation}
\label{qq_ty}
q'(\! -\! \lam_1+z)=\frac{\gam\! -\! 1}{\gam}\! -\! \frac{\gam\! -\! 1}{\gam^2\lam_1}\, z+\mathrm{o}(z), \quad\text{ and }\quad q''(\! -\! \lam_1+z) \! =\! - \frac{\gam\! -\! 1}{\gam^2\lam_1}+ \mathrm{o}(1), 
\end{equation}
as $z\! \to \! 0$ with ${\rm Re} (z) \! >\! 0$. We next derive from (\ref{fkvudc}) that for all $z\ino U_2$, 
\begin{eqnarray}
\label{diff_ddd}
&  & -zq^{(3)}(z)  =  \big( (\gam\! -\! 1) q(z)^{\gam -1}\! + 2\! -\! \tfrac{1}{\gam}\big)q^{\prime \prime} (z) +  (\gam\! -\! 1)^2 q^\prime (z)^2 q(z)^{\gam -2} \quad \textrm{and}\\
-zq^{(4)}(z) \!\!\!\!\!\!  & = &  \!\!\!\!\!\!   \big( (\gam\! -\! 1) q(z)^{\gam -1}\!\! + 3\! -\! \tfrac{1}{\gam}\big)q^{(3)} (z) +  3(\gam\! -\! 1)^2 q^{\prime \prime} (z)  q^\prime (z) q(z)^{\gam -2} \!\! + (\gam \! -\! 2) (\gam\! -\! 1)^2q^\prime (z)^3 q(z)^{\gam -3}\!\!\! . \nonumber
\end{eqnarray}
This entails that $\lim_{\lam \downarrow -\lam_1} q^{(3)} (\lam)\! = \! \infty$ and thus $-\lam_1$ is a singular point of $q$. Consequently $\lam_1\! = \! \llcr$. Moreover, (\ref{diff_ddd}) combined with (\ref{djxfbpl}) and (\ref{qq_ty}) entails (\ref{phi_ddd}) and (\ref{phi_d4}). 

It remains to prove that $q$ can be extended on an open subset containing $\overline{\bH}\backslash \{ -\llcr \}$. To that end, we recall that for any $t\ino \bbR\backslash \{ 0\}$, ${\rm Re}(q(\! -\! \llcr + it))\! >\! 0$. Thus, 
$(\! -\! \llcr + it, q(\! -\! \llcr + it))\ino V$ and Proposition \ref{ComplODE} implies that there exists $\rho_t\ino (0, \infty)$ and a unique analytic function $k_t\! : D(\! -\! \llcr  +it, \rho_t) \! \to \! \bbC\backslash D_-$ such that $k_t(\! -\! \llcr  +it)\! = \! q(\! -\! \llcr  +it)$ and $k^\prime_t(z)\! = \! F(z, k_t(z))$, for all $z\ino D(\! -\! \llcr  +it, \rho_t)$. Since $q$ satisfies the same differential equation on $\overline{\bH}\cap D(\! -\! \llcr  +it, \rho_t)$, we see that the function $x\ino [\! -\! \llcr, \! -\! \llcr+\rho_t)\! \mapsto \! q(x+it)$ and the function $x\ino [\! -\! \llcr, \! -\! \llcr+\rho_t)\! \mapsto \! k_t(x+it)$ satisfy the same (real time parameter) ODE, with the same initial condition. Since this ODE is locally Lipschitz in space, uniqueness of the solution in Picard-Lindel$\overset{..}{{\rm o}}$f Theorem entails that 
$k_t(x+it)\! = \! q(x+it)$, for all $x\ino [\! -\! \llcr, \! -\! \llcr+\rho_t)$. Since $k_t$ and $q$ are analytic on the connected open set $\bH \! \cap \! D(\! -\! \llcr  +it, \rho_t)$, the principle of isolated zeroes for analytic functions entails that $k_t$ and $q$ coincide on  $\bH \cap D(\! -\! \llcr  +it, \rho_t)$ and thus on $\overline{\bH}\cap D(\! -\! \llcr  +it, \rho_t)$. Let $t, t^\prime\ino \bbR \backslash \{ 0\}$ be such that $W\! :=\! 
 D(\! -\! \llcr  +it, \rho_t)\! \cap \! D(\! -\! \llcr  +it^\prime, \rho_{t^\prime})$ is non-empty. Since $k_t$ and $k_{t^\prime}$ are analytic on the connected open set $W$ and since they coincide with $q$ on the non-empty connected set $W\! \cap \! \bH$, the principle of isolated zeroes for analytic functions entails that $k_t$ and $k_{t^\prime}$ coincide on $W$. We now set $U\! := \! \bH \cup  
 \bigcup_{t\in \bbR \backslash \{ 0\}} D(\! -\! \llcr  +it, \rho_t)$. The previous arguments show that $q$ can be extended analytically on $U$ and obviously $U$ contains $\overline{\bH} \backslash \{ \! -\! \llcr\}$, 
which completes the proof of Lemma \ref{zeuccp}. \cqfd 

\paragraph{Proof of (\ref{dtcxtqcd}) in Theorem \ref{htdm0}.}
Next we want to apply Ikehara-Ingham Theorem that is recalled in Theorem \ref{IKEA} in Appendix. To that end, we next prove the following lemma. 
\begin{lem}
\label{sazerqsd} For all $z \ino \bbC$ such that $0\! < \! {\rm Re} (z) \! < \! \llcr$, we set
\begin{equation}
\label{poiuymlk}
G(z)\! := \!  \frac{q^{(3)} ( - \llcr+z)}{\llcr\! -\! z}  - \frac{(\gam\! -\! 1)^{\gam+2}}{\llcr^2 \gam^\gam}\, z^{\gam -2} \; .
\end{equation}
Then, for all $\theta \ino (0, \infty)$, 
\begin{equation}
\label{fructimus}
 \lam^{1-\gam} \int_{-\theta}^\theta \big|G(2\lam+i t)\! -\! G(\lam+it)\big| \, dt \; \underset{\lam \downarrow 0+}{-\!\!\! -\!\!\! \longrightarrow} \; 0 \; . 
\end{equation} 
\end{lem}
\noi
\textbf{Proof.} We fix $\lam \ino (0, \llcr/2)$, $\theta \ino (0, \infty)$ and $t\ino (-\theta , \theta)$. Observe that  
\begin{align*}
&G(2\lam+it)\! -\! G(\lam+it)=\int_\lam^{2\lam} \!\!\!\! \!\!\! du \, G^\prime (u+it) \\
&=\int_\lam^{2\lam} \!\!\!\!\! du \, \Big( \frac{q^{(4)}(\! -\! \llcr +u+it)}{\llcr \! -\! u \! -\! it}+\frac{q^{(3)}(\! -\! \llcr+u+it)}{(\llcr \! -\! u\! -\! it)^2} \! -\! \frac{(\gam \! -\! 1)^{\gam+2}(\gam \! -\! 2)}{\gam^{\gam}\llcr^2}(u+it)^{\gam-3}  \Big) . 
\end{align*}
By \eqref{phi_ddd} and \eqref{phi_d4} there are $C_1, C_2 ,\delta\ino (0, \infty)$ such that for all $u\ino (0, 2\delta)$ and all $t\ino (-\delta, \delta)$
\begin{align}
\label{hpvyrx}
&\bigg|  \frac{q^{(4)}(-\llcr+u+it)}{\llcr \! -\! u\! -\! it} \! -\! \frac{(\gam \! -\! 1)^{\gam+2}(\gam \! -\! 2)}{\gam^{\gam}\llcr^2}(u+it)^{\gam-3}\bigg| \le C_1|u+it|^{\gam-2},\\
\label{hpbarxv}
&\textrm{and} \quad \bigg| \frac{q^{(3)}(-\llcr+u+it)}{(\llcr\! -\! u \! -\! it)^2} \bigg| \le C_2|u+it|^{\gam-2}.
\end{align}
Next observe that 
$$ \int_{-\delta}^\delta \!\!\! dt \! \int_{\lam}^{2\lam} \!\!\!\!\!  \frac{\!\!\! du \; \; }{|u+it|^{2-\gam}} \! =\!  \lambda^{\gam } \!\! \int_{-\delta/\lam }^{\delta/\lam} \!\!\! \!\!\! \!\! ds   \int_{1}^{2} \!\!\!  \frac{\!\!\! dv \; \; }{|v+is|^{2-\gam}}
\leq 2\lam^{\gam} \!\! \int_0^{\delta/\lam} \!\!\!\!\!  \!\!\!  \frac{ds}{(1+s^2)^{\frac{2-\gam }{2}}} \leq  2\lam^{\gam} \!\! \int_0^{\delta/\lam} s^{\gam-2} 
=C_3 \lam \; , $$
where $C_3\! =\! 2\delta^{\gam -1}/ (\gam \! -\! 1)$.  
This implies 
\begin{equation}
\label{etccnhu}
\forall \lambda \ino (0, \delta)\, , \quad \lam^{1-\gam} \! \int_{-\delta}^\delta\!  |G(2\lam+i t)\! -\! G(\lam+it)|  dt \! \le \! C_3 (C_1+C_2) \lam^{2-\gam}
\end{equation} 
If $\theta \ino (0, \delta)$, then it implies (\ref{fructimus}). Suppose that $\theta \! \ge \! \delta$. By Lemma \ref{zeuccp}, $q$ is analytic on an open subset $U$ that contains $\overline{\bH}\backslash \{ \! -\! \llcr\}$. Thus, $G$ as defined in (\ref{poiuymlk}) is analytic on $\{ z\ino \bbC: {\rm Re} (z) \ino [0, \delta)\, ; \, \delta \! \le \! |{\rm Im} (z)| \! \le \! \theta \}$ and so is $G^\prime$. Then, we can 
set $C_4 \! :=\! \max \{ |G^\prime (u+it)| \,; \; u\ino [0, \delta), \,  \delta \! \le \! |t| \! \le \! \theta \}$ and by (\ref{etccnhu}), we get 
\begin{eqnarray*}
  \lam^{1-\gam} \! \int_{-\theta}^\theta \!\!  \big|G(2\lam+i t)\! -\! G(\lam+it)\big|  dt & \leq &  \lam^{1-\gam} \!  \int_{-\delta}^\delta \!\!\big|G(2\lam+i t)\! -\! G(\lam+it)\big|  dt + 2(\theta\! -\! \delta)C_4  \lam^{2-\gam} \\
 & \leq  & \big(C_3 (C_1+C_2)+ 2(\theta\! -\! \delta)C_4\big)  \lam^{2-\gam} \; , 
\end{eqnarray*} 
which implies (\ref{fructimus}). \cqfd 

\bigskip

We apply the variant of Ikehara-Ingham Theorem as recalled from Hu \& Shi \cite{HuShi97} in Theorem \ref{IKEA} (see Appendix). Here we take $\mu(dr)\! := \! \un_{(0, \infty)} (r) r^2 g_\Gamma (r) \, dr$, which is a finite measure since $ \int_0^\infty r^2 g_\Gamma (r) \, dr\! = \! q^{(3)} (0) $. More generally observe that for all $\lam\ino (0, \llcr)$, $\int_0^\infty e^{\lam r} r^2 g_\Gamma (r) \, dr \! = \! q^{(3)} (-\lam)\! <\!  \infty$. With the notation of Theorem \ref{IKEA}, $a\! := \! \llcr$ and $F(z)\! = \! q^{(3)} (-z)$ for all $z\ino \bbC$ such that $0\! < \! {\rm Re} (z) \! < \! \llcr$ and $G$ is as in (\ref{poiuymlk}) in Lemma \ref{sazerqsd}, with $b\! :=\! 2\! -\! \gamma$ and 
$c\! := \!  (\gam \! -\! 1)^{\gam +2} / (\llcr^2 \gam^{\gam})$. Thus Theorem \ref{IKEA} implies that 
\begin{equation}
\label{aesvvpt}
A(r) \! := \! \int_r^\infty \! \! \!\!\! u^2g_\Gamma (u) \, du \underset{r\rightarrow \infty}{\sim} K_1r^{1-\gam}e^{-\llcr r} \quad \textrm{where} \quad K_1\! :=\! \frac{(\gam \! -\! 1)^{\gam +2} }{\Gamma_{\! {\rm e}} (2\! -\! \gam) \llcr^2 \gam^{\gam} } \; .
 \end{equation}
 We next set $\phi (u)\! := \! c_\gam \nr  ( \Gam \! \le\!  u^{-\frac{\gam-1}{\gam}} )$, for all $u \ino (0, \infty)$ so that $A(r)\! := \! \int_r^\infty \! 
u^{2-\frac{1}{\gam}} \phi (u) \, du$ by the definition \eqref{def: g_ht} of $g_\Gam$. Note that $\phi$ is decreasing, thus, for all $r, s\ino (0, \infty)$, we get 
$$  \phi (r+s) \int_r^{r+s} \!\!\!\!\!\!\!\! du \, u^{2-\frac{1}{\gam}}\le A(r) \! -\! A(r+s) \le  \phi (r) \int_r^{r+s} \!\!\!\! \!\!\!\! du \, u^{2-\frac{1}{\gam}}\; . $$
To simplify notation we set $\alpha \! := \gam \! -\! 1$ and the previous inequalities implies that 
\begin{align*}
 (r+s)^{\alpha}e^{\llcr(r+s)}\phi(r+s) \int_r^{r+s} \!\!\!\!\!\!\!\! du \, u^{2-\frac{1}{\gam}} & \, \leq  \,   (1+s/r)^{\alpha} e^{\llcr s} r^{\alpha}e^{\llcr r} A(r)- (r+s)^{\alpha}e^{\llcr(r+s)}A(r+s) \\
\textrm{and} \quad  r^{\alpha}e^{\llcr r} \phi(r) \int_r^{r+s} \!\!\!\!\!\!\!\! du \, u^{2-\frac{1}{\gam}}& \geq  r^{\alpha}e^{\llcr r}A(r) \! -\! (1+s/r)^{-\alpha} e^{-\llcr s}(r+s)^{\alpha}e^{\llcr(r+s)}A(r+s) \; .
\end{align*}
As $s$ is fixed and $r\! \to \! \infty$,  $\int_r^{r+s}\! du \, u^{2-\frac{1}{\gam}}\! \sim \! sr^{2-\frac{1}{\gam}}$ and 
the right members of the previous inequalities respectively tend to $K_1 (e^{\llcr s} \! -\! 1)$ and $K_1 (1\! -\! e^{-\llcr s})$ 
by (\ref{aesvvpt}). This implies that for all $s\ino (0, \infty)$, 
$$ K_1 s^{-1} \big( 1\! -\! e^{-\llcr s} \big) \! \leq \!  \liminf_{r\to \infty} r^{\gam +1 -\frac{1}{\gam}} e^{\llcr r} \phi(r) \! \leq \! \limsup_{r\to \infty} r^{\gam +1 -\frac{1}{\gam}} e^{\llcr r} \phi(r) \! \leq \! K_1 s^{-1} \big( e^{\llcr s} \! -\! 1\big) \; .$$
This proves $\lim_{r\to \infty} r^{-\alpha} e^{\llcr r} \phi(r)\! = \! K_1\llcr$ by letting $s$ go to $0+$. Namely, 
$$ c_\gam \nr  \big( \Gam \! \le\!  r^{-\frac{\gam-1}{\gam}} \big) \underset{r \to \infty}{\sim} K_1 \llcr r^{\frac{1}{\gam} - 1 - \gam}e^{-\llcr r } \; , $$
which immediately implies (\ref{dtcxtqcd}) in Theorem \ref{htdm0}.

\paragraph{Proof of (\ref{frqvutgc}) in Theorem \ref{htdm0}.}
The proof of (\ref{frqvutgc}) is quite similar to that of (\ref{dtcxtqcd}). We set 
$$
\forall r\ino (0, \infty) \, , \quad g_D(r)\! := \! c_\gamma r^{-\frac{1}{\gamma}}\nr \big( D \! \le \! 2r^{-\frac{\gam-1}{\gam}} \big) \; .
$$
First note that $g_D(r) \! \le \!  g_\Gam(r)$ for all $r\ino (0, \infty)$, since $D\! \le \! 2\Gam$. Thus, 
\begin{equation}
\label{atfvnxv}
\forall \lam\ino [-\llcr, \infty), \quad \int_0^\infty \! \!\! e^{-\lam r} g_D(r) \, dr \! \le  \!  \int_0^\infty \! \!\! e^{-\lam r} g_\Gamma (r) \, dr \! = \! q^\prime (-\lam)\! < \! \infty \; .
\end{equation}
Next, we deduce from \eqref{echtscal} that
\begin{equation}
\label{lp: gd}
\forall \lam \ino (0, \infty) , \quad \int_0^\infty \!\!\! e^{-\lam r}g_D(r)\, dr=\int_0^\infty c_\gam r^{-1-\frac{1}{\gam}} r e^{-\lam r} \nr \big( r^{\frac{\gam-1}{\gam}}D \! \le  \! 2 \big)
\! = \! \bN \big[ \zeta e^{-\lam \zeta} \indi_{\{D\le 2\}}\big ]\; .
\end{equation}
On the other hand, \eqref{jointlawc} asserts that for all $\lam \ino (0, \infty)$, 
$$
\bN \big[ e^{-\lam \zeta}\indi_{\{D>2\}} \big] \! =\! {\rm L}_\lam(1, 0) \! =\! q(\lam)\! -\! \lam^{\frac{1}{\gam}}\! -\! \big( q(\lam)^\gam \! -\! \lam \big) q'(\lam) \; .
$$
Combining this with the fact that $\bN [1 \! -\! e^{-\lam\zeta}] \! =\! \lam^{1/\gam}$, we get for all $\lam \ino (0, \infty)$,  
$$
\bN \big[ 1\! -\! e^{-\lam\zeta}\indi_{\{D\le 2\}} \big] \! =\! q(\lam)\! -\! \big( q(\lam)^\gam\! -\! \lam \big) q'(\lam).
$$
By differentiating this equality, we deduce from \eqref{lp: gd} that 
$$\forall \lam\ino (0, \infty), \quad  \!\!\!\!  \int_0^\infty \! \!\! e^{-\lam r} g_D(r) \, dr \! = \! Q(\lam) \quad \textrm{where} \quad Q(\lam)\! :=\! 2q'(\lam)  - \gam q(\lam)^{\gam-1} q'(\lam)^2\! -\! \big( q(\lam)^\gam \! -\! \lam \big)q''(\lam).$$
By (\ref{atfvnxv}) the Laplace transform of $g_D$ can be extended analytically on $\bH$ and is continuous on $\overline{\bH}$. By Lemma \ref{zeuccp}, this is also the case of $Q$. The principle of isolated zeroes for analytic function then implies that $Q(z)\! = \! \int_0^\infty \! e^{-zr} g_D(r) \, dr$ for all $z\ino \overline{\bH}$. Moreover, by Lemma \ref{zeuccp}, $Q$ can be extended analytically on $U$ 
and $-\llcr$ is the only singular point of $Q$ in $U$. Thus, for all $z\ino \bH$ we get 
$$  \int_0^\infty \!\!\!\!\!  e^{-zr}  rg_D(r)  dr\! = \! -Q^\prime (z)\! = \! \gam(\gam\! -\! 1)q(z)^{\gam-2}q'(z)^3+3 \big(\gam q(z)^{\gam-1}q'(z)\! -\! 1) q''(z)+\big( q(z)^{\gam}\! -\! z \big) q^{(3)}(z).$$
For all $z\ino \bbC$ such that $0\! < \! {\rm Re} (z) \! < \! \llcr$, we set 
$$ F(z)\!  = \! -Q^\prime (-z)\! = \!  \int_0^\infty \!\!\!\!\!  e^{zr}  rg_D(r)  dr \quad \textrm{and} \quad G(z)\! = \!  \frac{- Q^\prime( - \llcr+z)}{\llcr  - z}  - \frac{2(\gam\! -\! 1)^{\gam+2}}{\gam^\gam\llcr}\, z^{\gam-2}\; .$$
Thanks to (\ref{phi_ddd}) and  (\ref{phi_d4}) in Lemma \ref{zeuccp}, the same arguments as in Lemma \ref{sazerqsd} imply that 
$$ \lam^{1-\gam} \int_{-\theta}^\theta \big|G(2\lam+i t)\! -\! G(\lam+it)\big| \, dt \; \underset{\lam \downarrow 0+}{-\!\!\! -\!\!\! \longrightarrow} \; 0 \; .  $$
We leave the details to the reader (the computations are long but straightforward). 
Then, the variant of Ikehara-Ingham Theorem recalled in Theorem \ref{IKEA} implies that   
\begin{equation}
\label{aqesvvpt}
 \int_r^\infty \! \! \!\!\! u g_D (u) \, du \underset{r\rightarrow \infty}{\sim} K_2r^{1-\gam}e^{-\llcr r} \quad \textrm{where} \quad K_2\! :=\! \frac{2(\gam \! -\! 1)^{\gam +2} }{\Gamma_{\! {\rm e}} (2\! -\! \gam) \llcr \gam^{\gam} } \; .
 \end{equation}
We next argue as in the proof of (\ref{dtcxtqcd}) to derive (\ref{frqvutgc}) from (\ref{aqesvvpt}). \cqfd

\appendix 

\section{Proof of Lemma \ref{mesdec2}.}
\label{pfseclem}
We first recall the following notation from Introduction: let $h\! \in \! \bC(\bbR_+, \bbR_+) $. For any $a\! \in \! [0, h(0)]$, set 
\begin{equation}
\label{defellr}
 \ell_a(h)= \inf \! \big\{ t\! \in \! \bbR_+ : h(t) \! = \! h(0) \! -\! a   \big\} \quad \textrm{and} \quad  r_a(h)=  \inf \! \big\{ t\! \in \! (0, \infty) : h(0)\! -\! a > h(t)  \big\}\wedge \zeta_h \; , 
\end{equation} 
with the convention that $\inf \emptyset \! = \! \infty$. Standard results on stopping times assert that $\ell_a (h)$ and $r_a(h)$ are $[0, \infty]$-valued 
Borel measurable functions of $h$: see for instance Revuz \& Yor \cite{yorbook}, Chapter I, Proposition 4.5 and Proposition 4.6, p.~43. Moreover, it is easy to check that for a fixed $h$, $a \mapsto \ell_a(h)$ is left continuous and that $a\mapsto r_a (h)$ is right continuous. By standard arguments, $ (a, h)\mapsto (\ell_a (h), r_a (h))$ 
is Borel measurable on the set $A\! := \! \{ (a,h) \! \in \! \bbR_+ \! \times \! \bC(\bbR_+, \bbR_+) : a \! \leq \! h(0)\}$. We next recall the following notation: for all $(a, h) \! \in \! A$, we set 
$$ \forall s \! \in \! \bbR_+ , \quad \cE_s (h, a) := h \big(  (\ell_a (h) +s)\! \wedge \! r_a(h) \big)-h(0) +a \; , $$
with the convention that $\cE (h,a)$ is the null function $\mathbf{0}$ if $\ell_a (h)\! = \! \infty$. 
The previous arguments entail that 
\begin{equation}
\label{Borextr}
\textrm{$(a, h) \! \in \!  A  \mapsto \cE(h,a) \! \in \! \bC(\bbR_+, \bbR_+)$ is Borel measurable.} 
\end{equation}

Recall from (\ref{defExc}) the definition of $\mathtt{Exc}$. Recall that $p_H\! : \! [0, \zeta_H ] \! \rightarrow \! \cT_H$ stands for the canonical projection and recall from (\ref{masmea}) that the mass measure $\bbm_H$ is the pushforward measure of the Lebesgue measure on $[0, \zeta_H]$ by $p_H$. Suppose that there exist $r,s \! \in \! (0, \zeta_H)$ such that 
$r \! < \! s$ and such that $H$ is constant on $(r, s)$. Thus $p_H((r, s))\! = \! \{ p_H(r)\}$ and $m_H (\{ p_H (r)\})\! \geq  
\! s\! -\! r \! >\! 0$, which contradicts the fact that $\bbm_H$ is diffuse.   
Recall from (\ref{brleafset}) the definition of the set of leaves $\mathtt{Lf} (\cT_H)$ of $\cT_H$. 
Suppose there exist $r,s \! \in \! (0, \zeta_H)$ such that 
$r \! < \! s$ and such that $H$ is strictly monotone on $(r,s)$.
It easily implies that $p_H((r, s))\! \subset \! \cT_H \backslash \mathtt{Lf}(\cT_H)$, but $\bbm_H (p_H((r,s))) \! \geq \! s\! -\! r \! >\! 0$, which contradicts the fact that $\bbm_H \big(\cT_H \backslash \mathtt{Lf}(\cT_H) \big) \! = \! 0$. Thus, we have proved the following.  
\begin{itemize}  
\item[($\ast$)] \textit{Let $H \! \in \! \mathtt{Exc}$. Let $r,s \! \in \! (0, \zeta_H)$ be such that 
$r \! < \! s$. Then on $(r, s)$, $H$ is not monotone.}  
\end{itemize}

Let $t\! \in \! (0, \infty)$ and $H \! \in \! \mathtt{Exc}$ be such that $\zeta_H \! >\! t $. Recall the following notation 
$$\forall s \! \in \! \bbR_+, \quad  H^-_s= H_{(t-s )_+} , \quad  H^{+}_s= H_{t+ s}, \quad  \overleftarrow{H}^a:= \cE (H^- \! , a) \quad \textrm{and} \quad  \overrightarrow{H}^a:= \cE (H^+\!  , a) , $$
for all $ a \in [0, H_t]$. Note that $H^-_0\! =\! H^+_0\! = \! H_t$. We also recall the following notation
\begin{equation}
\label{respindec0t}
 \cM_{0, t} (H)= \sum_{a \in \cJ_{0, t}} \delta_{(a, \overleftarrow{H}^a, \overrightarrow{H}^a)} \; , 
\end{equation} 
where $\cJ_{0,t}\! := \! \big\{ a\! \in \! [0, H_t] :  \textrm{either $\ell_a (H^-)  \! < \! r_a (H^-) $ or $\ell_a (H^+) \! < \! r_a (H^+) $} \big\} $, which is countable. Then, the definitions (\ref{defellr}) and ($\ast$) entail that 
\begin{equation} 
\label{closhyp}
\forall t\! \in \! (0, \infty), \; \forall H \! \in \! \mathtt{Exc} \; \textrm{such that} \; \zeta_H \! >\! t, \quad \textrm{the closure of $\cJ_{0,t}$ is $[0, H_t]$.}
\end{equation}
We next introduce the compact 
set $C_t\! :=\! \{s \! \in \! [0, \zeta_H \! -\! t]: H_{t+s}\! = \! \inf_{r\in [t, t+s]} H_r \}$, whose Lebesgue measure is denoted by $|C_t|$. 
We easily check that 
$p_H (C_t) \! \subset \! \{ \rho, p_H(t)\} \! \cup \! \big( \cT_H \backslash \mathtt{Lf}(\cT_H)\big)$. 
Since $\bbm_H$ is diffuse and supported by the set of leaves of $\cT_H$, we get $0\! = \! \bbm_H (p_H (C_t))\! \geq \! |C_t|$, which implies that $|C_t|\! = \! 0$. Then note that for all $a \! \in \! [0, H_t]$,       
$$ [0, \ell_a (H^+) ]\backslash C_t \subset \big\{  s \! \in \! [0, \ell_a (H^+)]: H_{t+s}\! > \! \!\!\!\!\! \inf_{\; r\in [t, t+s]} \!\!\!\!\!\! H_r \big\}\; \subset \!\!\!\!\!\!\!\!\!\!\!  \bigcup_{\quad b \in \cJ_{0, t} \cap [0, a)} \!\!\!\!\!\!\!\!\!\!\!  \big( \ell_b (H^+), r_b (H^+) \big) \subset [0, \ell_a (H^+) ]. $$
Since $|C_t|\! = \! 0$, this entails, 
$$\forall a \in [0, H_t], \quad \ell_a (H^+)=\sum_{b \in \cJ_{0,t}} \un_{[0, a)} (b)  \big( r_b (H^+)\! -\!   \ell_b (H^+)\big)= \sum_{b \in \cJ_{0,t}} \un_{[0, a)} (b) \zeta_{\overrightarrow{H}^b}\; .$$
Similar arguments imply that 
\begin{align}
\label{megablurp}
\forall a \in [0, H_t], \quad \ell_a (H^+)&=\sum_{b \in \cJ_{0,t}} \un_{[0, a)} (b) \zeta_{\overrightarrow{H}^b}, \quad   \ell_a (H^-)=\sum_{b \in \cJ_{0,t}} \un_{[0, a)} (b) 
\zeta_{\overleftarrow{H}^b}, \\
 r_a (H^+) &=\sum_{b \in \cJ_{0,t}} \un_{[0, a]} (b) \zeta_{\overrightarrow{H}^b}, \quad  r_a (H^-)=\sum_{b \in \cJ_{0,t}} \un_{[0,a ]} (b) \zeta_{\overleftarrow{H}^b}. \nonumber
\end{align}
Moreover, since $H$ is continuous with compact support, we immediately get 
\begin{equation}
\label{compthyp}
\forall \varepsilon , \eta \!  \in \! (0, \infty) , \quad 
\#   \big\{ a\!  \in\!  \cJ:  \;  \Gamma (\overleftarrow{H}^a) \! \vee \!  \Gamma (\overrightarrow{H}^a) \! >\! \eta \quad \textrm{or} \quad  \zeta_{\overleftarrow{H}^a} \! \!  \vee \! \zeta_{\overrightarrow{H}^a} \! > \! \varepsilon   \big\} < \infty \; . 
\end{equation}
Recall from Remark \ref{etcxfxs} that $\cT_{\overrightarrow{H}^a}$ can be identified with a subtree of $\cT_H$; therefore, 
up to this identification, the set of leaves of $\cT_{\overrightarrow{H}^a}$ is contained in the set of leaves of $\cT_H$ and 
$\bbm_{\overrightarrow{H}^{a}}$ is the restriction of $\bbm_H$ to $\cT_{\overrightarrow{H}^a}$. 
This implies that $\bbm_{\overrightarrow{H}^{a}}$ is diffuse and supported by the set of leaves of $\cT_{\overrightarrow{H}^a}$. Namely, $\overrightarrow{H}^{a}\! \in \! \mathtt{Exc}$. A similar argument show that $\overleftarrow{H}^a \! \in \! \mathtt{Exc}$. This fact combined with (\ref{closhyp}) and (\ref{compthyp}) implies the following: 
\begin{equation}
\label{atterfct}
\forall t\! \in \! (0, \infty), \; \forall H \! \in \! \mathtt{Exc} \; \textrm{such that} \;  \zeta_H \! >\! t, \quad \cM_{0, t} (H) \in \ccM_{{\rm pt}} (E)\; , 
\end{equation}
where $\ccM_{{\rm pt}} (E)$ is as in Definition \ref{defME}. Moreover (\ref{Borextr}) easily 
implies that $(a,t ,H)\!  \mapsto\!  (\overleftarrow{H}^a, \overrightarrow{H}^a)$ is Borel-measurable, which immediately 
implies Lemma \ref{mesdec2} ($i$).

\medskip

Let us prove  Lemma \ref{mesdec2} ($ii$). Recall from Definition \ref{defME} the definition of the sigma field $\cG$ on $\ccM_{{\rm pt}} (E)$. We next fix  $t\! \in \! (0, \infty)$ and $H \! \in \! \mathtt{Exc}$ 
such that $\zeta_H \! >\! t$. First note that (\ref{megablurp}) implies that $\ell_a (H^+)$ and $r_a (H^+)$ are 
$\ccB(\bbR_+) \! \otimes \! \cG$-measurable functions of $(a, \cM_{0, t} (H))$, where $\ccB(\bbR_+)$ stands for the Borel sigma field on $\bbR_+$. We then fix $s\! \in \! \bbR_+$ and we set   
$a(s)= \inf \{a \in \bbR_+: r_a (H^+) \! >\! s   \}$, with the convention that $\inf \emptyset \! = \! \infty$. The previous argument and the fact that $a \! \mapsto \! r_a (H^+)$ is right continuous entail that $a(s)$ can be viewed as a $\cG$-measurable function of $\cM_{0, t} (H)$. Note that if $a(s) \! < \! \infty$, then  
\begin{equation}
\label{rectps}
H_{t + s} = H^+_s = H_t -a(s)+ \overrightarrow{H}^{a(s)} \big( s\! -\! \ell_{a(s)} (H^+) \big) \; .
\end{equation}
Next, for all $a \! \in \! \bbR_+$, set 
$N_a\! =\! \sum_{b \in \cJ_{0, t}} \un_{(a,\infty)} (b) \un_{ \{\zeta_{\overrightarrow{H}^b }>0 \}}$. Recall that we previously proved 
that the closure of the set $\{ b \! \in \! \cJ_{0, t}: \ell_b (H^+) \! <\!  r_b (H^+)\}$ is $[0, H_t]$. Thus  
$H_t \! = \! \sup \{ a \! \in \! \bbR_+: N_a \! > \! 0 \}$, which proves that $H_t$ is a $\cG$-measurable function of $\cM_{0, t} (H)$. Moreover $(a, \cM_{0, t} (H)) \! \mapsto \! \overrightarrow{H}^a$ is $\ccB(\bbR_+) \! \otimes \! \cG$-measurable. Consequently, (\ref{rectps}) implies that $H^+_{s}$ is a $\cG$-measurable function of $\cM_{0, t} (H)$. Since the Borel sigma field on 
$\bC(\bbR_+, \bbR_+)$ is generated by coordinate applications, this implies that $H^+$ is a $\cG$-measurable function of $\cM_{0, t} (H)$. A similar argument shows that $H^-$ is also a $\cG$-measurable function of $\cM_{0, t} (H)$, which easily completes the proof of Lemma \ref{mesdec2} ($ii$). \cqfd

\section{Various results in complex analysis used in the proofs.}
\label{comanarec}    
In this section we briefly recall several results of complex analysis, without proof. 
Let $U$ be a non-empty open subset of $\bbC$ (or of $\bbR$); a function $f\! : U \! \rightarrow \! \bbC$ is called \textit{analytic} if it is locally given by a power serie expansion. 
We refer to the following result as to the \textit{principle of isolated zeroes}. 
\begin{itemize}
\item[] \textit{Let $U$ be a non-empty \textit{connected} open subset of $\bbC$ (or of $\bbR$) and let $f\! : U \! \rightarrow \! \bbC$ be analytic; 
if $f$ is not identically null, then $\{ z\ino U\! : f(z)\! = \! 0\}$ is discrete (namely it has no limit points).} 
\end{itemize}

  We use several times the following statement known as the \textit{Lagrange inversion formula} and whose proof can be found for instance 
in Dieudonn\'e \cite{Die71}, Chapter VIII, (7.3). Let $z_0\ino \bbC$ and $r\ino [0, \infty)$. We denote by $D(z_0, r)\! = \! \{ z\ino \bbC: |z\! -\! z_0| \! < \! r\}$ and by $\overline{D}(z_0, r)\! = \! \{ z\ino \bbC: |z\! -\! z_0| \! \le \! r\}$ respectively the open and the closed disks with centre $z_0$ and radius $r$. 
\begin{prop}
\label{Lagraninv} Let $r\ino (0, \infty)$. 
Let $U$ be a non-empty open subset of $\bbC$ that contains a closed disk $\overline{D}(0, r)$. Let $H \! : U \! \rightarrow \! \bbC$ be analytic. We set $m\! := \! \max_{x\in \overline{D}(0, r) } |H(x)| $. Then, for all 
$z \ino D(0, r/m)$, the equation $x\! =\!  z H(x)$ has a unique solution $x\! =: \! f(z)$ in $D(0, r)$. Moreover $f\! : D(0, r/m) \! \rightarrow \! \bbC$ is analytic and in a neighbourhood of $0$ the following power expansion holds true: 
$$f(z) = \sum_{n\geq 1} \frac{z^n}{n!} \Big(\frac{d^{n-1}}{dx^{n-1}} (H(x))^n\Big)_{|_{x=0}} \; .$$
\end{prop}

 Let $V$ be a non-empty open subset of $\bbC^2$. A function $F\! : V\! \rightarrow \bbC$ is called \textit{analytic in two variables} if for any $(z_0, v_0)\ino V$ there exists $\epp \ino (0, \infty)$ and an array of complex numbers $(a_{m,n})_{m,n\in \bN}$ such that for all $z, v\ino D(0, \epp)$, $(z_0 +z, v_0+v)\ino V$ and $F(z_0+z, v_0+v)\! =\!  \sum_{m,n\in \bN} a_{m,n} z^mv^n$, the sum being absolutely convergent. 
We shall also use a standard result for existence and uniqueness of solution to ordinary differential equation in a complex domain that is recalled as follows (for a proof, see for instance in Hille \cite{Hille} Theorem 2.2.1). 
\begin{prop}
\label{ComplODE} Let $V$ be a non-empty open subset of $\bbC^2$ and let $F\! : V\! \rightarrow \! \bbC$ be analytic in its two variables. Let $(z_0, v_0)\ino V$. Then, there exist $r\ino (0, \infty)$ 
and a unique analytic function $q\! : D(z_0, r) \! \rightarrow \! \bbC$ such that 
$$ \forall z \ino D(z_0, r) \, , \quad (z, q(z))\ino V\, , \quad q^\prime (z)= F(z, q(z)) \quad \textrm{and} \quad q(z_0)\! = \! v_0 \; .$$
\end{prop}

  In the proof of Theorem \ref{htdm0}, we shall use a variant of Ikehara-Ingham Theorem as stated in Hu \& Shi \cite{HuShi97} and whose proof closely follows the main steps of that of Theorem 11, page 234, in Tenenbaum \cite{Tenenbaum}. We recall this result here. To that end, we use the following notations: we set  
$D_-\! := \! \{ z\ino \bbC: \textrm{${\rm Re} (z) \! \le \! 0$ and ${\rm Im} (z)\! = \! 0$}\}$, the negative axis of the complex plane. For any $b\ino \bbC$, we use the following notation $z^b\! := \! \exp (b \log z )$, for all $z\ino \bbC \backslash D_-$, 
where $\log$ is the usual determination of the logarithm in $\bbC \backslash D_-$. Standard results in complex analysis assert that $z\! \mapsto \! z^b$ is analytic in the domain $\bbC\backslash D_-$. 
\begin{thm}
\label{IKEA}
Let $a,b,c\ino (0, \infty)$. Let $\mu$ be a finite measure on $\bbR_+$. Assume that $\int_{\bbR_+} \!\! e^{\lam r}\mu(dr) \! < \! \infty$ for all $\lam \! < \! a$. For all $z\ino \bbC$ such that $0\! < \! {\rm Re} (z) \! < \! a$, we set 

$$  \quad F(z)\! :=\! \int_{\bbR_+} \!\!\!\!\! e^{zr} \mu(dr) \quad \textrm{and} \quad G(z)\! :=\!  \frac{F(a\! -\! z)}{a\! -\! z}\! -\! cz^{-b} \; .$$
We next assume that 
\begin{equation}
\label{eqeta}
\forall \theta \ino (0, \infty) \, , \quad \eta(\lam, \theta):=\lam^{b-1}\! \int_{-\theta}^\theta \big|G(2\lam+i t)\! -\! G(\lam+it)\big| \, dt \; \underset{\lam \downarrow 0+}{-\!\!\! -\!\!\! \longrightarrow} \; 0 \; . 
\end{equation}
Then, there exist two constants $K_1, K_2 \ino (0, \infty)$ such that $K_1$ only depends on $a$, $K_2$ only depends on $a, b, c$ 
and such that for all sufficiently large $r\ino (0, \infty)$
\begin{equation}
\label{asyA}
\Big| e^{ar} r^{1-b} \mu \big( (r, \infty)\big) - \frac{c}{\Gamma_{\! {\rm e}}(b)} \Big| \leq K_2 \! \inf_{\theta\geq K_1 } \!\!  \big( \tfrac{1}{\theta} \! + \! \eta( \tfrac{1}{r} , \theta) \! +\!  (r\theta)^{-b} \big) \; \underset{r \rightarrow \infty}{-\!\!\! -\!\!\! \longrightarrow} \; 0\; .
\end{equation}
\end{thm}
\vspace{-5mm}

{\small

}

\end{document}